\documentclass[sn-mathphys-num]{sn-jnl}% Math and Physical Sciences Numbered Reference Style
\usepackage{graphicx}%
\usepackage{multirow}%
\usepackage{amsmath,amssymb,amsfonts}%
\usepackage{amsthm}%
\usepackage{mathrsfs}%
\usepackage[title]{appendix}%
\usepackage{xcolor}%
\usepackage{textcomp}%
\usepackage{manyfoot}%
\usepackage{booktabs}%
\usepackage{algorithm}%
\usepackage{algorithmicx}%
\usepackage{algpseudocode}%
\usepackage{listings}%
\usepackage{enumerate}
\usepackage{amsmath}
\theoremstyle{thmstyleone}%
\newtheorem{theorem}{Theorem}[section]%  meant for continuous numbers
\newtheorem{proposition}[theorem]{Proposition}%
\newtheorem{lemma}[theorem]{Lemma}
\newtheorem{corollary}[theorem]{Corollary}

\theoremstyle{thmstyletwo}%
\newtheorem{example}[theorem]{Example}%
\newtheorem{remark}[theorem]{Remark}%
\newtheorem{problem}[theorem]{Problem}

\theoremstyle{thmstylethree}%
\newtheorem{definition}[theorem]{Definition}%

\numberwithin{equation}{section}

\newcommand{\normmm}[1]{{\left\vert\kern-0.25ex\left\vert\kern-0.25ex\left\vert #1
		\right\vert\kern-0.25ex\right\vert\kern-0.25ex\right\vert}}
\newcommand{\normm}[1]{{\vert\kern-0.25ex\vert\kern-0.25ex\vert #1
		\vert\kern-0.25ex\vert\kern-0.25ex\vert}}

\allowdisplaybreaks[4]

\begin{document}

\title[The noncompact Schauder fixed point theorem in random normed modules and its applications]{The noncompact Schauder fixed point theorem in random normed modules and its applications}


\author[1]{\fnm{Tiexin} \sur{Guo}}\email{tiexinguo@csu.edu.cn}

\author[2]{\fnm{Yachao} \sur{Wang}}\email{wyachao@foxmail.com}

\author*[3,4]{\fnm{Hong-Kun} \sur{Xu}}\email{xuhk@hdu.edu.cn}

\author [5,6]{George Xianzhi Yuan}\email{george\_yuan99@suda.edu.cn $\&$ george\_yuan99@yahoo.com}

\author[7]{Goong Chen}\email{gchen@math.tamu.edu}
\footnotetext{The first author was supported in part by the National Natural Science Foundation of China (Grant Nos.
12371141,11971483) and the Provincial Natural Science Foundation of Hunan (Grant No. 2023JJ30642). The third author was supported in part by the National Natural Science Foundation of China (Grant No. U1811461) and the Australian Research Council/Discovery Project (Grant No. DP200100124). The fourth author was supported in part by the National Natural Science Foundation of China (Grant Nos. U1811462,71971031).
}


\affil[1]{School of Mathematics and Statistics, Central South University, Changsha 410083, China}

\affil[2]{School of Mathematics and Statistics, Yancheng Teachers University, Yancheng 224002, China}

\affil*[3]{School of Science, Hangzhou Dianzi University, Hangzhou 310018, China}

\affil[4]{College of Mathematics and Information Science,  Henan Normal University,
	Xinxiang, 453007, China}

\affil[5]{School of Business, Sun Yat--sen University, Guangzhou 510275, China}

\affil[6]{Center for Financial Engineering, Soochow University, Suzhou 215006, China}

\affil[7]{Department of Mathematics, Texas A\&M University, College Station, TX 77843, USA}


\abstract{Motivated by the randomized version of the classical Bolzano--Weierstrass theorem, in this paper we first introduce the notion of a random sequentially compact set in a random normed module and develop the related theory systematically. From these developments, we prove the corresponding Schauder fixed point theorem: let $E$ be a random normed module and $G$ a random sequentially compact $L^0$--convex set of $E$, then every $\sigma$--stable continuous mapping from $G$ to $G$ has a fixed point, which unifies all the previous random generalizations of the Schauder fixed point theorem. As one of the applications of the theorem, we prove the existence of Nash equilibrium points in the context of conditional information. It should be pointed out that the main challenge in this paper lies in overcoming noncompactness since a random sequentially compact set is generally noncompact.}




\maketitle

\setcounter{secnumdepth}{0}
\section{Introduction and main results}\label{intro}
Throughout, $\mathbb{R}^d$ denotes the $d$--dimensional Euclidean space. The classical Brouwer fixed point theorem \cite{Brou} states that each compact convex set $G$ of $\mathbb{R}^d$ has the fixed point property (briefly, FPP), namely, every continuous mapping from $G$ to $G$ has a fixed point. It is well known that the Brouwer fixed point theorem and its extensions have played a central role in analysis, optimization and economic theory \cite{GD,Bor,Nash}, and numerous others.

\par
In 1930, Schauder \cite{Sch} first established an infinite--dimensional generalization of the Brouwer fixed point theorem, now well known as the Schauder fixed point theorem. It states that each compact convex set of a normed space has the FPP. As noted in \cite{ET21,GD,Mau}, the Schauder fixed point theorem is one of the most powerful tools in dealing with nonlinear problems in analysis and, in particular, it has played a major role in the development of fixed point theory and the theory of differential equations. Due to the fact that the Schauder fixed point theorem has fundamental importance, the theorem has been generalized in various directions by different methods. These generalizations can be divided into two kinds: one is purely topological, see, for instance \cite{AMO,Cau,Cau2,Cau3,Cau4,CT,ET21,Fan,KRAS,Park08,Re73,Ty,Zei} and the references therein; the other is of probabilistic interest in connection with stochastic analysis and stochastic finance, see for instance \cite{Bhar,DKKS} and the references therein.

\par
This paper is organized as follows. Section 1 provides some prerequirites and survey. Section \ref{section2} is devoted to a systematic study of random sequential compactness and random total boundedness, in particular the random Hausdorff theorem (Theorem \ref{theorem2.3}) is established. In Section \ref{section3}, we first establish the equivalence between random sequential continuity and $\mathcal{T}_c$--continuity for a $\sigma$--stable mapping, see Lemma \ref{lemma3.3}, then define random Schauder projection and further finish the proofs of Theorem \ref{theorem1.12} and \ref{theorem1.15}. Section \ref{section4} is devoted to the proof of Theorem \ref{theorem1.19}. Finally, Section \ref{section5} concludes this paper with some important remarks.

\setcounter{secnumdepth}{3}
\section{The Schauder conjecture. Prerequirites.}
\par
In this section, we provide a survey of the recent advances toward solving the Schauder conjecture. Appropriate prerequirites are also introduced.

\subsection{The FPP and compactness of a convex set in a metrizable linear topological space}\label{subsection1.1}
\par
Following \cite{Sch}, Schauder conjectured that every compact convex set in an arbitrary Hausdorff linear topological space has the FPP, which is well known as Problem 54 in \cite{Mau}. In 1935, Tychonoff \cite{Ty} proved that the Schauder conjecture holds in a Hausdorff locally convex space. In 1964, Fan \cite{Fan} proved that the conjecture holds in a Hausdorff linear topological space admitting sufficiently many continuous linear functionals. In 2001, Cauty \cite{Cau} claimed that he completely solved the conjecture. However, it was later found out that Cauty's proof given in \cite{Cau} as well as its elaboration in \cite{Do} contained a gap, see \cite[p.131]{Mau} and \cite[p.91]{Isac}. Thereafter, Cauty \cite{Cau2,Cau3,Cau4} filled in the gap so that he can affirmatively resolve the Schauder conjecture. In 2021, Ennassik and Taoudi \cite{ET21} provided a nice alternative proof to the Schauder conjecture.

\par
On the other hand, in the course of solving the Schauder conjecture the problem of whether the FPP of a convex set in a Hausdorff linear topological space also implies its compactness was first studied by Klee \cite{Kl}, where an affirmative answer was given for the case of a metrizable locally convex space, and eventually solved by Dobrowoski and Marciszewski \cite{DM}, where it was proved that the FPP of a convex set in a metrizable linear topological space implies its compactness and it was also illustrated that the FPP of a convex set in a Hausdorff locally convex space does not always imply its compactness. Therefore,  the results in \cite{Cau2} and \cite{DM} can combined to yield that a convex set in a metrizable linear topological space has the FPP iff the convex set is compact.

\subsection{Two topologies for a random normed module and the notion of a $\sigma$--stable set}\label{subsection1.2}
\par
The development of the field of random functional analysis is based on the idea of randomizing space theory in traditional functional analysis. Its basic framework consists of random normed modules and random locally convex modules. Thus, random functional analysis may also be aptly called functional analysis over random normed modules and random locally convex modules, see \cite{Guo3,GZWG,GZWW,GZWYYZ} for some important advances in the field. Although general random locally convex modules have been deeply studied in \cite{Guo3,GZWYYZ,GZZ1,Wu}, the work in this paper involves essentially only random normed modules. To save space, in what follows we only introduce some needed concepts related to random normed modules.

\par
The notion of a random normed module is a generalization by randomizing that of an ordinary normed space. It is well known that the norm topology for a normed space is unique, but the $L^0$--norm of a random normed module can induce two topologies---the $(\varepsilon,\lambda)$--topology and the locally $L^0$--convex topology, which are both frequently employed in the theoretical development of random normed modules and their applications in finance. Since the Schauder fixed point theorem is of topological nature, generalizing the theorem in random normed modules necessarily involves their complicated topological structures.

\par
To introduce the two types of topologies for random normed modules, we first recall some frequently used terminologies and notation as follows.

\par
Throughout this paper, $(\Omega, \mathcal{F},P)$ denotes a given probability space, $\mathbb{N}$ the set of positive integers, $\mathbb{K}$ either the scalar field $\mathbb{R}$ of real numbers or $\mathbb{C}$ of complex numbers, $L^0(\mathcal{F}, \mathbb{K})$ the algebra of equivalence classes of $\mathbb{K}$--valued $\mathcal{F}$--measurable random variables on $(\Omega, \mathcal{F},P)$ under the usual algebraic operations, $L^0(\mathcal{F}): = L^0(\mathcal{F}, \mathbb{R})$ and $\bar{L}^0(\mathcal{F})$ the set of equivalence classes of extended real--valued $\mathcal{F}$--measurable random variables on $(\Omega, \mathcal{F},P)$.

\par
Proposition \ref{proposition1.1} below can be aptly called the randomized version of the order completeness of $\mathbb{R}$, which is frequently used in this paper as well as in the development of random functional analysis.

\begin{proposition}[\cite{DS58}]\label{proposition1.1}
	$\bar{L}^0(\mathcal{F})$ is a complete lattice under the partial order: $\xi\leq \eta$ iff $\xi^0(\omega) \leq \eta^0(\omega)$ for almost all $\omega$ in $\Omega$ (briefly, $\xi^0(\omega) \leq \eta^0(\omega)$ a.s. or $\xi^0 \leq \eta^0$ a.s.), where $\xi^0$ and $\eta^0$ are arbitrarily chosen representatives of $\xi$ and $\eta$ in $\bar{L}^0(\mathcal{F})$, respectively. For a nonempty subset $H$ of $\bar{L}^0(\mathcal{F})$, $\bigvee H$ and $\bigwedge H$ stands for the supremum and infimum of $H$, respectively. Furthermore, the lattice $\bar{L}^0(\mathcal{F})$ possesses the following nice properties:
	\begin{enumerate}[(1)]
		\item There exist two sequences $\{a_n, n\in \mathbb{N}\}$ and $\{b_n, n\in \mathbb{N}\}$ in $H$ such that $\bigvee_{n\geq 1} a_n= \bigvee H$ and $\bigwedge_{n\geq 1} b_n= \bigwedge H$.
		\item If $H$ is directed upwards (resp., downwards), namely there exists some $h_3\in H$ for any two elements $h_1$ and $h_2$ in $H$ such that $h_1\bigvee h_2\leq h_3~(resp., h_1\bigwedge h_2\geq h_3)$, then $\{a_n, n\in \mathbb{N} \}~(resp., \{b_n, n\in \mathbb{N}\})$ in (1) can be chosen as nondecreasing (resp., nonincreasing).
		\item $(L^0(\mathcal{F}), \leq)$ is a Dedekind complete lattice.
	\end{enumerate}
\end{proposition}

\par
As usual, $\xi< \eta$ means $\xi\leq \eta$ and $\xi\neq \eta$ for any $\xi$ and $\eta$ in $\bar{L}^0(\mathcal{F})$, whereas, for any $A\in \mathcal{F}$, $\xi< \eta$ on $A$ $(\xi \leq \eta$ on $A)$ means $\xi^0(\omega)< \eta^0(\omega)~(\xi^0(\omega) \leq \eta^0(\omega))$ for almost all $\omega$ in $A$, where $\xi^0$ and $\eta^0$ are arbitrarily chosen representatives of $\xi$ and $\eta$, respectively.

\par
The following is a set of special notation used in this paper:
\par
$L^0_+(\mathcal{F})= \{\xi\in L^0(\mathcal{F}): \xi\geq 0 \};$
\par
$L^0_{++}(\mathcal{F})= \{\xi\in L^0(\mathcal{F}): \xi> 0~ \text{on} ~\Omega \};$
\par
$\bar{L}^0_+(\mathcal{F})= \{\xi \in \bar{L}^0(\mathcal{F}): \xi\geq 0 \};$
\par
$\bar{L}^0_{++} (\mathcal{F})=\{\xi \in \bar{L}^0(\mathcal{F}): \xi> 0~ \text{on} ~\Omega \}.$
\par
$\tilde{I}_A$ always stands for the equivalence class of $I_A$ for any $A\in \mathcal{F}$, where $I_A$ is the characteristic function of $A$, namely $I_A(\omega)=1$ if $\omega\in A$ and $I_A(\omega)=0$ otherwise.

\par
The theory of random normed modules was independently introduced and developed by Guo \cite{Guo1,Guo2} in connection with the theory of random normed spaces \cite{SS}, and by Gigli \cite{Gigli1,Gigli} in connection with nonsmooth differential geometry on metric measure spaces, where the notion of an $L^{0}$-normed $L^{0}$-module amounts to that of a random normed module. The development of nonsmooth differential geometry on metric measure spaces \cite{Gigli1,Gigli,LP,LPV,GMT} has showed that the theory of random normed modules under the $(\varepsilon,\lambda)$-topology is playing an essential role in the field of nonsmooth differential geometry. It should be pointed out that the theory of random normed modules under the locally $L^{0}$-convex topology was independently developed by Filipovi\'{c}, et al \cite{FKV1} in connection with conditional risk measures.

\par
Definition \ref{definition1.2} here is adopted from \cite{Guo1,Guo2} by following the traditional nomenclature of random metric spaces and random normed spaces (see \cite[Chapters 9 and 15]{SS}).

\begin{definition}[\cite{Guo1,Guo2}]\label{definition1.2}
	An ordered pair $(E, \|\cdot\|)$ is called a random normed module (briefly, an $RN$ module) over $\mathbb{K}$ with base $(\Omega, \mathcal{F}, P)$ if $E$ is a left module over the algebra $L^0(\mathcal{F}, \mathbb{K})$ (briefly, an $L^0(\mathcal{F}, \mathbb{K})$--module) and $\|\cdot\|$ is a mapping from $E$ to $L^0_+(\mathcal{F})$ such that the following conditions are satisfied:
	\begin{enumerate}[(1)]
		\item $\|\xi x\|= |\xi| \|x\|$ for any $\xi \in L^0(\mathcal{F}, \mathbb{K})$ and any $x\in E$, where $\xi x $ stands for the module multiplication of $x$ by $\xi$ and $|\xi|$ is the equivalence class of $|\xi^0|$ for an arbitrarily chosen representative $\xi^0$ of $\xi$ (of course, $|\xi^0|$ is the function defined by $|\xi^0|(\omega)= |\xi^0(\omega)|$ for any $\omega\in \Omega$);
		\item $\|x+y\|\leq \|x\|+ \|y\|$ for any $x$ and $y$ in $E$;
		\item $\|x\|=0$ implies $x=\theta$ (the null element of $E$).
	\end{enumerate}
	As usual, $\|\cdot\|$ is called the $L^0$--norm on $E$.
\end{definition}
\par
A mapping $\|\cdot\|$ from an $L^0(\mathcal{F}, \mathbb{K})$--module $E$ to $L^0_+(\mathcal{F})$ is called an $L^0$--seminorm if $\|\cdot\|$ only satisfies conditions (1) and (2) in Definition \ref{definition1.2}. Similarly, one can have the notion of a random inner product module (briefly, an $RIP$ module), see \cite{Guo3} for details.

\par
When $(\Omega, \mathcal{F}, P)$ is trivial, namely $\mathcal{F}= \{\Omega, \emptyset\}$, it is clear that an $RN$ module with base $(\Omega, \mathcal{F}, P)$ just reduces to an ordinary normed space and, thus, the former is a generalization by randomness of the latter.

\par
The algebra $L^0(\mathcal{F}, \mathbb{K})$ is the simplest $RN$ module with the $L^0$--norm $\|\cdot\|$ defined by $\|\xi\|= |\xi|$ for any $\xi \in L^0(\mathcal{F}, \mathbb{K})$ (thereafter, the $L^0$--norm $\|\cdot\|$ is always denoted by $|\cdot|$). It is well known that $L^0(\mathcal{F}, \mathbb{K})$, endowed with the topology of convergence in probability, is a metrizable linear topological space (in fact, also a metrizable topological algebra). This topology has a local base $\{N_{\theta}(\varepsilon,\lambda): \varepsilon>0, 0<\lambda <1\}$, where $N_{\theta}(\varepsilon,\lambda)= \{x\in L^0(\mathcal{F}, \mathbb{K}): P\{\omega\in \Omega: |x|(\omega)< \varepsilon \}>1-\lambda \}$ for any positive numbers $\varepsilon$ and $\lambda$ with $0<\lambda<1$, this observation leads Schweizer and Sklar to introduce the $(\varepsilon,\lambda)$--topology for an abstract probabilistic metric or normed space in 1960, see \cite{SS} for details. Similarly, we can usually introduce the $(\varepsilon,\lambda)$--topology for an $RN$ module as we have made in \cite{Guo3,GP}: let $(E, \|\cdot\|)$ be an $RN$ module over $\mathbb{K}$ with base $(\Omega, \mathcal{F}, P)$, let $N_{\theta}(\varepsilon,\lambda)= \{x\in E: P\{\omega\in \Omega~|~ \|x\|(\omega)< \varepsilon \}>1-\lambda \}$. Then it is easy to see that the family $\{N_{\theta}(\varepsilon,\lambda): \varepsilon>0, 0<\lambda<1 \}$ forms a local base of some metrizable linear topology for $E$, called the $(\varepsilon,\lambda)$--topology induced by $\|\cdot\|$. From now on, we always denote by $\mathcal{T}_{\varepsilon,\lambda}$ the $(\varepsilon,\lambda)$--topology for an arbitrary $RN$ module $(E,\|\cdot\| )$ if no confusion arises. With regard to the $(\varepsilon,\lambda)$--topology, we have the following.

\begin{proposition}[\cite{Guo3,GP}]\label{proposition1.3}
	Let $(E,\|\cdot\| )$ be an $RN$ module over $\mathbb{K}$ with base $(\Omega, \mathcal{F}, P)$. Then $(E, \mathcal{T}_{\varepsilon,\lambda})$ is a metrizable topological module over the topological algebra $(L^0(\mathcal{F},$ $ \mathbb{K}), \mathcal{T}_{\varepsilon,\lambda})$, namely $\mathcal{T}_{\varepsilon,\lambda}$ is metrizable and makes the module multiplication $\cdot: (L^0(\mathcal{F},$ $ \mathbb{K}), \mathcal{T}_{\varepsilon,\lambda}) \times (E, \mathcal{T}_{\varepsilon,\lambda}) \rightarrow (E, \mathcal{T}_{\varepsilon,\lambda})$ jointly continuous.
\end{proposition}

\par
It is well known that the conjugate space (or topological dual) of $(L^0(\mathcal{F},$ $ \mathbb{K}), \mathcal{T}_{\varepsilon,\lambda})$ is trivial when $\mathcal{F}$ is atomless and, thus, please bear in mind that the $(\varepsilon,\lambda)$--topology for an $RN$ module is often nonlocally convex. For an $RN$ module over $\mathbb{K}$ with base $(\Omega, \mathcal{F}, P)$, let $E^*_{\varepsilon,\lambda}$ be the $L^0(\mathcal{F}, \mathbb{K})$--module of continuous module homomorphisms from $(E, \mathcal{T}_{\varepsilon,\lambda}) $ to $(L^0(\mathcal{F},$ $ \mathbb{K}), \mathcal{T}_{\varepsilon,\lambda})$, called the random conjugate space of $(E, \mathcal{T}_{\varepsilon,\lambda}) $ under $\mathcal{T}_{\varepsilon,\lambda}$. It is well known that a linear operator $f$ from $E$ to $L^0(\mathcal{F}, \mathbb{K})$ belongs to $E^*_{\varepsilon,\lambda}$ iff $f$ is a.s. bounded, namely there exists $\xi\in L^0_+(\mathcal{F})$ such that $|f(x)|\leq \xi \|x\|$ for any $x\in E$. The theory of random conjugate spaces has played an essential role in the development of $RN$ modules, see \cite{Guo3,GL} for details.

\par
Before the year 2009, $RN$ modules were all developed under the $(\varepsilon,\lambda)$--topology in \cite{Guo3}. In 2009, Filipovi\'{c}, et al \cite{FKV1} developed the theory of a locally $L^0$--convex module in order to establish a theory of more general convex analysis (called $L^0$--convex analysis or random convex analysis) suitable for the study of conditional convex risk measures. Besides applications in financial mathematics, it was shown in \cite{Guo3,WG,Zapa} that the work on topological modules in \cite{FKV1} also leads to another topology---the locally $L^0$--convex topology for an $RN$ module.

\par
We first recall that a nonempty subset $G$ of an $L^0(\mathcal{F}, \mathbb{K})$--module $E$ is said to be $L^0$--convex if $\xi x+ (1-\xi) y\in G$ for any $x$ and $y$ in $G$ and $\xi\in L^0_+(\mathcal{F})$ with $0\leq \xi \leq 1$, $L^0$--absorbent if for any $x\in E$ there exists $\delta\in L^0_{++} (\mathcal{F})$ such that $\lambda x\in G$ for any $\lambda\in L^0(\mathcal{F}, \mathbb{K})$ with $|\lambda|\leq \delta$, and $L^0$--balanced if $\lambda x\in G$ for any $x\in G$ and $\lambda\in L^0(\mathcal{F}, \mathbb{K})$ with $|\lambda|\leq 1$. For $L^0(\mathcal{F}, \mathbb{K})$, Filipovi\'{c}, et al \cite{FKV1} first considered the ring topology $\mathcal{T}_c$ rather than its $(\varepsilon,\lambda)$--topology: a subset $G$ of $L^0(\mathcal{F}, \mathbb{K})$ is said to be $L^0$--open if for any $x\in G$ there exists some $\varepsilon\in L^0_{++} (\mathcal{F})$ such that $x+ V_{\theta} (\varepsilon)\subset G$, where $V_{\theta}(\varepsilon)= \{x\in L^0(\mathcal{F}, \mathbb{K}): |x|< \varepsilon ~\text{on}~\Omega \}$. The family $\mathcal{T}_c$ of $L^0$--open sets of $L^0(\mathcal{F}, \mathbb{K})$ then forms a Hausdorff topology such that $(L^0(\mathcal{F}, \mathbb{K}), \mathcal{T}_c)$ is a topological ring, namely the addition and ring multiplication operations are both jointly continuous. (A caution here is that since $\mathcal{T}_c$ is too strong, the scalar multiplication sending $\alpha\in \mathbb{K}$ to $\alpha\cdot 1$ is not necessarily continuous.) Further, in \cite{FKV1} a topological space $(E,\mathcal{T})$ is called a topological $L^0$--module if $E$ is an $L^0(\mathcal{F}, \mathbb{K})$-module and if $(E,\mathcal{T})$ is a topological module over the topological ring $(L^0(\mathcal{F}, \mathbb{K}), \mathcal{T}_c)$. A topological $L^0$--module $(E,\mathcal{T})$ is said to be a locally $L^0$--convex module if it has a local base consisting of $L^0$--convex, $L^0$--balanced and $L^0$--absorbent sets, and $\mathcal{T}$ is called the locally $L^0$--convex topology for $E$. On $RN$ modules, Filipovi\'{c}, et al \cite{FKV1} gave the following.

\begin{proposition}[\cite{FKV1}]\label{proposition1.4}
	Let $(E, \|\cdot\|)$ be an $RN$ module over $\mathbb{K}$ with base $(\Omega, \mathcal{F}, P)$, then the family $\{V_{\theta}(\varepsilon): \varepsilon\in L^0_{++}(\mathcal{F}) \}$ forms a local base of some Hausdorff locally $L^0$--convex topology for $E$, called the locally $L^0$--convex topology induced by $\|\cdot\|$ and denoted by $\mathcal{T}_c$, where $V_{\theta}= \{x\in E: \|x\|< \varepsilon~\text{on}~\Omega\}$.
\end{proposition}

\par
From now on, for an arbitrary $RN$ module $(E, \|\cdot\|)$ we always denote by $\mathcal{T}_c$ the locally $L^0$--convex topology induced by $\|\cdot\|$ if no ambiguites occur. Further, if $(E, \|\cdot\|)$ is an $RN$ module over $\mathbb{K}$ with base $(\Omega, \mathcal{F}, P)$, the $L^0(\mathcal{F}, \mathbb{K})$-module $E^*_c$ of continuous module homomorphisms from $(E, \mathcal{T}_c)$ to $(L^0(\mathcal{F}, \mathbb{K}), \mathcal{T}_c)$ is called the random conjugate space of $(E, \|\cdot\|)$ under $\mathcal{T}_c$. It is also well known that $E^*_c=E^*_{\varepsilon,\lambda}$, see \cite{Guo3} for details. Throughout this paper, we always use $E^*$ for $E^*_{\varepsilon,\lambda}$ or $E^*_c$.

\par
For an $RN$ module $(E, \|\cdot\|)$, it is obvious that $\mathcal{T}_c$ is much stronger than $\mathcal{T}_{\varepsilon,\lambda}$. It has been shown that the two types of topologies have their respective advantages and disadvantages whether for the theoretical investigations or for financial mathematical applications, see \cite{Guo5,GZWW,GZWYYZ,GZZ1,GZZ3} for details. In fact, it often requires one to combine the advantages of both topologies to solve the problems in the theory of $RN$ modules and their applications. In the process it is important to establish the connection between the two types of topologies. The following notion of $\sigma$--stability plays a key role in establishing this connection.

\begin{definition}[\cite{Guo3}]\label{definition1.5}
	Let $E$ be an $L^0(\mathcal{F}, \mathbb{K})$-module and $G$ a nonempty subset of $E$. $G$ is said to be stable if $\tilde{I}_A x+\tilde{I}_{A^c} y\in G$ for any $x$ and $y$ in $G$ and any $A\in \mathcal{F}$, where $A^c=\Omega \backslash A$. $G$ is said to be $\sigma$--stable (or, to have the countable concatenation property in the original terminology of \cite{Guo3}) if for any sequence $\{x_n, n\in \mathbb{N} \}$ in $G$ and any countable partition $\{A_n, n\in \mathbb{N} \}$ of $\Omega$ to $\mathcal{F}$ (namely each $A_n\in \mathcal{F}, A_i\cap A_j=\emptyset$ for any $i\neq j$, and $\bigcup_{n\in \mathbb{N}} A_n= \Omega$) there exists $x\in G$ such that $\tilde{I}_{A_n} x= \tilde{I}_{A_n} x_n$ for each $n\in \mathbb{N}$.
\end{definition}

\par
It is known from \cite{Guo3,WGL,WZZ} that when $(E, \|\cdot\|)$ is an $RN$ module or a more general regular $L^0$--module, $x$ in Definition \ref{definition1.5} is unique and that $x$ can be written as $\sum_{n\in \mathbb{N}} \tilde{I}_{A_n} x_n$ or $\sum_n \tilde{I}_{A_n} x_n$.

\par
Throughout this paper, let $(E, \|\cdot\|)$ be an $RN$ module over $\mathbb{K}$ with base $(\Omega, \mathcal{F}, P)$, $G$ a $\sigma$--stable subset of $E$ and $H$ a nonempty subset of $G$, $\sigma(H)$ always stands for the $\sigma$--stable hull of $H$, namely $\sigma(H)= \{\sum_n \tilde{I}_{A_n} x_n: \{A_n, n\in \mathbb{N}\}$ is a countable partition of $\Omega$ to $\mathcal{F}$ and $\{x_n, n\in \mathbb{N}\}$ is a sequence in $H \}$. That is to say, $\sigma(H)$ is the smallest $\sigma$--stable subset containing $H$.

\par
It is clear that an $L^0(\mathcal{F}, \mathbb{K})$--module $E$ is itself stable and an $L^0$--convex set of $E$ is also stable. It is already pointed out in \cite{Guo3} that $E^*$ is $\sigma$-stable for an $RN$ module $(E, \|\cdot\|)$.

\par
Now, the connection between some basic results derived from both types of topologies can be summarized in Proposition \ref{proposition1.6} for easier future reference.

\begin{proposition}[\cite{Guo3,GZZ1}]\label{proposition1.6}
	Let $(E, \|\cdot\|)$ be an $RN$ module, then the following statements hold:
	\begin{enumerate}[$(1)$]
		\item $\overline{G}_{\varepsilon,\lambda}= \overline{G}_c$ for a $\sigma$--stable set $G$ of $E$, where $\overline{G}_{\varepsilon,\lambda}$ and $\overline{G}_c$ denote the closures of $G$ under $\mathcal{T}_{\varepsilon,\lambda}$ and $\mathcal{T}_c$, respectively.
		\item $E$ is $\mathcal{T}_{\varepsilon,\lambda}$--complete iff both $E$ is $\mathcal{T}_c$--complete and $E$ is $\sigma$--stable.
	\end{enumerate}
\end{proposition}

\begin{remark}\label{remark1.7}
	Statement (2) in Proposition \ref{proposition1.6} is merely a restatement of Theorem 3.18 of \cite{Guo3} since $E$ is stable and it is easy to prove that $E$ must be $\sigma$--stable when $E$ is $\mathcal{T}_{\varepsilon,\lambda}$--complete. But for a nonempty subset $G$ of an $RN$ module, by the same method of proof for Theorem 3.18 of \cite{Guo3} one can easily see that $G$ is always $\mathcal{T}_c$--complete whenever $G$ is $\mathcal{T}_{\varepsilon,\lambda}$--complete. It is also easy to see that $G$ must be $\mathcal{T}_{\varepsilon,\lambda}$--complete whenever $G$ is both $\mathcal{T}_c$--complete and $\sigma$--stable.
\end{remark}

\par
In what follows, we give two important examples of $RN$ modules used in this paper.

\par
For the first example, let us recall that a mapping $V$ from $(\Omega, \mathcal{F}, P)$ to a metric space $(M,d)$ is called a random element if $V^{-1}(G):= \{\omega\in \Omega: V(\omega)\in G \}\in \mathcal{F}$ for each open subset $G$ of $M$; a random element $V$ is said to be simple if the range of $V$ is a finite subset of $M$; a random element $V$ is said to be strong if $V$ is a pointwise limit of a sequence of simple random elements.

\par
As usual, $L^0(\mathcal{F},M)$ denotes the set of equivalence classes of strong random elements from $(\Omega, \mathcal{F}, P)$ to a metric space $(M,d)$.

\begin{example}[\cite{Guo6}]\label{example1.8}
	When the metric space $(M,d)$ is a normed space $(X,\|\cdot\|)$ over the scalar field $\mathbb{K}$, $L^0(\mathcal{F},X)$ naturally forms an $RN$ module over $\mathbb{K}$ with base $(\Omega, \mathcal{F}, P)$ in the following manner. First, $L^0(\mathcal{F},X)$ is a linear space over $\mathbb{K}$ under the ordinary addition and scalar multiplication operations on equivalence classes. Further, the module multiplication is defined by $\xi p: =$ the equivalence class of $\xi^0 \cdot p^0$, where $\xi^0$ and $p^0$ stand for arbitrarily chosen representatives of $\xi\in L^0(\mathcal{F}, \mathbb{K})$ and $p\in L^0(\mathcal{F},X)$, respectively, while $(\xi^0\cdot p^0)(\omega)= \xi^0(\omega) \cdot p^0(\omega)$ for any $\omega\in \Omega$. Finally, for any $p\in L^0(\mathcal{F},X)$, the $L^0$--norm of $p$ is defined as the equivalence of $\|p^0\|$, where $\|p^0\|(\omega):= \|p^0(\omega)\|$ for any $\omega\in \Omega$ and $p^0$ is an arbitrarily chosen representative of $p$. As usual, we still use $\|\cdot\|$ for the $L^0$--norm on $L^0(\mathcal{F},X)$.
\end{example}

\par
It is easy to check that, for any nonempty subset $V$ of a normed space $(X, \|\cdot\|)$, $L^0(\mathcal{F}, V)$ is a $\sigma$--stable subset of $L^0(\mathcal{F},X)$. Especially, when $V$ is a convex subset of $X$, $L^0(\mathcal{F},V)$ is also an $L^0$--convex subset of $L^0(\mathcal{F},X)$.

\par
Example \ref{example1.9} below has played a crucial role in the theory of conditional convex risk measures (equivalently, conditional concave utility functions) since it has been employed as standard model spaces for the latter, see \cite{HR,FKV1,FKV2,GZZ3} for details.

\begin{example}[\cite{HR,FKV1,FKV2}]\label{example1.9}
	Let $(\Omega, \mathcal{E}, P)$ be a probability space and $\mathcal{F}$ a $\sigma$--subalgebra of $\mathcal{E}$. For any extended positive number $r\in [1,+\infty]$, define the mapping $\normm{\cdot}_r: L^0(\mathcal{E})\rightarrow \bar{L}^0_+(\mathcal{F})$ as follows: for any $x\in L^0(\mathcal{E})$,
	\begin{eqnarray}\nonumber
		\normm{x}_r=\left\{
		\begin{aligned}
			\lim_{n\rightarrow \infty} E[|x|^r \wedge n~|\mathcal{F}~]^{\frac{1}{r}} &, & when~ r\in [1,+\infty); \\
			\wedge\{\eta\in \bar{L}^0_+(\mathcal{F}): |x|\leq \eta \} &, & when~ r=+\infty.
		\end{aligned}
		\right.
	\end{eqnarray}
	Again, let $L^r_{\mathcal{F}} (\mathcal{E})=\{x\in L^0(\mathcal{E}): \normm{x}_r\in L^0_+(\mathcal{F}) \}$, then it is known that $L^r_{\mathcal{F}} (\mathcal{E})$ = $L^0(\mathcal{F})\cdot L^r(\mathcal{E}):=$ $\{\xi\cdot x: \xi\in L^0(\mathcal{F})~\text{and}~ x\in L^r(\mathcal{E})\}$, where $L^r(\mathcal{E})$ denotes the usual Banach space of equivalence classes of $r$--integrable or essentially bounded real--valued $\mathcal{E}$--measurable functions on $\Omega$. Namely, $L^r_{\mathcal{F}}(\mathcal{E})$ is the $L^0(\mathcal{F})$--module generated by $L^r(\mathcal{E})$. It is also easy to see that $(L^r_{\mathcal{F}}(\mathcal{E}), \normm{\cdot}_r)$ is an $RN$ module over $\mathbb{R}$ with base $(\Omega, \mathcal{F}, P)$.
\end{example}

\par
It is also known from \cite{Guo6,GZZ3} that $L^r_{\mathcal{F}} (\mathcal{E})$ is $\sigma$--stable $RN$ module.

\subsection{The Schauder fixed point theorem for a $\sigma$-stable continuous \\
selfmapping on a random sequentially compact $L^0$--convex set of an $RN$ module}\label{subsection1.3}
\par
As convex sets and convex functions play useful roles in functional analysis and its applications to single--period mathematical finance, $L^0$--convex sets and $L^0$--convex functions do likewise in random functional analysis and its applications to multi--period mathematical finance. See, for example, \cite{FKV1,FKV2,Guo3,GZWYYZ,GZZ3,HR} for details. However, closed $L^0$--convex sets frequently encountered in $RN$ modules and their financial mathematical applications are not compact in general under the $(\varepsilon,\lambda)$--topology. (Of course, they are not compact under the stronger locally $L^0$--convex topology.) In fact, it is known from \cite{Guo4} that when $(M,d)$ is a metric space containing at least two points, $L^0(\mathcal{F}, M)$ is compact under the topology of convergence in probability measure iff both $(M, d)$ is compact and $\mathcal{F}$ is essentially generated by $P$--atoms, see \cite[Theorem 2.2]{Guo4} for details. Hence for a compact convex set $V$ of a normed space $(B, \|\cdot\|)$, $L^0(\mathcal{F}, V)$ is not compact in general under the $(\varepsilon,\lambda)$--topology as a closed $L^0$--convex subset of the $RN$ module $L^0(\mathcal{F}, B)$. Consequently, although a Schauder conjecture has been completely solved by Cauty \cite{Cau,Cau2,Cau3} and Ennassik and Taoudi \cite{ET21} as stated in Subsection \ref{subsection1.1}, their fixed point theorems do not apply to continuous selfmappings on closed $L^0$--convex sets like $L^0(\mathcal{F}, V)$. Just as pointed out in \cite{Guo6,GZWG}, applications of $RN$ modules to stochastic equations and stochastic finance require an appropriate generalization of Schauder fixed point theorem for a noncompact closed $L^0$--convex set in $RN$ modules.

\par
The randomized version of the classical Bolzano--Weierstrass theorem offers a hope of generalizing the Schauder fixed point theorem to $RN$ modules in a noncompact way. The randomized Bolzano--Weierstrass theorem, as presented and used in \cite{KS}, considerably simplifies the proof of no--arbitrage criteria in mathematical finance. It states that for any a.s. bounded sequence $\{\xi^0_n, n\in \mathbb{N} \}$ of random variables from $(\Omega, \mathcal{F}, P)$ to $\mathbb{R}^d$ (namely $\sup \{|\xi^0_n(\omega)|: n\in \mathbb{N} \}< +\infty$ a.s.) there exists a sequence $\{n_k, k\in \mathbb{N}\}$ of positive integer--valued random variables on $(\Omega, \mathcal{F}, P)$ such that $\{n_k(\omega), k\in \mathbb{N} \}$ is strictly increasing for each $\omega\in \Omega$ and $\{\xi^0_{n_k}, k\in \mathbb{N} \}$ converges a.s., where $\xi^0_{n_k} (\omega)= \xi^0_{n_k(\omega) } =\sum^{\infty} _{l=1} I_{(n_k=l)} (\omega) \xi^0_l(\omega)$ for any $k\in \mathbb{N}$ and any $\omega\in \Omega$ (as usual, $(n_k=l)$ denotes the set $\{\omega\in \Omega:n_k(\omega)=l \}$ for any $k$ and $l$ in $\mathbb{N}$). If we use $\xi_i$ for the equivalence class of $\xi^0_i$ for each $i\in \mathbb{N}$ (similarly, $\xi_{n_k}$ for the equivalence class of $\xi^0_{n_k}$ for each $k\in \mathbb{N}$), then $\xi_{n_k} = \sum^{\infty} _{l=1} \tilde{I}_{(n_k=l)} \xi_l$ for each $k\in \mathbb{N}$. Here, recall that a set $G$ of an $RN$ module $(E, \|\cdot\|)$ with base $(\Omega, \mathcal{F}, P)$ is said to be a.s. bounded (or, bounded with respect to the $L^0$--norm $\|\cdot\|$) if there exists $\xi \in L^0_+(\mathcal{F})$ such that $\|x\|\leq \xi$ for any $x$ in $G$ (namely, $\bigvee \{\|x\|: x\in G \}\in L^0_+(\mathcal{F})$). Therefore, the randomized Bolzano--Weierstrass theorem can be restated as: any a.s. bounded nonempty subset of $L^0(\mathcal{F}, \mathbb{R}^d)$ is random relatively sequentially compact in the sense of Definition \ref{definition1.10}, below.

\par
In fact, with the randomized Bolzano--Weierstrass theorem stated above it is not difficult for one to present the notion of a random sequentially compact set in an $RN$ module. For example, in \cite{Guo6} Guo presented such a notion and conjectured the following Theorem \ref{theorem1.12}, which is the best possible result we can prove regarding the theorem in this paper.

\begin{definition}\label{definition1.10}
	Let $(E, \|\cdot\|)$ be an $RN$ module over $\mathbb{K}$ with base $(\Omega, \mathcal{F}, P)$ and $G$ a nonempty subset such that $G$ is contained in a $\sigma$--stable subset $H$ of $E$. Given a sequence $\{x_n, n\in \mathbb{N} \}$ in $G$, a sequence $\{y_k, k\in \mathbb{N} \}$ in $H$ is called a random subsequence of $\{x_n, n\in \mathbb{N} \}$ if there exists a sequence $\{n_k, k\in \mathbb{N} \}$ of positive integer--valued random variables on $(\Omega, \mathcal{F}, P)$ such that the following two conditions are satisfied:
	\begin{enumerate}[(1)]
		\item $n_k(\omega)< n_{k+1} (\omega)$ for any $k\in \mathbb{N}$ and any $\omega\in \Omega$;
		\item $y_k=x_{n_k}:= \sum^{\infty}_{l=1} \tilde{I}_{(n_k=l)} x_l$ for each $k\in \mathbb{N}$.
	\end{enumerate}
	Further, $G$ is said to be random (resp., relatively) sequentially compact if there exists a random subsequence $\{y_k, k\in \mathbb{N}\}$ of $\{x_k, k\in \mathbb{N} \}$ for any sequence $\{x_k, k\in \mathbb{N} \}$ in $G$ such that $\{y_k, k\in \mathbb{N}\}$ converges in the $(\varepsilon,\lambda)$--topology to some element in $G$ (resp., in $E$).
\end{definition}

\par
In Definition \ref{definition1.10}, to ensure that the random subsequence $\{x_{n_k}, k\in \mathbb{N} \}$ is well defined, we require that $G$ be a subset of some $\sigma$--stable subset $H$. The two cases, to which Definition \ref{definition1.10} is frequently applied, are: (1) $E$ itself is $\sigma$--stable, in this case we take $E$ as $H$ and $G$ as any nonempty set; (2) $G$ itself is $\sigma$--stable, in this case we take $H=G$. It is easy to check that for a $\sigma$--stable set $G$, $G$ is random relatively sequentially compact iff $\overline{G}_{\varepsilon,\lambda}$ is random sequentially compact.

\par
For a compact convex set $V$ of a normed space $B$, although $L^0(\mathcal{F}, V)$ is generally not a compact $L^0$--convex subset of the $RN$ module $L^0(\mathcal{F}, B)$, in Section \ref{section2} of this paper one will see that $L^0(\mathcal{F}, V)$ is always a random sequentially compact $L^0$--convex subset of $L^0(\mathcal{F}, B)$. See Corollary \ref{corollary2.6}.

\par
One of the central results of this paper is Theorem \ref{theorem1.12} below. To state Theorem \ref{theorem1.12}, we first give Definition \ref{definition1.11}.

\begin{definition}\label{definition1.11}
	Let $(E_1, \|\cdot\|_1)$ and $(E_2, \|\cdot\|_2)$ be two $RN$ modules over the same scalar field $\mathbb{K}$ with the same base $(\Omega, \mathcal{F}, P)$, and $G_1$ and $G_2$ two nonempty subsets of $E_1$ and $E_2$, respectively, and $T$ a mapping from $G_1$ to $G_2$. $T$ is said to be:
	\begin{enumerate}[(1)]
		\item $\mathcal{T}_{\varepsilon,\lambda}$--continuous if $T$ is a continuous mapping from $(G_1, \mathcal{T}_{\varepsilon,\lambda})$ to $(G_2, \mathcal{T}_{\varepsilon,\lambda})$.
		\item $\mathcal{T}_c$--continuous  if $T$ is a continuous mapping from $(G_1, \mathcal{T}_c)$ to $(G_2, \mathcal{T}_c)$
		\item a.s. sequentially continuous at $x_0\in G_1$ if $\{\|T(x_n)-T(x_0)\|_2, n\in \mathbb{N} \}$ converges a.s. to $0$ for any sequence $\{x_n, n\in \mathbb{N}\}$ in $G_1$ such that $\{\|x_n-x_0\|, n\in \mathbb{N} \}$ converges a.s. to $0$. Further, $T$ is said to be a.s. sequentially continuous if $T$ is a.s. sequentially continuous at any point in $G_1$.
		\item $\sigma$--stable if both $G_1$ and $G_2$ are $\sigma$--stable and $T(\sum^{\infty}_{k=1} \tilde{I}_{A_k} x_k) = \sum^{\infty}_{k=1} \tilde{I}_{A_k} T(x_k)$ for any countable partition $\{A_k, k\in \mathbb{N} \}$ of $\Omega$ to $\mathcal{F}$ and for any sequence $\{x_n, n\in \mathbb{N} \}$ in $G_1$.
		\item random sequentially continuous at $x_0\in G_1$ if $G_1$ is $\sigma$--stable and if for any sequence $\{x_n, n\in \mathbb{N} \}$ in $G_1$ convergent in $\mathcal{T}_{\varepsilon,\lambda}$ to $x_0$ there exists a random subsequence $\{x_{n_k}, k\in \mathbb{N}\}$ of $\{x_n, n\in \mathbb{N} \}$ such that $\{T(x_{n_k}), k\in \mathbb{N} \}$ converges in $\mathcal{T}_{\varepsilon,\lambda}$ to $T(x_0)$. Further, $T$ is said to be random sequentially continuous if $T$ is random sequentially continuous at any point in $G_1$.
	\end{enumerate}
\end{definition}

\par
For (1) and (2) in Definition \ref{definition1.11}, $(G_1, \mathcal{T}_{\varepsilon,\lambda})$ and $(G_1, \mathcal{T}_c)$ denote the topological spaces when $G_1$ is endowed with the relative topologies with respect to the $(\varepsilon,\lambda)$--topology and the locally $L^0$--convex topology for $E_1$, respectively. Likewise, one can understand $(G_2, \mathcal{T}_{\varepsilon,\lambda})$ and $(G_2, \mathcal{T}_c)$. Thus $\mathcal{T}_{\varepsilon,\lambda}$--continuity and $\mathcal{T}_c$--continuity are both a kind of continuity in the sense of standard topology. Although a.s. convergence does not correspond to a topology, a.s. sequentially continuity is already used in \cite{Bhar,DKKS} in an implicit or even explicit way, see Corollaries \ref{corollary1.13} and \ref{corollary1.14} in the following. It is clear that a.s. sequential continuity is stronger than $\mathcal{T}_{\varepsilon,\lambda}$--continuity, and further when $G_1$ is $\sigma$--stable it is also obvious that $\mathcal{T}_{\varepsilon,\lambda}$--continuity of a mapping $T$ from $G_1$ to $G_2$ necessarily implies random sequential continuity of $T$.

\par
Random sequential continuity may appear somewhat strange at first glance, but, to our surprise, in Section \ref{section3} of this paper we will prove that when $G_1$ and $G_2$ are both $\sigma$--stable, a $\sigma$--stable mapping $T$ from $G_1$ to $G_2$ is random sequentially continuous iff $T$ is $\mathcal{T}_c$--continuous. With this equivalence, although the locally $L^0$--convex topology $\mathcal{T}_c$ is too strong to treat fixed point problems for $\mathcal{T}_c$--continuous mappings, we can always convert a fixed point problem of a $\sigma$--stable $\mathcal{T}_c$--continuous mapping to one of a random sequentially continuous mapping. Consequently, the Schauder fixed point theorem on a random sequentially compact $L^0$--convex set---Theorem \ref{theorem1.12} holds true for both a $\sigma$--stable $\mathcal{T}_c$--continuous mapping and a $\sigma$--stable $\mathcal{T}_{\varepsilon,\lambda}$--continuous mapping. Besides, since the $(\varepsilon,\lambda)$--topology $\mathcal{T}_{\varepsilon,\lambda}$ for an $RN$ module is metrizable and a random sequentially compact $L^0$--convex set is not compact in general, from the end of Subsection \ref{subsection1.1} one can not rely on the fact that Theorem \ref{theorem1.12} is always available for every $\mathcal{T}_{\varepsilon,\lambda}$--continuous mapping. Thus Theorem \ref{theorem1.12} is the best possible in this sense.

\par
As a reminder, we note here that in Section \ref{section2} we will prove that a random sequentially compact $L^0$-convex set is always $\sigma$--stable and, thus, a random sequentially continuous mapping on it is well defined.

\begin{theorem}\label{theorem1.12}
	Let $(E, \|\cdot\|)$ be an $RN$ module with base $(\Omega, \mathcal{F}, P)$ and $G$ a random sequentially compact $L^0$--convex subset of $E$. Then every $\sigma$--stable and random sequentially continuous mapping $T$ from $G$ to $G$ has a fixed point.
\end{theorem}

\par
Clearly, when $(\Omega, \mathcal{F}, P)$ is trivial, namely $\mathcal{F}= \{\Omega, \emptyset\}$, $(E, \|\cdot\|)$ reduces to an ordinary normed space, $G$ to a compact convex set of $E$ and $T$ to an ordinary continuous mapping, and hence the classical Schauder fixed point theorem is a special case of Theorem \ref{theorem1.12}.

\par
Historically, the first famous random generalization of the classical Schauder fixed point theorem is Corollary \ref{corollary1.13} below:

\begin{corollary}[{\cite[Theorem 10]{Bhar}}]\label{corollary1.13}
	Let $X$ be a compact convex subset of a normed space $B$ and $T: \Omega \times X\rightarrow X$ a sample--continuous random operator. Then $T$ has a random fixed point.
\end{corollary}

\par
In 1976, Bharucha--Reid and Mukheriea proved Corollary \ref{corollary1.13} by measurable selection theorems of multi--valued measurable functions. In \cite{Bhar}, they define $A: \Omega\rightarrow 2^X$ as follows:
$$A(\omega)= \{x\in X: T(\omega,x)=x\} ~\text{for each}~\omega\in \Omega.$$
Then $A(\omega)\neq \emptyset$ for each $\omega\in \Omega$ by the classical Schauder fixed point theorem. They further prove that $A$ is measurable and construct a measurable selection $v$ of $A$ so that $v$ is a random fixed point of $T$. In fact, such a $v$ always exists by Theorem 4.1 of \cite{Wagn}.

\par
However, in this paper we can use Theorem \ref{theorem1.12} to immediately give Corollary \ref{corollary1.13} a simpler proof independent of any measurable selection theorems: since any $X$--valued random element $g^0: \Omega\rightarrow X$ must be a strong random element by the separability of $X$ and that $T(\omega, \cdot): X\rightarrow X$ is continuous for each $\omega\in \Omega$; $T(\cdot, g^0(\cdot)): \Omega\rightarrow X$ is clearly also an $X$--valued strong random element, and thus this induces a mapping $\hat{T}: L^0(\mathcal{F}, X)\rightarrow L^0(\mathcal{F}, X)$ as follows:
\par
$\hat{T}(g)=$ the equivalence class of $T(\cdot, g^0(\cdot))$ for any $g\in L^0(\mathcal{F}, X)$, where $g^0$ is an arbitrarily chosen representative of $g$.\\
Then $L^0(\mathcal{F}, X)$ is a random sequentially compact $L^0$--convex subset of the $RN$ module $L^0(\mathcal{F}, B)$ (see Corollary \ref{corollary2.6}) and it is also easy to check that $\hat{T}$ is a $\sigma$--stable and a.s. sequentially continuous mapping. So, by Theorem \ref{theorem1.12} $\hat{T}$ has a fixed point $g$ in $L^0(\mathcal{F}, X)$, an arbitrarily chosen representative $g^0$ of $g$ is, of course, a random fixed point of $T$, namely $T(\omega, g^0(\omega))= g^0(\omega)$ a.s..

\par
In 2013, Drapeau, et al \cite{DKKS} first gave a random generalization of the classical Brouwer fixed point theorem in the random Euclidean space $L^0(\mathcal{F}, \mathbb{R}^d)$, namely they obtained Corollary \ref{corollary1.14} below, which is merely a restatement of Proposition 3.1 of \cite{DKKS}.

\begin{corollary}[{\cite[Proposition 3.1]{DKKS}}]\label{corollary1.14}
	Let $G$ be an a.s. bounded $\mathcal{T}_{\varepsilon,\lambda}$--closed $L^0$--convex subset of $L^0(\mathcal{F}, \mathbb{R}^d)$ (clearly, where ``~$\mathcal{T}_{\varepsilon,\lambda}$--closed'' is just closed in the topology of convergence in probability measure). Then every $\sigma$--stable a.s. sequentially continuous mapping $T$ from $G$ to $G$ has a fixed point.
\end{corollary}

\par
According to the randomized Bolzano--Weierstrass theorem \cite{KS}, $G$ in Corollary \ref{corollary1.14} is random sequentially compact, while an a.s. sequentially continuous mapping is obviously random sequentially continuous. So Corollary \ref{corollary1.14} is also a special case of Theorem \ref{theorem1.12}.

\par
Besides Corollary \ref{corollary1.14}, another nice contribution of \cite{DKKS} is the definition of an $L^0$--simplex (or, a conditional simplex in the terminology in \cite{DKKS}) and its random labeling function. The work in \cite{DKKS} enables a further study of its barycentric subdivision and the establishment of the $L^0$--version of Sperner's lemma (see the first half of the proof of Theorem 2.3 of \cite{DKKS}). By modifying the second half of the proof of Theorem 2.3 of \cite{DKKS}, we can prove that Corollary \ref{corollary1.14} remains valid for a random sequentially continuous mapping, namely, we can obtain a complete Brouwer fixed point theorem in $L^0(\mathcal{F}, \mathbb{R}^d)$ (see Lemmas \ref{lemma3.6} and \ref{lemma3.7} of this paper), which facilitates the first step in the proof of Theorem \ref{theorem1.12}.

\par
To outline the key second step in the proof of Theorem \ref{theorem1.12}, first recall the idea of the proof of the classical Schauder fixed point theorem. The classical Hausdorff theorem states that a nonempty subset in a normed space (or more generally, a metric space) is sequentially compact iff it is both totally bounded and complete. It is the total boundedness of a compact convex set that allows Schauder to define the Schauder projection on a finite $\varepsilon$--net of the convex set and further define the Schauder approximation for a compact mapping so that the classical Brouwer fixed point theorem can be employed to prove the classical Schauder fixed point theorem, see \cite[p.117]{GD} for details. Since a random sequentially compact set in an $RN$ module is, generally, not compact, the classical notion of a total boundedness is impossible to fit together with the notion of random sequential compactness in order for us to establish a random version of the Hausdorff theorem. For this, we introduce the notion of random total boundedness in Section \ref{section2}, and further prove that a $\sigma$--stable set is random sequentially compact iff it is both random totally bounded and $\mathcal{T}_{\varepsilon,\lambda}$--complete (see Theorem \ref{theorem2.3} of this paper). Based on the random Hausdorff theorem, we can first decompose the mapping $T$ in Theorem \ref{theorem1.12} into a sequence $\{T_n, n\in \mathbb{N} \}$ of local mappings according to some countable partition $\{A_n, n\in \mathbb{N}\}$ of $\Omega$ to $\mathcal{F}$, then construct a random Schauder approximation for each $T_n$ in a subtle manner, namely, in a manner rather different from the classical Schauder approximation. See the construction before Lemma \ref{lemma3.8} of this paper so that we can obtain an approximate fixed point $x_n$ for each $T_n$, at the same time $\sum_n \tilde{I}_{A_n} x_n$ is exactly an approximate fixed point of $T$, and eventually achieve the proof of Theorem \ref{theorem1.12}.

\par
As usual, Theorem \ref{theorem1.12} also leads directly to the following random generalization of the classical Krasnoselskii fixed point theorem in \cite{KRAS}.

\begin{theorem}\label{theorem1.15}
	Let $(E, \|\cdot\|)$ be a $\mathcal{T}_{\varepsilon,\lambda}$--complete $RN$ module with base $(\Omega, \mathcal{F}, P)$, $G$ a $\mathcal{T}_{\varepsilon,\lambda}$--closed $L^0$--convex subset of $E$, $S$ and $T$ two mappings from $G$ to $E$ such that the following three conditions are satisfied:
	\begin{enumerate}[(1)]
		\item $S$ is random contractive, namely there exists $\alpha\in L^0_+(\mathcal{F})$ satisfying $\alpha<1$ on $\Omega$ such that $\|S(x)-S(y)\|\leq \alpha \|x-y\|$ for all $x$ and $y$ in $G$ ($S$ must be $\sigma$--stable by Lemma 2.11 of \cite{GZWG});
		\item $T$ is a $\sigma$--stable and random sequentially continuous mapping such that $T(G)$ is random relatively sequentially compact;
		\item $S(x)+T(y)\in G$ for any $x$ and $y$ in $G$.
	\end{enumerate}
	Then $S+T$ has a fixed point.
\end{theorem}

\par
At the end of this subsection, we would also like to point out an interesting fact: both Theorem \ref{theorem1.12} and Corollary \ref{corollary1.13} are a random generalization of the classical Schauder fixed point theorem, but the former comes from the viewpoint of randomizing space theory and the latter originates from random operator theory. It is also obvious that the former is much more general and deeper than the latter!

\subsection{The existence of Nash equilibrium points in the framework of conditional information}\label{subsection1.4}
\par
The classical Nash theorem on the existence of equilibrium points consider $n$--fold affine (or, linear) functions defined on $n$--fold Cartesian product of simplexes, where $n$ stands for the number of players. See \cite{Nash} for details. The finite (non--coorperative) game considered in \cite{Nash} is essentially static. In most dynamic models, information accumulates over time and players in a finite game naturally make a mixed strategy contingent on information available at the time a mixed strategy is made, so the accumulated information becomes imbedded in mixed strategies and payoff functions of players. Just motivated by the work of \cite{HR} on the conditional mean--variance frontier and that of \cite{DS05,FKV2,GZZ3,FM14} on the conditional risk measures, in this paper we consider the finite game under the conditional framework (or simply, the finite conditional game). As an interesting application of Theorem \ref{theorem1.12} we obtain the conditional version of the classical existence theorem of Nash equilibrium points, namely Theorem \ref{theorem1.19} below, which is a natural generalization of Theorem 1 of \cite{Nash}.

\par
To describe the finite conditional game, we first introduce the following standard terminologies and notation.

\begin{definition}\label{definition1.16}
	Let $(E, \|\cdot\|)$ be an $RN$ module over $\mathbb{K}$ with base $(\Omega, \mathcal{F}, P)$ and $G$ a nonempty subset of $E$. The set $Conv_{L^0} (G):= \{\sum^n_{i=1} \xi_i x_i: n\in \mathbb{N}, \{x_1,x_2,\cdots, x_n \}\subset G~\text{and}~\{\xi_1,\xi_2,\cdots,\xi_n \}\subset L^0_+(\mathcal{F})~\text{such that}~\sum^n_{i=1} \xi_i=1 \}$ is called the $L^0$--convex hull of $G$. The finite subset $\{x_1,x_2,\cdots, x_n \}$ of $G$ is said to be $L^0$--independent if $\sum^n_{i=1} \xi_i x_i=0$ for some $\xi_1,\xi_2,\cdots$ and $\xi_n$ in $L^0(\mathcal{F}, \mathbb{K})$ implies $\xi_i=0$ for each $i=1\sim n$. Further, $G$ is said to be $L^0$--independent if each finite nonempty subset of $G$ is $L^0$--independent, while the $L^0$--submodule generated by $G$ is said to be $L^0$--free. Finally, the finite subset $\{x_1,x_2,\cdots, x_n \}$ of $E$ is said to be $L^0$--affinely independent if either $n=1$ or if $n>1$ and $\{x_i-x_1: i=2\sim n \}$ is $L^0$--independent, while $Conv_{L^0} (\{x_1, x_2,\cdots, x_n \})$ is called an $(n-1)$ dimensional $L^0$--simplex in $E$.
\end{definition}

\par
Obviously, for an $(n-1)$ dimensional $L^0$--simplex $S:= Conv_{L^0} (\{x_1, x_2,\cdots, x_n \})$, $\sum^n_{i=1} \xi_i x_i = \sum^n_{i=1} \eta_i x_i$ for $\{\xi_i: i=1\sim n \}$ and
$\{\eta_i: i=1\sim n \}\subset L^0_+(\mathcal{F})$ such that $\sum^n_{i=1} \xi_i = \sum^n_{i=1} \eta_i=1$ iff $\xi_i= \eta_i$ for each $i=1\sim n$, so each element $x\in S$ can be uniquely represented as $x= \sum^n_{i=1} \xi_i x_i$ with $\xi_i\in L^0_+(\mathcal{F})$ such that $\sum^n_{i=1} \xi_i =1$.

\par
In what follows, we always use $L^0(\mathcal{F},[0,1] )$ for the set of equivalence classes of random variables form $(\Omega, \mathcal{F}, P)$ to $[0,1]$. Similarly, we do for $L^0(\mathcal{F},(0,1) )$. It is obvious that $L^0(\mathcal{F},[0,1] )$ and $L^0(\mathcal{F},(0,1) )$ are both $L^0$--convex subsets of $L^0(\mathcal{F})$.

\begin{definition}\label{definition1.17}
	Let $(E, \|\cdot\|)$ be an $RN$ module over $\mathbb{K}$ with base $(\Omega, \mathcal{F}, P)$ and $G$ an $L^0$--convex subset of $E$. A mapping $f: G\rightarrow L^0(\mathcal{F})$ is called an $L^0$--affine (or, $L^0$--linear) function if $f(\lambda g_1+ (1-\lambda)g_2)= \lambda f(g_1)+ (1-\lambda) f(g_2)$ for any $\lambda\in L^0(\mathcal{F},[0,1] )$ and any $g_1$ and $g_2$ in $G$.
\end{definition}

\par
Let $(E_1, \|\cdot\|_1), (E_2, \|\cdot\|_2), \cdots $ and $(E_n, \|\cdot\|_n)$ be $n$ $RN$ modules over $\mathbb{K}$ with base $(\Omega, \mathcal{F}, P)$. It is easy to see that $E= \Pi^n_{i=1} E_i$ (namely, the product space of $E_1,E_2,\cdots$ and $E_n$) is still an $RN$ module over $\mathbb{K}$ with base $(\Omega, \mathcal{F}, P)$ under the $L^0$--norm defined by $\|x\|= \sqrt{\sum^n_{i=1}\|x_i\|_i^2 }$ for any $x=(x_1,x_2,\cdots, x_n)\in E$. For any $x\in E$ and any $y_i\in E_i$, we use the abbreviated notation $(x,y_i)$ for $(x_1, x_2,\cdots, x_{i-1}, y_i, x_{i+1}, \cdots, x_n)$ for each $i=1\sim n$, in case no ambiguities arise.

\begin{definition}\label{definition1.18}
	Let $(E_1, \|\cdot\|_1), (E_2, \|\cdot\|_2), \cdots $ and $(E_n, \|\cdot\|_n)$ be $RN$ modules over $\mathbb{K}$ with base $(\Omega, \mathcal{F}, P)$ and $G_i$ an $L^0$--convex subset of $E_i$ for each $i=1\sim n$. A mapping $f: G= \Pi^n_{i=1} G_i$ (the product space of $G_1,G_2, \cdots$ and $G_n$) $\rightarrow L^0(\mathcal{F})$ is called an $n$--fold $L^0$--affine (or, $L^0$--linear) function if for each $i= 1 \sim n$ and any fixed $x\in G$ the mapping $f(x_1, x_2, \cdots, x_{i-1}, y_i, x_{i+1}, \cdots, x_n)$ is an $L^0$--affine function of $y_i\in G_i$, namely $f(x_1, x_2, \cdots, x_{i-1},\cdot, x_{i+1}, \cdots, x_n): G_i \rightarrow L^0(\mathcal{F})$ is an $L^0$--affine function for each $i=1\sim n$.
\end{definition}

\par
Now, we can consider the $n$--person conditional game as follows: the probability space $(\Omega, \mathcal{F}, P)$ stands for the random model, where $\mathcal{F}$ stands for the information available at the time the $n$ players make their mixed strategies; for each player $i$ the $RN$ module $(E_i, \|\cdot\|_i)$ over $\mathbb{K}$ with base $(\Omega, \mathcal{F}, P)$ stands for the space of his random payoffs (for example, $E_i$ is the conditional $L^r$--space---$L^r_{\mathcal{F}} (\mathcal{E})$ (cf. Example \ref{example1.9}) or an $L^0$--submodule of $L^r_{\mathcal{F}} (\mathcal{E})$), the $(n_i-1)$--dimensional $L^0$--simplex $S_i: = Conv _{L^0} (\{\pi_{i\alpha}: \alpha=1\sim n_i \})$ of $E_i$ for the set of his mixed strategies with the $L^0$--affinely independent set $\{\pi_{i\alpha}: \alpha=1\sim n_i \}$ as the set of his pure strategies, and an $n$--fold $L^0$--affine function $p_i: S=\Pi^n_{i=1} S_i \rightarrow L^0(\mathcal{F})$ for his $L^0$--payoff function. We call $s= (s_1, s_2, \cdots, s_n)\in S$ a Nash equilibrium point of the $n$--person conditional game if $p_i(s_1, s_2, \cdots, s_n) \geq p_i(s_1, s_2, \cdots, s_{i-1}, t_i, s_{i+1}, \cdots s_n)$ for each $i=1\sim n$ and each $t_i\in S_i$, namely $p_i(s_1, s_2, \cdots, s_n)=\max \{p_i(s_1, s_2, \cdots, s_{i-1}, t_i, s_{i+1}, \cdots s_n): t_i\in S_i \}$ for each $i=1\sim n$.

\par
As a natural generalization of Theorem 1 of \cite{Nash}, we have the following.

\begin{theorem}\label{theorem1.19}
	Every $n$--person conditional game has a Nash equilibrium point.
\end{theorem}

\section{Random sequential compactness and random total boundedness}\label{section2}
\par
To proceed forward, we introduce the following notation for brevity
\begin{enumerate}[$\Pi_{\mathcal{F}}:$]
	\item the set of countable partitions of $\Omega$ to $\mathcal{F}$.
\end{enumerate}

\begin{enumerate}[$G_{(\varepsilon)}$]
	\item $:=\bigcup_{x\in G} B(x,\varepsilon)$ for a nonempty subset $G$ of some $RN$ module $(E,\|\cdot\|)$ with base $(\Omega, \mathcal{F},P)$, where $\varepsilon\in L^0_{++}(\mathcal{F})$ and $B(x,\varepsilon)=\{y\in E: \|y-x\|< \varepsilon$ on $\Omega\}$.
\end{enumerate}
\par
As an aid for Definition \ref{definition2.1} below, assume that $(E, \|\cdot\|)$ is an $RN$ module with base $(\Omega,\mathcal{F},P)$ and $H$ a nonempty finite subset of $E$. It is easy to see, by letting $H=\{x_1, x_2,\cdots, x_n\}$, that the set $\{\sum^n_{k=1} \tilde{I}_{A_k} x_k: \{A_1, A_2, \cdots, A_n\}$ is an n-partition of $\Omega$ to $\mathcal{F}$ $\}$ is always a $\sigma$-stable subset of $E$ and it is just $\sigma(H)$. So $\sigma(H)$ is always well-defined for any nonempty finite subset of an $RN$ module. Hence, $\tilde{I}_A G=\{\tilde{I}_A\cdot g: g\in G\}$ for any $A\in \mathcal{F}$ and $G\subset E$.

\begin{definition}\label{definition2.1}
	Let $(E,\|\cdot\|)$ be an $RN$ module with base $(\Omega,\mathcal{F},P)$ and $G$ a nonempty subset of $E$. $G$ is said to be random totally bounded if for any $\varepsilon\in L^0_{++}(\mathcal{F})$ there exist $\{A_n,n\in \mathbb{N}\}$ in $\Pi_{\mathcal{F}}$ and a sequence $\{G_n,n\in \mathbb{N}\}$ of nonempty finite subsets of $G$ such that $\tilde{I}_{A_n} G\subset \tilde{I}_{A_n} (\sigma(G_n))_{(\varepsilon)}$ for each $n\in \mathbb{N}$.
\end{definition}

\begin{remark}\label{remark2.2}
	In Definition \ref{definition2.1}, if $(E,\|\cdot\|)$ is a $\sigma$-stable $RN$ module and if we define $d:E\times E\rightarrow L^0_+(\mathcal{F})$ by $d(x,y)=\|x-y\|$, then $(E,d)$ is a $d$--$\sigma$--stable random metric space, see \cite{GWYZ} for the related concept and discussions. Further, if we denote $\sum_n \tilde{I}_{A_n}\sigma(G_n)$ $:=\{\sum_n \tilde{I}_{A_n} y_n: y_n\in \sigma(G_n)$ for each $n\in \mathbb{N} \}$by $G'$ for any fixed $\{A_n, n\in \mathbb{N}\}$ in $\Pi_{\mathcal{F}}$, then it is easy to see that $\sum_n \tilde{I}_{A_n}(\sigma(G_n))_{(\varepsilon)}=\bigcup_{x\in G'}B(x,\varepsilon)$. So the notion of random total boundedness can also be defined for a nonempty subset of a $d$--$\sigma$--stable random metric space. In fact, all the results in this section also hold for $d$--$\sigma$--stable random metric spaces. Although a more general notion, called conditional total boundedness, was first introduced in \cite{DJKK} for conditional metric spaces and further discussed in \cite{JZ} for $L^0$--modules, these papers employed some abstract notions such as conditional join operations and conditional finite sets so that it is not easy for the reader to grasp the notion of conditional total boundedness. This paper is only interested in $RN$ modules and our formulation in Definition \ref{definition2.1} is, we believe, both more accessible and more convenient for use since it is well defined for any nonempty subsets.
\end{remark}
\par
The central result of this section is the following.
\begin{theorem}\label{theorem2.3}
	Let $(E,\|\cdot\|)$ be an $RN$ module with base $(\Omega,\mathcal{F},P)$ and $G$ a $\sigma$--stable subset of $E$. Then $G$ is random sequentially compact iff it is both random totally bounded and $\mathcal{T}_{\varepsilon,\lambda}$--complete.
\end{theorem}

\par
Just as pointed out in Remark \ref{remark2.2}, Theorem \ref{theorem2.3} remains valid for a $d$--$\sigma$--stable random metric space. When $(\Omega,\mathcal{F},P)$ is trivial Theorem \ref{theorem2.3} reduce to the classical Hausdorff theorem. It is quite obvious that the notion of random total boundedness is an outgrowth of the locally $L^0$--convex topology, whereas the notion of random sequential compactness is an outgrowth of the $(\varepsilon,\lambda)$--topology. Such results as Theorem \ref{theorem2.3} are most useful in random functional analysis since it establishes a connection between some basic results derived from the two kinds of topologies, and in particular, they will play a crucial role in constructing the Schauder projection and in establishing Schauder fixed point theorem in $RN$ modules by simultaneously considering the two kinds of topologies, see the proofs of Lemma \ref{lemma3.8}, Lemma \ref{lemma3.9} and Theorem \ref{theorem1.12} for details.

\par
In the literature of stochastic finance or probability theory, there is usually an equivalent variant of Proposition \ref{proposition1.1}: let $\bar{\mathcal{L}}^0(\mathcal{F})$ be the set of extended real--valued $\mathcal{F}$--measurable functions on $(\Omega,\mathcal{F},P)$ and $G$ a nonempty subset of $\bar{\mathcal{L}}^0(\mathcal{F})$, $\eta\in \bar{\mathcal{L}}^0(\mathcal{F})$ is called an essential upper bound of $G$ if $\xi(\omega)\leq \eta(\omega)$ a.s. for each $\xi\in G$. Further, an essential upper bound $\eta$ of $G$ is called an essential supremum of $G$ if $\eta(\omega)\leq \eta'(\omega)$ a.s. for each essential upper bound $\eta'$ of $G$. Clearly, the essential supremum of $G$ is unique in the sense of almost sure equality, denoted by $esssup(G)$. Likewise, one can define an essential infimum for $G$, denoted by $essinf(G)$. The so--called essential supremum (resp., infimum) theorem states that every nonempty subset of $\bar{\mathcal{L}}^0(\mathcal{F})$ has an essential supremum (resp., infimum), and the essential supremum (resp., infimum) has the similar properties as in Proposition \ref{proposition1.1}. For a subfamily $\mathcal{A}$ of $\mathcal{F}$, there exists $B\in \mathcal{F}$ such that $I_B=esssup\{I_A: A\in \mathcal{A}\}$ a.s., and $B$ is called an essential supremum of $\mathcal{A}$, denoted by $esssup(\mathcal{A})$. Likewise, one can have the notion of an essential infimum of $\mathcal{A}$.

\par
It should be pointed out that in this paper we will strictly distinguish the
equivalence class of a measurable function from this function itself since we need to consider some applications of $RN$ modules to random operators, which means that we will necessarily give up the convention that measurable functions equal a.s. are identified, and for this we only need to give some proper interpretations for some terminologies and propositions about equivalence classes: for example, we say that a sequence $\{\xi_n, n\in \mathbb{N}\}$ in $L^0(\mathcal{F})$ converges uniformly to $\xi_0\in L^0(\mathcal{F})$ on some $A\in \mathcal{F}$ if for an arbitrarily chosen representative $\xi^0_n$ of $\xi_n$ for each $n\geq 0$ there exists a measurable subset $A'$ of $A$ with $P(A\backslash A')=0$ such that the resulting sequence $\{\xi^0_n: n\in \mathbb{N}\}$ converges uniformly to $\xi^0_0$ on $A'.$
\par
For any $A\in \mathcal{F}$, let $\tilde{A}=\{B\in \mathcal{F}: P(A\bigtriangleup B)=0\}$, where $A\bigtriangleup B=(A\backslash B)\bigcup(B\backslash A)$, we also often use $I_{\tilde{A}}$ for $\tilde{I}_A$. Let $\xi$ and $\eta$ be any two elements in $\bar{L}^0(\mathcal{F})$, for arbitrarily respectively chosen representatives $\xi^0$ and $\eta^0$ of $\xi$ and $\eta$. Let $A=\{\omega\in \Omega: \xi^0(\omega)< \eta^0(\omega)\}$; we usually denote $\tilde{A}$ by $(\xi< \eta)$. Similarly, one can also write $(\xi\leq \eta)$. Besides, let $(E,\|\cdot\|)$ be an $RN$ module, $G$ a nonempty subset of $E$ and $x\in E$, $d(x,G):=\bigwedge\{\|x-y\|: y\in G\}$, called the random distance from $x$ to $G$. As in the classical case, it is easy to verify that $|d(x,G)-d(y,G)|\leq \|x-y\|$ for any $x$ and $y\in E$.
\par
For the proof of Theorem \ref{theorem2.3}, we first give the following lemmas \ref{lemma2.4} and \ref{lemma2.5}.
\begin{lemma}\label{lemma2.4}
	Let $(E,\|\cdot\|)$ be an $RN$ module with base $(\Omega,\mathcal{F},P)$ and $G$ a nonempty subset. Then the following statements hold:
	\begin{enumerate}[$(1)$]
		\item $d(x,G)=d(x,\sigma (G))$ for any $ x\in E$ when $G$ is contained in some $\sigma$--stable subset of $E$.\\
		\item $G_{(\varepsilon)}=\{x:d(x,G)<\varepsilon $ on $\Omega\}$ for any $\varepsilon\in L^0_+(\mathcal{F})$ when $G$ is $\sigma$--stable.\\
		\item When $G=\{x_1,x_2,x_3,\cdots,x_n\}$, namely $G$ is a finite set, for any given $x\in E$ there exists a finite partition $\{A_1,A_2,\cdots,A_n\}$ of $\Omega$ to $\mathcal{F}$ such that $d(x,G)=\|x-\sum^n_{k=1} \tilde{I}_{A_k}x_k\|.$\\
		\item When $G$ is a random sequentially compact $\sigma$--stable subset and $f:G\rightarrow L^0(\mathcal{F})$ is $\sigma$--stable and random sequentially continuous, there exist $x_0$ and $y_0$ in $G$ such that $f(x_0)=\bigvee f(G)$ and $f(y_0)=\bigwedge f(G).$\\
		\item If $G$ is random sequentially compact and $\sigma$--stable, then $G$ is a.s. bounded.\\
		\item If $G$ is random totally bounded, then $G$ is a.s. bounded.
		
	\end{enumerate}
\end{lemma}
\begin{proof}
	(1). It is obvious that $d(x,G)\geq d(x,\sigma (G))$. Conversely, for any sequence $\{x_n,n\in \mathbb{N}\}$ in $G$ and any $\{A_n,n\in \mathbb{N}\}$ in $\Pi_{\mathcal{F}}$, $\|x-\sum_n\tilde{I}_{A_n}x_n\|=\sum_n\tilde{I}_{A_n}\|x-x_n\|\geq \sum_n\tilde{I}_{A_n} d(x,G)=d(x,G)$, so $d(x,\sigma(G))\geq d(x,G)$.
	\par
	(2). For any $x\in G_{(\varepsilon)}$, there exists $y\in G$ such that $\|x-y\|<\varepsilon$ on $\Omega$, which in turn implies that $d(x,G)\leq \|x-y\|<\varepsilon$ on $\Omega$. Conversely, if $x\in E$ is such that $d(x,G)<\varepsilon$ on $\Omega$, then by Proposition \ref{proposition1.1} there exists a sequence $\{y_n, n\in \mathbb{N}\}$ in $G$ such that $\{\|x-y_n\|, n\in \mathbb{N}\}$ converges to $d(x,G)$ in a nonincreasing way since $\{\|x-y\|: y\in G\}$ is directed downwards. In fact, for any $y_1$ and $y_2$ in $G$, let $A=(\|x-y_1\|<\|x-y_2\|)$, then it is easy to see that $\|x-y_1\|\bigwedge\|x-y_2\|=I_A\|x-y_1\|+(1-I_A)\|x-y_2\|=\|x-(I_Ay_1+(1-I_A)y_2)\|=\|x-y_3\|$, where $y_3=I_Ay_1+(1-I_A)y_2\in G$. By the Egoroff theorem there exists an $F_k\in \mathcal{F}$ for each $k\in \mathbb{N}$ such that $\{\|x-y_n\|: n\in \mathbb{N}\}$ converges uniformly to $d(x,G)$ on $F_k$ and $P(\Omega\backslash F_k)<\frac{1}{k}$. Thus there exists $n_k\in \mathbb{N}$ for each $k\in \mathbb{N}$ such that $\|x-y_{n_k}\|<\varepsilon$ on $F_k$. Since $P(\bigcup^{\infty}_{k=1} F_k)=1$, put $A_k=F_k\backslash (\bigcup^{k-1}_{i=1} F_i)$ for any $k\geq 2$ and $A_1=F_1$, then $\{A_k,k\in \mathbb{N}\}$ forms a countable partition of $\Omega$ to $\mathcal{F}$. Let $y=\sum^{\infty}_{k=1}\tilde{I}_{A_k}y_{n_k}$, then $y\in G$ since $G$ is $\sigma$--stable, and it is obvious that $\|x-y\|<\varepsilon$ on $\Omega$, namely $x\in B(y,\varepsilon)\subset \bigcup_{z\in G}B(z,\varepsilon)$, that is, $x\in G_{(\varepsilon)}$.
	\par
	(3). When $n=2, d(x,G)=\|x-y_1\|\bigwedge \|x-y_2\|$. Let $A$ be a representative of $(\|x-y_1\|<\|x-y_2\|)$, then it is obvious that $\|x-y_1\|\bigwedge \|x-y_2\|=\|x-(\tilde{I}_Ay_1+\tilde{I}_{A^c}y_2)\|$. By induction, assume that (3) holds for $n=k$, namely, there exists a finite partition $\{A_i: i=1\sim k\}$ of $\Omega$ to $\mathcal{F}$ such that $\bigwedge\{\|x-y_i\|: i=1 \sim k\}=\|x-(\sum^k_{i=1}\tilde{I}_{A_i}x_i)\|$, then $\bigwedge\{\|x-y_i\|: i=1 \sim k+1\}=\|x-(\sum^k_{i=1}\tilde{I}_{A_i}x_i)\|\bigwedge \|x-x_{k+1}\|$, by the proved case when $n=2$ one can have that there exists $B\in \mathcal{F}$ such that $\|x-(\sum^k_{i=1}\tilde{I}_{A_i}x_i)\|\bigwedge \|x-x_{k+1}\|=\|x-(I_B(\sum^k_{i=1}\tilde{I}_{A_i}x_i)+\tilde{I}_{B^c} x_{k+1})\|=\|x-\sum^{k+1}_{i=1}\tilde{I}_{D_i}x_i\|$, where $D_i=B\bigcap A_i$ for each $i=1 \sim k$ and $D_{k+1}=B^c$, it is also clear that $\{D_i,i=1\sim k+1\}$  forms a finite partition of $\Omega$ to $\mathcal{F}$, which ends the proof of (3).
	\par
	(4). We only need to prove the existence of $x_0$ since the existence of $y_0$ is similar. We first prove that $f(G)$ is directed upwards as follows: for any $y_1$ and $y_2$ in $G$, let $A=(f(y_1)>f(y_2))$ and $y_3=I_A y_1+(1-I_A)y_2$, then $y_3\in G$ and it is easy to check that $f(y_3)=f(y_1)\bigvee f(y_2)$. Thus, by Proposition \ref{proposition1.1} there exists a sequence $\{x_n, n\in \mathbb{N}\}$ in $G$ such that $\{f(x_n), n\in \mathbb{N}\}$ converges to $\bigvee f(G)$ in a nondecreasing way. Since $G$ is $\sigma$--stable and random sequentially compact, there exist $x_0\in G$ and a random subsequence $\{y_k, k\in \mathbb{N}\}$ in $G$ of $\{x_n, n\in \mathbb{N}\}$ such that $\{y_k, k\in \mathbb{N}\}$ converges in $\mathcal{T}_{\varepsilon,\lambda}$ to $x_0$, where $y_k=x_{n_k}$ for each $k\in \mathbb{N}$. We can, without loss of generality, assume that $n_k(\omega)\geq k$ for each $k\in \mathbb{N}$ and $\omega\in \Omega$. Then by $\sigma$--stability of $f$, $f(y_k)=f(x_{n_k})=f(\sum^{\infty}_{l=k}\tilde{I}_{(n_k=l)}x_l)=\sum^{\infty}_{l=k}\tilde{I}_{(n_k=l)} f(x_l)\geq f(x_k)$ for each $k\in \mathbb{N}$. Since $f$ is random sequentially continuous, there exists a random subsequence $\{z_k, k\in \mathbb{N}\}$ in $G$ of $\{y_k, k\in \mathbb{N}\}$ such that $\{f(z_k), k\in \mathbb{N}\}$ converges in probability to $f(x_0)$, where $z_k=y_{n'_k}$ for each $k\in \mathbb{N}$ and we can also assume that $n'_k(\omega)\geq k$ for each $k\in \mathbb{N}$ and $\omega\in \Omega$ so that we can similarly have that $f(z_k)\geq f(y_k)$ for each $k\in \mathbb{N}$. Eventually, we have $f(x_0)=\bigvee f(G)$.
	\par
	(5). Taking $f=\|\cdot\|$ in (4) yields $x_0\in G$ such that $\|x_0\|=\bigvee \{\|x\|: x\in G\}$.
	\par
	(6). Let $\varepsilon,\{G_n,n\in \mathbb{N}\}$ and $\{A_n,n\in \mathbb{N}\}$ be the same as in Definition \ref{definition2.1}. Taking $\varepsilon=1$ yields that $\tilde{I}_{A_n}\|x\|\leq \bigvee \tilde{I}_{A_n}(\{\|y\|: y\in \sigma(G_n)\}+1)$ for each $x\in G$ and each $n\in \mathbb{N}$. Since it is also obvious that $\bigvee\{\|y\|: y\in \sigma(G_n)\}= \bigvee\{\|y\|: y\in G_n\}$, then $\xi_n: =\bigvee\{\|y\|: y\in G_n\}\in L^0_+(\mathcal{F})$ since $G_n$ is a finite subset of $G$, which in turn implies that $\|x\|\leq \sum_n \tilde{I}_{A_n}\xi_n:= \xi$ for each $x\in G$, namely $G$ is a.s. bounded by noticing $\xi\in L^0_+(\mathcal{F})$.
\end{proof}

\begin{lemma}\label{lemma2.5}
	Let $(E,\|\cdot\|)$ be an $RN$ module with base $(\Omega, \mathcal{F}, P)$ and $G$ a nonempty subset of $E$, then the following hold:
	\begin{enumerate}[$(1)$]
		\item If $G$ is a stable and $\mathcal{T}_{\varepsilon,\lambda}$--closed subset of some $\sigma$-stable subset of $E$, then $G$ is $\sigma$--stable.
		\item If a sequence $\{x_n, n\in \mathbb{N}\}$ in some $\sigma$-stable subset of $E$ converges in $\mathcal{T}_{\varepsilon,\lambda}$ to $x_0\in E$, then any random subsequence $\{x_{n_k}, k\in \mathbb{N}\}$ of $\{x_n, n\in \mathbb{N}\}$ also converges in $\mathcal{T}_{\varepsilon,\lambda}$ to $x_0$.
		\item $G$ is random totally bounded iff $\bar{G}_c$ is random totally bounded.
		\item If $E$ is $\sigma$-stable and $G$ is random totally bounded, then $\sigma$(G) is random totally bounded.
		\item If $G$ is stable and contained in a $\sigma$-stable subset of $E$, then $\bar{G}_{\varepsilon,\lambda}=\overline{(\sigma(G))}_c$.
		\item If $G$ is random sequentially compact, then $G$ is both $\mathcal{T}_{\varepsilon,\lambda}$--complete and $\mathcal{T}_c$--complete.
	\end{enumerate}
\end{lemma}

\begin{proof}
	(1). We can assume, without loss of generality, that $\theta\in G$. For any given $\{A_n,n\in \mathbb{N}\}$ in $\Pi_{\mathcal{F}}$ and any given sequence $\{x_n,n\in \mathbb{N}\}$ in $G$, by the hypothesis on $G$ there exists $x\in E$ such that $\tilde{I}_{A_n} x= \tilde{I}_{A_n} x_n$ for each $n\in \mathbb{N}$. Since $G$ is stable, $p_n:= \sum^n_{k=1} \tilde{I}_{A_k} x_k= \sum^n_{k=1} \tilde{I}_{A_k} x_k+ \tilde{I}_{\bigcup^{\infty}_{i=k+1}I_{A_i}} \cdot \theta \in G$, then $\|x-p_n\|= \|\sum^{\infty}_{i=k+1} \tilde{I}_{A_i} x_i\|= \|\sum^{\infty}_{i=k+1} \tilde{I}_{A_i} x\|= (\sum^{\infty}_{i=k+1} \tilde{I}_{A_i})\|x\| $ converges in probability to 0 since $\sum^{\infty}_{i=k+1} P(A_i)\rightarrow 0$ (as $k\rightarrow +\infty$). So $x\in G$.
	\par
	(2). Denote $z_k=x_{n_k}$ for each $k\in \mathbb{N}$. To prove $\{z_k,k\in \mathbb{N}\}$ converges in $\mathcal{T}_{\varepsilon,\lambda}$ to $x_0$, namely $\{\|z_k-x_0\|, k\in \mathbb{N}\}$ converges in probability to 0, it suffices to prove that any subsequence $\{z_{k_l}: l\in \mathbb{N}\}$ of $\{z_k, k\in \mathbb{N}\}$ has a subsequence $\{z_{k_{l_m}}: m\in \mathbb{N}\}$ such that $\{\|z_{k_{l_m}}-x_0\|, m\in \mathbb{N}\}$ converges a.s. to 0. Since $\{x_{k_l}: l\in \mathbb{N}\}$ still converges in $\mathcal{T}_{\varepsilon,\lambda}$ to $x_0$, there exists a subsequence $\{x_{k_{l_m}}: m\in \mathbb{N}\}$ such that $\{\|x_{k_{l_m}} -x_0\|, m\in \mathbb{N}\}$ converges a.s. to 0. We can assume, without loss of generality, that the sequence $\{n_k, k\in \mathbb{N}\}$ of random indexes for $\{x_{n_k}, k\in \mathbb{N}\}$ satisfies $n_1(\omega)>1$ for each $\omega\in \Omega$. Then $n_k(\omega)> k$ for each $k\in \mathbb{N}$ and each $\omega\in \Omega,$ and hence $\{z_{k_{l_m}}: m\in \mathbb{N}\}$ also converges a.s. to $x_0$.
	
	\par
	(3). Assume that $G$ is random totally bounded. For any given $\varepsilon\in L^0_{++}(\mathcal{F})$, there exist $\{A_n,n \in \mathbb{N}\}$ in $\Pi_{\mathcal{F}}$ and a sequence $\{G_n,n \in \mathbb{N}\}$ of nonempty finite subsets of $G$ such that $\tilde{I}_{A_n}G\subset \tilde{I}_{A_n} (\sigma(G_n))_{(\frac{\varepsilon}{2})}$ for each $n\in \mathbb{N}$. Again since $\bar{G}_c\subset G_{(\frac{\varepsilon}{2})}$, $\tilde{I}_{A_n} \bar{G}_c \subset \tilde{I}_{A_n} G_{(\frac{\varepsilon}{2})} \subset \tilde{I}_{A_n} (\sigma(G_n))_{(\varepsilon)}$, $\bar{G}_c$ is random totally bounded since each $G_n\subset G\subset \bar{G}_c$. Conversely, if $\bar{G}_c$ is random totally bounded, for any given $\varepsilon\in L^0_{++}(\mathcal{F})$, there exist $\{A_n,n \in \mathbb{N}\}$ in $\Pi_{\mathcal{F}}$ and a sequence $\{H_n,n \in \mathbb{N}\}$ of nonempty finite subsets of $\bar{G}_c$ such that $\tilde{I}_{A_n} \bar{G}_c\subset \tilde{I}_{A_n} (\sigma(H_n))_{(\frac{\varepsilon}{2})}$ for each $n\in \mathbb{N}$. Let $H_n=\{\bar{x}_{ni}: i=1,2,\cdots, m_n\}$ for each $n\in \mathbb{N}$. We can choose $x_{ni}\in G$ such that $\|x_{ni}- \bar{x}_{ni}\|< \frac{\varepsilon}{2}$ on $\Omega$ since each $\bar{x}_{ni}\in \bar{G}_c$. Since $\sigma(H_n)= \{\sum_{i=1}^{m_n} \tilde{I}_{A_i} \bar{x}_{ni}: \{A_i, i=1\sim m_n\}$ is an $m_n$--partition of $\Omega$ to $\mathcal{F}\}$ by noticing that $H_n$ consists of $m_n$ elements, it is obvious that $(\sigma(H_n))_{(\frac{\varepsilon}{2})}\subset (\sigma(G_n))_{(\varepsilon)}$, where $G_n=\{x_{ni}: i=1\sim m_n\}$, namely, $\tilde{I}_{A_n}G\subset \tilde{I}_{A_n}\bar{G}_c\subset \tilde{I}_{A_n} (\sigma(H_n))_{(\frac{\varepsilon}{2})}\subset \tilde{I}_{A_n} (\sigma(G_n))_{(\varepsilon)},$ which has shown that $G$ is random totally bounded.
	\par
	(4). Let $\varepsilon,G,\{G_n,n\in \mathbb{N}\}$ and $\{A_n,n\in \mathbb{N}\}$ be as in Definition \ref{definition2.1} such that $\tilde{I}_{A_n} G \subset \tilde{I}_{A_n} (\sigma(G_n))_{(\varepsilon)}$ for each $n\in \mathbb{N}$. Since $\sigma(\tilde{I}_{A_n} H)= \tilde{I}_{A_n} \sigma(H)$ for any subset $H$ of $E$ and $(\sigma(G_n))_{(\varepsilon)}$ is $\sigma$--stable, $\tilde{I}_{A_n} \sigma(G)=\sigma(\tilde{I}_{A_n}G)\subset \sigma(\tilde{I}_{A_n} (\sigma(G_n))_{(\varepsilon)})= \tilde{I}_{A_n} \sigma(\sigma(G_n)_{(\varepsilon)})= \tilde{I}_{A_n} (\sigma(G_n))_{(\varepsilon)}$, so $\sigma(G)$ is random totally bounded by observing $G_n\subset G \subset \sigma(G)$.
	\par
	(5). We can assume, without loss of generality, that $\theta\in G$. Let $x$ be any element in $\sigma(G)$ and $x=\sum^{\infty}_{n=1} \tilde{I}_{A_n} x_n$ for some $\{A_n, n\in \mathbb{N}\}$ in $\Pi_{\mathcal{F}}$ and some sequence $\{x_n, n\in \mathbb{N}\}$ in $G$. Then $p_n: = \sum ^n_{k=1} \tilde{I}_{A_k} x_k= \sum ^n_{k=1} \tilde{I}_{A_k} x_k+ I_{\bigcup^{\infty}_{k=n+1}A_k}\theta\in G$ for each $n\in \mathbb{N}$ and it is clear that $\{p_n,n\in \mathbb{N}\}$ converges in $\mathcal{T}_{\varepsilon,\lambda}$ to $x$, so $x\in \bar{G}_{\varepsilon,\lambda}$. Thus, $\sigma(G)\subset \bar{G}_{\varepsilon,\lambda}$, namely $\bar{G}_{\varepsilon,\lambda}= \overline {(\sigma(G))}_{\varepsilon,\lambda}$, and further one can have $\bar{G}_{\varepsilon,\lambda}= \overline {(\sigma(G))}_c$ by Theorem 3.12 of \cite{Guo3}.
	\par
	(6). Let $\{x_n,n\in \mathbb{N}\}$ be a $\mathcal{T}_{\varepsilon,\lambda}$--Cauchy sequence in $G$. We only need to prove that $\{x_n,n\in \mathbb{N}\}$ has a subsequence convergent in $G$. It is obvious that $\{x_n,n\in \mathbb{N}\}$ has a subsequence $\{x_{n_k},k\in \mathbb{N}\}$ such that $\normm{x_{n_{k+1}}- x_{n_k}}< \frac{1}{2^k}$ for each $k\in \mathbb{N}$, where $\normm{\cdot}$ is the quasinorm on $E$ defined by $\normm{x}=\int_{\Omega} \frac{\|x\|}{1+\|x\|} dP$ for each $x\in E$. We can assume, without loss of generality, that $\{x_{n_k},k\in \mathbb{N}\}$ is just $\{x_n,n\in \mathbb{N}\}$ itself, namely $\normm{x_{k+1}- x_k}< \frac{1}{2^k}$ for each $k\in \mathbb{N}$. Since $G$ is random sequentially compact, there exists a random subsequence $\{x_{n_k},k\in \mathbb{N}\}$ of $\{x_n,n\in \mathbb{N}\}$ such that $\{x_{n_k},k\in \mathbb{N}\}$ converges in $\mathcal{T}_{\varepsilon,\lambda}$ to an element $x_0\in G$. Once more, we can assume that $n_k(\omega)> k$ for each $\omega\in \Omega$ and each $k\in \mathbb{N}$. Since $\{x_{n_k},k\in \mathbb{N}\}$  is a $\mathcal{T}_{\varepsilon,\lambda}$--Cauchy sequence, for any given positive number $\varepsilon$ there exists some $k_0\in \mathbb{N}$ such that $\normm{x_{n_k}-x_{n_{k_0}}}< \frac{\varepsilon}{2}$ for any $k>k_0$. Further, let $M$ be a sufficiently large positive integer such that $M>k_0$ and $P(n_{k_0}\geq M)< \frac{\varepsilon}{4}$. Then for any given $k\in \mathbb{N}$,
	\begin{align*}
		&\normm{ \sum^{\infty}_{l=M} \tilde{I}_{(n_{k_0}=l)} (x_k-x_l) } \\
		&=\int_{\Omega} \frac{\|\sum^{\infty}_{l=M} \tilde{I}_{(n_{k_0}=l)} (x_k-x_l)\|}{1+\|\sum^{\infty}_{l=M} \tilde{I}_{(n_{k_0}=l)} (x_k-x_l)\|}dP\\
		&\leq \sum^{\infty}_{l=M} \int_{\Omega} \frac{\tilde{I}_{(n_{k_0}=l)}\|x_k-x_l\| }{1+ \tilde{I}_{(n_{k_0}=l)}\|x_k-x_l\|}dP\\
		&=\sum^{\infty}_{l=M} \int_{(n_{k_0}=l)} \frac{\tilde{I}_{(n_{k_0}=l)}\|x_k-x_l\| }{1+ \tilde{I}_{(n_{k_0}=l)}\|x_k-x_l\|}dP\\
		&\leq \sum^{\infty}_{l=M} P(n_{k_0}=l)\\
		&=P(n_{k_0}\geq M)\\
		&<\frac{\varepsilon}{4}.
	\end{align*}
	\par
	We can assume, without loss of generality, that $k_0$ is also so large that $\frac{1}{2^{k_0-1}}< \frac{\varepsilon}{4}$. Then when $k>M$ (so $k>k_0$ since $M>k_0$), we have the following:
	\begin{align*}
		\normm{ x_k-x_{n_{k_0}} }&=\normm{x_k- \sum^{\infty}_{l=k_0+1} \tilde{I}_{(n_{k_0}=l) }x_l}\\
		&=\normm{\sum^{\infty}_{l=k_0+1} \tilde{I}_{(n_{k_0}=l) }(x_k-x_l)}~~~(\text{note}: \sum^{\infty}_{l=k_0+1} \tilde{I}_{(n_{k_0}=l) }=1)\\
		&=\normm{\sum^{M-1}_{l=k_0+1} \tilde{I}_{(n_{k_0}=l) }(x_k-x_l)}+ \normm{\sum^{\infty}_{l=M} \tilde{I}_{(n_{k_0}=l) }(x_k-x_l)}\\
		&\leq \sum^{M-1}_{l=k_0+1}\normm{x_k-x_l}+ \frac{\varepsilon}{4}\\
		&<\sum^{M-1}_{l=k_0+1} \frac{1}{2^{l-1}}+ \frac{\varepsilon}{4}  \\
		&< \frac{1}{2^{k_0-1}}+ \frac{\varepsilon}{4}\\
		&<\frac{\varepsilon}{4}+ \frac{\varepsilon}{4}\\
		&=\frac{\varepsilon}{2}.
	\end{align*}
	\par
	To sum up, when $k>M$ , we can prove that $\normm{x_k-x_{n_k}}\leq \normm{x_k-x_{n_{k_0}}}+ \normm{x_{x_{k_0}}- x_{n_k}}< \varepsilon$, which shows that $lim_{k\rightarrow \infty} \normm{x_k-x_{n_k}}=0$. So $\{x_k,k\in \mathbb{N}\}$ must converge in $\mathcal{T}_{\varepsilon,\lambda}$ to $x_0$. Thus $G$ is $\mathcal{T}_{\varepsilon,\lambda}$--complete. Further, it immediately follows from Remark \ref{remark1.7} that $G$ is also $\mathcal{T}_c$--complete.
\end{proof}

We are now ready to prove Theorem \ref{theorem2.3} as follows.

\begin{proof}[\textbf{Proof of Theorem \ref{theorem2.3}}]
\textbf{Necessity}. We only need to prove that $G$ is random total bounded since it must be $\mathcal{T}_{\varepsilon,\lambda}$-complete from (6) of Lemma \ref{lemma2.5}. Given an $\varepsilon\in L^0_{++}(\mathcal{F})$, we will look for an element $\{A_n,n\in \mathbb{N}\}$ in $\Pi_{\mathcal{F}}$ and a sequence $\{G_n,n\in \mathbb{N}\}$ of nonempty finite subsets of $G$ such that $\tilde{I}_{A_n}G\subset \tilde{I}_{A_n}(\sigma(G_n))_{(\varepsilon)}$ for each $n\in \mathbb{N}$, namely, $\{A_n,n\in \mathbb{N}\}$ and $\{G_n,n\in \mathbb{N}\}$ satisfy Definition \ref{definition2.1}. We first look for $G_1$ and $A_1$ as follows.
\par
For this, let us first arbitrarily choose an element $x_1\in G$ and define $f_1:G\rightarrow L^0(\mathcal{F})$ by $f_1(x)=\|x-x_1\|$ for any $x\in G$. Then $f_1$ is $\mathcal{T}_{\varepsilon,\lambda}$--continuous (and, naturally, also random sequentially continuous). Denote $\bigvee f_1(G)$ by $\xi_1$ and let $B$ be an arbitrarily chosen representative of $(\xi_1<\varepsilon)$. If $P(B)>0$, then we take $G_1=\{x_1\}$ and $A_1=B$. Otherwise, $\xi_1=\bigvee f_1(G)\geq \varepsilon$. By (4) of Lemma \ref{lemma2.4} there exists $x_2\in G$ such that $f_1(x_2)=\|x_2-x_1\|=\xi_1\geq \varepsilon$. Further we define $f_2:G\rightarrow L^0(\mathcal{F})$ by $f_2(x)=\|x-x_1\|\bigwedge \|x-x_2\|$ for any $x\in G$. Denote $\bigvee f_2(G)$ by $\xi_2$ and let $B_1$ be an arbitrarily chosen representative of $(\xi_2<\varepsilon)$. If $P(B_1)>0,$ then we take $G_1=\{x_1,x_2\}$ and $A_1=B_1$. Since $f_2(x)=d(x,G_1)=d(x,\sigma (G_1))$ for any $x\in G$, it is obvious that $\tilde{I}_{A_1}G\subset \tilde{I}_{A_1}(\sigma(G_1))_{(\varepsilon)}$. If $P(B_1)=0,$ then $\xi_2=\bigvee f_2(G)\geq \varepsilon$. Again by (4) of Lemma \ref{lemma2.4} there exists $x_3\in G$ such that $\|x_3-x_1\|\bigwedge\|x_3-x_2\|\geq \varepsilon$. We can continue to look for $G_1$ and $A_1$ as above. Then by the induction method we can certainly obtain, after some finite steps, for example, in $n_1$ steps, a finite subset $G_1=\{x_1,x_2,\cdots,x_{n_1}\}$ of $G$ such that $P(\xi_{n_1}<\varepsilon)>0$, where $\xi_{n_1}=\bigvee\{d(x,G_1):x\in G\}$. Then by taking $A_1=$ an arbitrarily chosen representative of $(\xi_{n_1}<\varepsilon)$ we have $\tilde{I}_{A_1}G\subset \tilde{I}_{A_1}(\sigma(G_1))_{(\varepsilon)}$ by (1) and (2) of Lemma \ref{lemma2.4}. Otherwise, if, we can not find $G_1$ and $A_1$ after any finite steps, then there exists a sequence $\{x_n, n\in \mathbb{N}\}$ in $G$ such that $\bigwedge\{\|x_{k+1}-x_i\|: 1\leq i\leq k\}\geq \varepsilon$ for each $k\in \mathbb{N}$ so $\{x_n, n\in \mathbb{N}\}$ does not admit any random subsequence convergent in $\mathcal{T}_{\varepsilon,\lambda}$, which is a contradiction to the fact that $G$ is random sequentially compact.
\par
Once  we have found $G_1$ and $A_1$, if $P(A_1)=1$, we can take $A_1=\Omega$. Then $\{A_1\}$ has formed an element in $\Pi_{\mathcal{F}}$, at which time $G$ is automatically random totally bounded. If $P(A_1)<1,$ we continue to look for $G_2$ and $A_2$ as follows.
\par
Let $\Omega'=\Omega\backslash A_1,\mathcal{F}'=\Omega'\bigcap\mathcal{F}$ and $P': \mathcal{F}'\rightarrow [0,1]$ be defined by $P'(A)=P(A)/ P(\Omega')$ for each $A\in \mathcal{F}'$. Then $(\Omega', \mathcal{F}',P')$ is a probability space and $(E',\|\cdot\|')$ can be regarded as an $RN$ module with base $(\Omega', \mathcal{F}',P')$, where $E'=\tilde{I}_{\Omega'}E:= \{\tilde{I}_{\Omega'} x: x\in E\}$ and $\|\cdot\|'$ is the limitation of $\|\cdot\|$ to $E'$. Further, it is also obvious that $\tilde{I}_{\Omega'}G:=\{\tilde{I}_{\Omega'} x: x\in G\}$ is a random sequentially compact $\sigma$--stable subset of $E'$. Since $\xi_{n_1}=\bigvee \{d(x,G_1): x\in G\}\geq \varepsilon$ on $\Omega'$, namely $\bigvee\{d(x,G_1'):x\in \tilde{I}_{\Omega'}G\}\geq \tilde{I}_{\Omega'}\varepsilon$, where $G_1'=\tilde{I}_{\Omega'}G_1$. Then there exists $x_{n_1+1}\in G$ such that $d(x_{n_1+1},G_1)\geq \varepsilon$ on $\Omega'$. As in the process of constructing $G_1$ and $A_1$, we can add some finite subset $\{x_{n_1+1},\cdots,x_{n_2}\}$ of $G$ to $G_1$ so that we can obtain $G_2=\{x_1,x_2,\cdots,x_{n_1},x_{n_1+1},\cdots,x_{n_2}\}$ and some $A_2\in \mathcal{F}$ such that $A_2\subset \Omega', P(A_2)>0$ and $\bigvee\{d(x,G_2): x\in G\}<\varepsilon$ on $A_2$. In fact, $A_2$ is exactly $\Omega'\bigcap B$, where $B$ is a chosen representative of $(\bigvee\{d(x,G_2): x\in G\}<\varepsilon)$, namely, $G_2$ and $A_2$ are just desired.
\par
Once we have found $G_2$ and $A_2$, if $P(A_1\bigcup A_2)=1$, we have finished the proof of Necessity since $\{A_1,A_2\}$ forms an element in $\Pi_{\mathcal{F}}$ and $G$ is automatically random totally bounded. If $P(A_1\bigcup A_2)<1,$ we will continue to construct $G_3$ and $A_3,$ and necessarily for $G_4$ and $A_4$, if necessary, etc. By the induction method we can obtain the desired sequence $\{G_n,n\in \mathbb{N}\}$ and $\{A_n,n\in \mathbb{N}\}$ in $\Pi_{\mathcal{F}}$ after at most countably many steps, which completes the proof of Necessity.
\par
\textbf{Sufficiency}. Assume that $G$ is random totally bounded and $\mathcal{T}_{\varepsilon,\lambda}$--complete. For $\varepsilon=\frac{1}{2}$, let $\{G^{(1)}_n,n\in \mathbb{N}\}$ and $\{A^{(1)}_n,n\in \mathbb{N}\}$ be as in Definition \ref{definition2.1} such that $\tilde{I}_{A^{(1)}_n}G\subset \tilde{I}_{A^{(1)}_n}(\sigma(G^{(1)}_n))_{(\frac{1}{2})}$ for each $n\in \mathbb{N}$. We can assume, without loss of generality, that $P(A^{(1)}_n)>0$ for each $n\in \mathbb{N}$, and further we put $G^{(1)}_n=\{x^{(1)}_{ni}: i=1,2,\cdots, m^{(1)}_n\}$ for each $n\in \mathbb{N}$. Let $\{y_k,k\in \mathbb{N}\}$ be a given sequence in $G$, we will construct a random subsequence $\{y_{k^{(1)}_p},p\in \mathbb{N}\}$ in $G$ of $\{y_k,k\in \mathbb{N}\}$ such that the random diameter of $\{y_{k^{(1)}_p},p\in \mathbb{N}\}$, namely, $\bigvee\{\|y_{k^{(1)}_p}-y_{k^{(1)}_q}\|: p,q\in \mathbb{N}\}\leq 1$ by first constructing $\{k^{(1)}_p, p\in \mathbb{N}\}$ on each $A^{(1)}_n$.
\par
First, we consider this construction on $A^{(1)}_1$. Since $\tilde{I}_{A^{(1)}_1}\{y_k,k\in \mathbb{N}\}\subset \tilde{I}_{A^{(1)}_1}G\subset \tilde{I}_{A^{(1)}_1}(\sigma(G^{(1)}_1))_{(\frac{1}{2})}$, then by (1) and (2) of Lemma \ref{lemma2.4} $d(y_k, \sigma(G_1^{(1)}))=d(y_k, G_1^{(1)})= \bigwedge\{\|y_k-x^{(1)}_{1i}\|: 1\leq i\leq m_1^{(1)}\}< \frac{1}{2}$ on $A_1^{(1)}$ for each $k\in \mathbb{N}$. Since $\mathbb{N}\times \{i: 1\leq i\leq m_1^{(1)}\}$ is countable, we can choose a representative $\|y_k-x^{(1)}_{1i}\|^0$ for each  $\|y_k-x^{(1)}_{1i}\|$ such that $min\{\|y_k-x^{(1)}_{1i}\|^0(\omega): 1\leq i\leq m_1^{(1)}\}< \frac{1}{2}$ for each $\omega\in A_1^{(1)}$ and each $k\in \mathbb{N}$. For each $i\in \{1,2,\cdots, m_1^{(1)}\}$, let $B_i=\{\omega\in A_1^{(1)}: \sum^{\infty}_{k=1}I_{B_{ki}}(\omega)=+\infty\}$, where $B_{ki}=\{\omega\in \Omega: \|y_k-x^{(1)}_{1i}\|^0(\omega)<\frac{1}{2}\}$, then it is easy to check that $\bigcup_{i=1}^{m_1^{(1)}}B_i=A_1^{(1)}$ since for each $\omega\in A_1^{(1)}$ there always exists some $i$ such that $\|y_k-x^{(1)}_{1i}\|^0(\omega)<\frac{1}{2}$ for infinitely many $k^,$s. Now, let $D_1=B_1$ and $D_i=B_i\backslash(\bigcup^{i-1}_{l=1}B_l)$ for each $i\in \{2,3,\cdots, m_1^{(1)}\}$, then $\{D_i: 1\leq i\leq m_1^{(1)}\}$ forms a finite partition of $A_1^{(1)}$ to $\mathcal{F}$. Further we define $l_1^{(1)}: \Omega\rightarrow \{1,2,\cdots,m_1^{(1)}\}$ by $l_1^{(1)}(\omega)=1$ if $\omega\in D_1\bigcup (A_1^{(1)})^c$ and $l_1^{(1)}(\omega)=i$ if $\omega\in D_i$ for some $i\in \{2,3,\cdots,m_1^{(1)}\}$. Then $l_1^{(1)}$ is a random variable with values in $\{1,2,\cdots,m_1^{(1)}\}$ and $x^1_{1l_1^{(1)}}:= \sum_{i=1}^{m_1^{(1)}}\tilde{I}_{(l_1^{(1)}=i)}x^{(1)}_{1i}\in G$. Since $\|y_k-x^{(1)}_{1l_1^{(1)}}\|=\sum_{i=1}^{m_1^{(1)}}\tilde{I}_{(l_1^{(1)}=i)}\|y_k-x^{(1)}_{1i}\|$ for each $k\in \mathbb{N}$ and $\mathbb{N}\times \{1,2,\cdots, m_1^{(1)}\}$ is countable, we can choose a representative $\|y_k-x^{(1)}_{1l_1^{(1)}}\|^0$ of $\|y_k-x^{(1)}_{1l_1^{(1)}}\|$ for each $k\in \mathbb{N}$ such that for each $\omega\in A_1^{(1)}$ there exist infinitely many $k^,$s such that $\|y_k-x^{(1)}_{1l_1^{(1)}}\|^0(\omega)<\frac{1}{2}$.
\par
Similarly, for each $n\geq 2$  there exist an element $x^{(1)}_{nl_n^{(1)}}$ in $G$ and a representative $\|y_k-x^{(1)}_{nl_n^{(1)}}\|^0$ of $\|y_k-x^{(1)}_{nl_n^{(1)}}\|$ for each $k\in \mathbb{N}$ such that for each $\omega\in A_n^{(1)}$ there always exist infinitely many $k^,$s such that $\|y_k-x^{(1)}_{nl_n^{(1)}}\|^0(\omega)<\frac{1}{2}$. Now, let $k^1_0(\omega)\equiv 1$, we can recursively define a sequence $\{k_p^{(1)}, p\in \mathbb{N}\}$ of random variables with values in $\mathbb{N}$ as follows:
\par
$k_p^{(1)}(\omega)=min\{k\in \mathbb{N}: k>k^{(1)}_{p-1}(\omega)$ and $\|y_k-x^{(1)}_{nl_n^{(1)}}\|^0(\omega)<\frac{1}{2}\}$ when $\omega\in A_n^{(1)}$ for some $n\in \mathbb{N}$.
\par
For any $n$ and $m$ in $\mathbb{N}$, let $B_{nm}=[\{\omega\in \Omega: k^{(1)}_{p-1}(\omega)<m$ and $\|y_m-x^{(1)}_{nl_n^{(1)}}\|^0(\omega)<\frac{1}{2}\}\backslash (\bigcup^{m-1}_{r=1}\{\omega\in\Omega: k^{(1)}_{p-1}(\omega)<r $ and $\|y_r-x^{(1)}_{nl_n^{(1)}}\|^0(\omega)<\frac{1}{2}\})]\bigcap A_n^{(1)}$, then by the property of $x^{(1)}_{nl_n^{(1)}}$ as above it is obvious that $A_n^{(1)}=\bigcup ^{\infty}_{m=1}B_{nm}$ for each $n\in \mathbb{N}$. Noting that $\{B_{nm}: m\in \mathbb{N}\}$ for each fixed $n\in \mathbb{N}$ is a disjoint family, one has the following:
\par
\begin{equation}\label{equation2.1}
	\tilde{I}_{A_n^{(1)}}=\sum_{m=1}^{\infty}\tilde{I}_{B_{nm}} ~~{\text {for~ each}} ~n\in \mathbb{N}.
\end{equation}
\par
Since $\|y_m-x^{(1)}_{nl_n^{(1)}}\|^0(\omega)<\frac{1}{2}$ for each $\omega\in B_{nm}$, one has the following:
\begin{equation}\label{equation2.2}
	\tilde{I}_{B_{nm}}\|y_m-x^{(1)}_{nl_n^{(1)}}\|<\frac{1}{2}~\text{ on}~B_{nm}~\text{for~any}~n~\text{and}~m~\text{in}~\mathbb{N}.
\end{equation}
\par
Since $\{\omega\in \Omega:k_p^{(1)}(\omega)=m\}=\bigcup^{\infty}_{n=1}B_{nm}$ for each $m\in \mathbb{N}$ and $\{B_{nm}: n\in \mathbb{N}\}$ is also a disjoint family for each fixed $m\in \mathbb{N}$, one has the following:
\begin{equation}\label{equation2.3}
	\tilde{I}_{(k_p^{(1)}=m)}=\sum_{n=1}^{\infty}\tilde{I}_{B_{nm}}~\text{for~each}~m\in \mathbb{N}.
\end{equation}

\par
Now, let $x^{(1)}=\sum_{n=1}^{\infty}\tilde{I}_{A_n^{(1)}}x^{(1)}_{nl_n^{(1)}}$. Then $x^{(1)}\in G$ and

\begin{align}\label{equation2.4}
&\|y_{k_p^{(1)}}-x^{(1)}\|\nonumber\\
	&=\|\sum_{m=1}^{\infty}\tilde{I}_{(k_p^{(1)}=m)}y_m-\sum_{n=1}^{\infty}\tilde{I}_{A_n^{(1)}}x^{(1)}_{nl_n^{(1)}}\|\nonumber\\
	&=\|\sum_{m=1}^{\infty} \sum_{n=1}^{\infty} \tilde{I}_{B_{nm}}y_m- \sum_{n=1}^{\infty} \sum_{m=1}^{\infty} \tilde{I}_{B_{nm}} x^{(1)}_{nl_n^{(1)}}\|~~~\text{(by (\ref{equation2.1}) and (\ref{equation2.3}))} \\
	&=\sum_{n=1}^{\infty} \sum_{m=1}^{\infty} \tilde{I}_{B_{nm}}\|y_m-x^{(1)}_{nl_n^{(1)}}\|\nonumber\\
	&<(\sum_{n=1}^{\infty} \sum_{m=1}^{\infty} \tilde{I}_{B_{nm}})\cdot \frac{1}{2}~\text{on}~\Omega~~~\text{(by (\ref{equation2.2}))}\nonumber\\
	&=\frac{1}{2}~\text{on}~\Omega~~~(\text{note:}~\sum_{n=1}^{\infty} \sum_{m=1}^{\infty}B_{nm}= \sum_{n=1}^{\infty} A_n^{(1)}=\Omega).\nonumber
\end{align}

\par
Thus
\begin{equation}\label{equation2.5}
	\bigvee\{\|y_{k_p^{(1)}}- y_{k_q^{(1)}}\|: p,q\in \mathbb{N}\}\leq 1.
\end{equation}
\par
For $\varepsilon=\frac{1}{4}$, we next construct a random subsequence $\{y_{k_p^{(2)}}, p\in \mathbb{N}\}$ of $\{y_k, k\in \mathbb{N}\}$ such that the following two properties are satisfied:
\begin{enumerate}[$(i)$]
	\item $\bigvee\{\|y_{k_p^{(2)}}- y_{k_q^{(2)}}\|: p,q\in \mathbb{N}\}\leq \frac{1}{2}$.
	\item $\bigvee\{\|y_{k_p^{(1)}}- y_{k_q^{(2)}}\|: p,q\in \mathbb{N}\}\leq 1$.
\end{enumerate}
\par
Let $\{A_n^{(2)}, n\in \mathbb{N}\}$ in $\Pi_{\mathcal{F}}$ and $\{G_n^{(2)}, n\in \mathbb{N}\}$ satisfy Definition \ref{definition2.1} for $\varepsilon=\frac{1}{4}$, namely, $\tilde{I}_{A_n^{(2)}}G\subset \tilde{I}_{A_n^{(2)}} (\sigma(G_n^{(2)}))_{(\frac{1}{4})}$ for each $n\in \mathbb{N}$. Let $G_n^{(2)}= \{x^{(2)}_{ni}: i=1,2,\cdots, m^{(2)}_n\}$ for each $n\in \mathbb{N}$. By replacing $\{y_k, k\in \mathbb{N}\}$ with $\{y_{k_p^{(1)}}, p\in \mathbb{N}\}$ and similarly defining $l_n^{(2)}, x^{(2)}_{nl_n^{(2)}}$ and $\|y_{k_p^{(1)}}-x^{(2)}_{nl_n^{(2)}}\|^0$ as in the process of defining $l_n^{(1)}, x^{(1)}_{nl_n^{(1)}}$ and $\|y_k-x^{(1)}_{nl_n^{(1)}}\|^0$, we can have that for each $\omega\in A_n^{(2)}$ there exist infinitely many $p^,$s such that $\|y_{k_p^{(1)}}-x^{(2)}_{nl_n^{(2)}}\|^0(\omega)<\frac{1}{4}$.
\par
Now, let $s_0(\omega)\equiv 1$, again we can recursively define a sequence $\{s_p,p\in \mathbb{N}\}$ of random variables with values in $\mathbb{N}$ as follows:
\par
$s_p(\omega)= min\{r\in \mathbb{N}: r>s_{p-1}(\omega)$ and $\|y_{k_r^{(1)}}- x^{(2)}_{nl_n^{(2)}}\|^0(\omega)< \frac{1}{4}\}$ when $\omega\in A_n^{(2)}$ for some $n\in \mathbb{N}$.
\par
Let $x^{(2)}=\sum^{\infty}_{n=1}\tilde{I}_{A_n^{(2)}} x^{(2)}_{nl_n^{(2)}}$  and $k_p^{(2)}= k^{(1)}_{s_p}$ for each $p\in \mathbb{N}$, where $k^{(1)}_{s_p}(\omega)= k^{(1)}_{s_p(\omega)}(\omega)$ for each $\omega\in \Omega$, then $x^{(2)}\in G$. We can similarly prove the following as in the proofs of (\ref{equation2.4}) and (\ref{equation2.5}):
\begin{equation}\label{equation2.6}
	\|y_{k_p^{(2)}}- x^{(2)}\|<\frac{1}{4}~~\text{on}~\Omega.
\end{equation}
\begin{equation}\label{equation2.7}
	\bigvee\{\|y_{k_p^{(2)}}-y_{k_q^{(2)}}\|: p,q\in \mathbb{N} \}\leq \frac{1}{2},
\end{equation}
which is just our desired result (i).
\par
As for (ii), for any $p$ and $q\in \mathbb{N}$, since $k_q^{(2)}= k^{(1)}_{s_q}$ and
\begin{align*}
	y_{k_q^{(2)}} & = y_{k^{(1)}_{s_q}} \\
	& = \sum^{\infty}_{m=1} \tilde{I}_{(s_q=m)} y_{k_m^{(1)}},
\end{align*}
we have
\begin{align*}
	\|y_{k_p^{(1)}}- y_{k_q^{(2)}}\| & = \|y_{k_p^{(1)}}- \sum^{\infty}_{m=1} \tilde{I}_{(s_q=m)} y_{k_m^{(1)}}\| \\
	& = \sum^{\infty}_{m=1} \tilde{I}_{(s_q=m)} \|y_{k_p^{(1)}}- y_{k_m^{(1)}}\| \\
	& \leq 1,
\end{align*}
which gives (ii).
\par
We can inductively define a random subsequence $\{y_{k_p^{(n)}}, p\in \mathbb{N}\}$ of $\{y_k, k\in \mathbb{N}\}$ for each $n\in \mathbb{N}$ with the following two properties:

(iii)~~~$\bigvee\{\|y_{k_p^{(n)}}- y_{k_q^{(n)}}\|: p,q\in \mathbb{N}\}\leq \frac{1}{n}$ for each $n\in \mathbb{N}$.

(iv)~~~$\bigvee\{\|y_{k_p^{(n-1)}}- y_{k_q^{(n)}}\|: p,q\in \mathbb{N}\}\leq \frac{1}{n-1}$ for each $n\geq 2$.
\par
Finally, as usual, by taking a diagonal subsequence, let $k_p=k_p^{(p)}$ for each $p\in \mathbb{N}$. Then $\{y_{k_p}: p\in \mathbb{N}\}$ is still a random subsequence of $\{y_k, k\in \mathbb{N}\}$, which is obviously a $\mathcal{T}_{\varepsilon,\lambda}$--Cauchy sequence in $G$ by (iv) and thus converges in $\mathcal{T}_{\varepsilon,\lambda}$ to an element $y_0$ in $G$ by the $\mathcal{T}_{\varepsilon,\lambda}$--completeness of $G$.
\end{proof}

\par
From  the proof of the sufficiency part of Theorem \ref{theorem2.3}, one can see that a $\sigma$--stable set $G$ is random sequentially compact iff both $G$ is $\mathcal{T}_{\varepsilon,\lambda}$--complete and for any positive number $\varepsilon$ there exist $\{A_n, n\in \mathbb{N} \}$ in $\Pi_{\mathcal{F}}$ and a sequence $\{G_n, n\in \mathbb{N}\}$ of nonempty finite subsets of $G$ such that $\tilde{I}_{A_n} G \subset \tilde{I}_{A_n} (\sigma(G_n))_{(\varepsilon)}$ for each $n\in \mathbb{N}$. That is to say, to verify that a $\sigma$--stable $\mathcal{T}_{\varepsilon,\lambda}$--complete set $G$ is random totally bounded, it suffices to check that for any positive $\varepsilon$ there exist $\{A_n, n\in \mathbb{N} \}$ in $\Pi_{\mathcal{F}}$ and a sequence $\{G_n, n\in \mathbb{N}\}$ of nonempty finite subsets of $G$ such that Definition \ref{definition2.1} is satisfied. Such a fact will be used in the proof of Corollary \ref{corollary2.6} below. Since a metric space can be isometrically embedded into a normed space, Corollary \ref{corollary2.6} also shows that $L^0(\mathcal{F}, M)$ is random sequentially compact when $(M,d)$ is a compact metric space.
\begin{corollary}\label{corollary2.6}
	Let $(X,\|\cdot\|)$ be a normed space and $M$ a nonempty compact set of $X$. Then $L^0(\mathcal{F}, M)$ is a random sequentially compact subset of the $RN$ module $L^0(\mathcal{F}, X)$.
\end{corollary}
\begin{proof}
	It is easy to see that both $L^0(\mathcal{F}, X)$ and $L^0(\mathcal{F}, M)$ are $\sigma$--stable and $L^0(\mathcal{F}, M)$ is also $\mathcal{T}_{\varepsilon,\lambda}$--complete. By Theorem \ref{theorem2.3} we only need to prove that $L^0(\mathcal{F}, M)$ is random totally bounded. For any given positive number $\varepsilon$, there exists a finite subset $\{x_1,x_2,\cdots, x_n\}$ of $M$ such that $M\subset \bigcup^n_{i=1} B(x_i,\varepsilon)$, where $B(x_i, \varepsilon)= \{x\in X: \|x-x_i\|< \varepsilon\}$. Let $p\in L^0(\mathcal{F},M)$ be any given element and arbitrarily choose a representative $p^0$ of $p$. Then $D_i:= (p^0)^{-1}(C_i)= \{\omega\in \Omega: p^0(\omega)\in C_i\}\in \mathcal{F}$, where $C_1=B(x_1, \varepsilon)$ and $C_i=B(x_i, \varepsilon)\backslash (\bigcup_{k=1}^{i-1} B(x_k, \varepsilon))$ for any $i=2\sim n.$ It is obvious that $\|p^0(\omega)-\sum^n_{i=1} I_{D_i}(\omega)x_i\|< \varepsilon$ for each $\omega\in \Omega$, namely, $\|p-\sum^n_{i=1} \tilde{I}_{D_i}x_i\|< \varepsilon$ on $\Omega$. Since $\{D_i, i=1\sim n\}$ forms a finite partition of $\Omega$ to $\mathcal{F}$, $p\in B(\sum^n_{i=1} \tilde{I}_{D_i}x_i, \varepsilon)\subset (\sigma(G_1))_{(\varepsilon)}$, where $G_1=\{x_1,x_2,\cdots,x_n\}$. This shows that $L^0(\mathcal{F}, M)\subset (\sigma(G_1))_{(\varepsilon)}$. We then have $L^0(\mathcal{F}, M)$, of course, random totally bounded by taking $\{A_n, n\in \mathbb{N}\}=\{\Omega\}$ and $\{G_n, n\in \mathbb{N}\}=\{G_1\}$ in Definition \ref{definition2.1}.
\end{proof}
\begin{corollary}\label{corollary2.7}
	Let $(E,\|\cdot\|)$ be a $\mathcal{T}_{\varepsilon,\lambda}$--complete $RN$ module with base $(\Omega,\mathcal{F},P)$ and $G$ a random sequentially compact $\sigma$--stable subset of $E$. Then $\overline{[Conv_{L^0}(G)]}_{\varepsilon,\lambda}$ is random sequentially compact.
\end{corollary}
\begin{proof}
	Since ~\ $Conv_{L^0}(G)$ ~\ is $L^0$--convex, ~\ it is obvious that ~\ $Conv_{L^0}(G)$ ~\ is stable (but not necessarily $\sigma$--stable), which in turn also implies that $\overline{[Conv_{L^0}(G)]}_{\varepsilon,\lambda}= \overline{[\sigma(Conv_{L^0}(G))]}_c $ by (5) of Lemma \ref{lemma2.5}. Further, $\overline{[Conv_{L^0}(G)]}_{\varepsilon,\lambda}$ is also $\sigma$--stable by (1) of Lemma \ref{lemma2.5} since $\overline{[Conv_{L^0}(G)]}_{\varepsilon,\lambda}$ is obviously stable. We only need to prove that $\overline{[\sigma(Conv_{L^0}(G))]}_c$ is random totally bounded. Then by Properties (3) and (4) of Lemma \ref{lemma2.5} we only need to prove that $Conv_{L^0}(G)$ is random totally bounded, as follows.
	\par
	By Theorem \ref{theorem2.3}, $G$ is random totally bounded. Then for any $\varepsilon\in L^0_{++}(\mathcal{F})$ there exist $\{A_n,n\in \mathbb{N}\}$ in $\Pi_{\mathcal{F}}$ and a sequence $\{G_n,n\in \mathbb{N}\}$ of nonempty finite subsets of $G$ such that $\tilde{I}_{A_n} G\subset \tilde{I}_{A_n} (\sigma(G_n))_{(\frac{\varepsilon}{2})}$ for each $n\in \mathbb{N}$. But, since each $G_n$ is finite, for example, $G_n= \{x^{(n)}_{ni}: i=1,2,\cdots, m_n\}$, $\sigma(G_n)= \{\sum_{i=1}^{m_n} \tilde{I}_{D_i} x^{(n)}_{ni}: \{D_i: i=1\sim m_n\}$ is an $m_n$--partition of $\Omega$ to $\mathcal{F} \}\subset Conv_{L^0}(G_n)= \{\sum_{i=1}^{m_n} \xi_i x^{(n)}_{ni}: \xi_i\in L^0_+(\mathcal{F})$ for each $i= 1\sim m_n$ such that $\sum_{i=1}^{m_n}\xi_i=1\}$. For any fixed $n\in \mathbb{N}$, $\{(\xi_1, \xi_2, \cdots, \xi_{m_n}): \sum_{i=1}^{m_n} \xi_i= 1$ and $\xi_i\in L^0_+(\mathcal{F})$ for each $i=1\sim m_n\}$ is, obviously, an a.s. bounded set in $L^0(\mathcal{F}, \mathbb{R}^{m_n})$, and hence also random sequentially compact, from which one can have that $Conv_{L^0}(G_n)$ is random sequentially compact. And, thus, it is also clearly $\mathcal{T}_{\varepsilon,\lambda}$--closed. Further, it is obvious that $Conv_{L^0}(G_n)$ is $\sigma$--stable. Again by Theorem \ref{theorem2.3} $Conv_{L^0}(G_n)$ is random totally bounded. Then for each $n\in \mathbb{N}$ there exist $\{H_k, k\in \mathbb{N}\}$ in $\Pi_{\mathcal{F}}$ and a sequence $\{G_{nk}: k\in \mathbb{N}\}$ of nonempty finite subsets of $Conv_{L^0}(G_n)$ such that $\tilde{I}_{H_k} Conv_{L^0}(G_n)\subset \tilde{I}_{H_k} (\sigma(G_{nk}))_{(\frac{\varepsilon}{2})}$ for each $k\in \mathbb{N}$. Thus,
	\begin{equation}\label{equation2.8}
		\aligned
		\tilde{I}_{A_n\bigcap H_k} (Conv_{L^0}(G_n))_{(\frac{\varepsilon}{2})} & = \tilde{I}_{A_n\bigcap H_k} (Conv_{L^0}(G_n)+ B(\theta, \frac{\varepsilon}{2}))  \\
		& =\tilde{I}_{A_n\bigcap H_k} Conv_{L^0}(G_n)+ \tilde{I}_{A_n\bigcap H_k} B(\theta, \frac{\varepsilon}{2})   \\
		& \subset \tilde{I}_{A_n\bigcap H_k} (\sigma(G_{nk}))_{(\frac{\varepsilon}{2})} + \tilde{I}_{A_n\bigcap H_k} B(\theta, \frac{\varepsilon}{2})  \\
		& \subset \tilde{I}_{A_n\bigcap H_k} \sigma(G_{nk})+ \tilde{I}_{A_n\bigcap H_k} B(\theta, \varepsilon)  \\
		& = \tilde{I}_{A_n\bigcap H_k} (\sigma(G_{nk}))_{(\varepsilon)}.
		\endaligned
	\end{equation}
	
	%\begin{align*}
	%  \tilde{I}_{A_n\bigcap H_k} (Conv_{L^0}(G_n))_{(\frac{\varepsilon}{2})} & = \tilde{I}_{A_n\bigcap H_k} (Conv_{L^0}(G_n)+ B(\theta, \frac{\varepsilon}{2}))  \\
	%   & =\tilde{I}_{A_n\bigcap H_k} Conv_{L^0}(G_n)+ \tilde{I}_{A_n\bigcap H_k} B(\theta, \frac{\varepsilon}{2})   \\
	%   & \subset \tilde{I}_{A_n\bigcap H_k} (\sigma(G_{nk}))_{(\frac{\varepsilon}{2})} + \tilde{I}_{A_n\bigcap H_k} B(\theta, \frac{\varepsilon}{2})  \\
	%   & \subset \tilde{I}_{A_n\bigcap H_k} \sigma(G_{nk})+ \tilde{I}_{A_n\bigcap H_k} B(\theta, \varepsilon)  \\
	%   & = \tilde{I}_{A_n\bigcap H_k} (\sigma(G_{nk}))_{(\varepsilon)}. \tag{1}
	%\end{align*}
	%\hfill (1) \tag{1}
	
	\par
	Now, from the fact that $\tilde{I}_{A_n}G\subset \tilde{I}_{A_n} (\sigma(G_n))_{(\frac{\varepsilon}{2})}$ for each $n\in \mathbb{N}$, one can have that $\tilde{I}_{A_n\bigcap H_k}G\subset \tilde{I}_{A_n\bigcap H_k} (\sigma(G_n))_{(\frac{\varepsilon}{2})}= \tilde{I}_{A_n\bigcap H_k} (\sigma(G_n)+ B(\theta, \frac{\varepsilon}{2}))$. It then immediately follows that
	\begin{align*}
		\tilde{I}_{A_n\bigcap H_k} Conv_{L^0}(G) & = Conv_{L^0}(\tilde{I}_{A_n\bigcap H_k} G) \\
		& \subset Conv_{L^0} (\tilde{I}_{A_n\bigcap H_k} (\sigma(G_n)+ B(\theta, \frac{\varepsilon}{2}))) \\
		& = \tilde{I}_{A_n\bigcap H_k} Conv_{L^0} (\sigma(G_n)+ B(\theta, \frac{\varepsilon}{2})) \\
		& \subset\tilde{I}_{A_n\bigcap H_k} (Conv_{L^0}(\sigma(G_n))+ Conv_{L^0} (B(\theta, \frac{\varepsilon}{2}))) \\
		& = \tilde{I}_{A_n\bigcap H_k} (Conv_{L^0}(G_n)+ B(\theta, \frac{\varepsilon}{2})) \\
		& = \tilde{I}_{A_n\bigcap H_k} (Conv_{L^0}(G_n))_{(\frac{\varepsilon}{2})} \\
		& \subset \tilde{I}_{A_n\bigcap H_k} (\sigma(G_{nk}))_{(\varepsilon)},
	\end{align*}
	by (\ref{equation2.8}).
	\par
	To sum up, for any $\varepsilon\in L^0_{++}(\mathcal{F})$ there exist $\{A_n\bigcap H_k: n,k\in \mathbb{N}\}$ in $\Pi_{\mathcal{F}}$ and a sequence $\{G_{nk}: n,k\in \mathbb{N}\}$ of nonempty finite subsets of $Conv_{L^0}(G)$ such that $\tilde{I}_{A_n\bigcap H_k} Conv_{L^0}(G)\subset \tilde{I}_{A_n\bigcap H_k} (\sigma(G_{nk}))_{(\varepsilon)} $ for each $n$ and $k\in \mathbb{N}$, namely, $Conv_{L^0}(G)$ is random totally bounded.
\end{proof}
~\
\par
Let $(E,\|\cdot\|)$ be an $RN$ module with base $(\Omega,\mathcal{F},P)$ and $G$ a $\mathcal{T}_{\varepsilon,\lambda}$--closed $L^0$--convex subset of $E$. $G$ is said to be $L^0$--convexly compact \cite{GZWW} if every family of $\mathcal{T}_{\varepsilon,\lambda}$--closed $L^0$--convex subsets of $G$ has a nonempty intersection whenever it has the finite intersection property. $G$ is said to have a random normal structure \cite{GZWG} if for each a.s. bounded  $\mathcal{T}_{\varepsilon,\lambda}$--closed $L^0$--convex subset $H$ of $G$ such that $D(H):= \bigvee \{\|x-y\|: x,y\in H\}> 0$, there exists a nondiametral point $h$ of $H$, namely $\bigvee\{\|h-x\|: x\in H\}< D(H)$ on $(D(H)> 0)^0$, where $D(H)$ is called the random diameter of $H$ and $(D(H)>0)^0$ is an arbitrarily chosen representative of $(D(H)>0)$. We will end this section with Corollary \ref{corollary2.8} below.
\begin{corollary}\label{corollary2.8}
	Let $(E,\|\cdot\|)$ be an $RN$ module over $\mathbb{K}$ with base $(\Omega,\mathcal{F},P)$ and $G$ a random sequentially compact $L^0$--convex subset of $E$. Then $G$ is $L^0$--convexly compact and has a random normal structure.
\end{corollary}
\begin{proof}
	$G$ is stable since it is $L^0$-convex. It then immediately follows from (6) of Lemma \ref{lemma2.5} that $G$ is $\mathcal{T}_{\varepsilon,\lambda}$--complete, so we can, without loss of generality, assume that $E$ is $\mathcal{T}_{\varepsilon,\lambda}$--complete. (Otherwise, we can consider the $\mathcal{T}_{\varepsilon,\lambda}$--completion of $E$ and note that $G$ is invariant in the process of $\mathcal{T}_{\varepsilon,\lambda}$--completion.) Let $f: (E, \mathcal{T}_{\varepsilon,\lambda}) \rightarrow (L^0(\mathcal{F},K), \mathcal{T}_{\varepsilon,\lambda})$ be a continuous module homomorphism, namely, $f\in E^*$. Then by (4) of Lemma \ref{lemma2.4} there exists $x_0\in G$ such that  $Re(f(x_0))= \bigvee \{Re(f(x)): x\in G\}$ since $G$ is $\sigma$--stable and $Re(f(\cdot)): G\rightarrow L^0(\mathcal{F})$ is $\mathcal{T}_{\varepsilon,\lambda}$--continuous and $\sigma$--stable. Hence $G$ is $L^0$--convexly compact by Theorem 2.21 of \cite{GZWW}.
	
	\par
	If $G$ does not have a random normal structure, then there exists a $\mathcal{T}_{\varepsilon,\lambda}$--closed $L^0$--convex subset $H$ of $G$ such that $H$ does not have a nondiametral point. We can assume, without loss of generality, that $H$ is just $G$ itself. Then for an arbitrarily chosen $x_1\in G$ it always holds that $\bigvee\{\|x_1-y\|: y\in G\}=D(G)$, denoted by $d$. By (4) of Lemma \ref{lemma2.4} there exists some $x_2\in G$ such that $\|x_1-x_2\|=d$. Further there exists some $x_3\in G$ such that $\|x_3- \frac{x_1+x_2}{2}\|=d$ since $\frac{x_1+x_2}{2}$ is still a diametral point of $G$, which necessarily implies that $\|x_3-x_i\|=d$ for each $i\in \{1,2\}$. We can inductively define a sequence $\{x_n, n\in \mathbb{N}\}$ in $G$ such that $\|x_{n+1}- \frac{x_1+x_2+\cdots+x_n}{n}\|=d$ for each $ n\in \mathbb{N}$ and $\|x_i-x_j\|=d$ for any $i,j\in \mathbb{N}$ with $i\neq j$. Let $\{x_{n_k}: k\in \mathbb{N}\}$ be any random subsequence of $\{x_n, n\in \mathbb{N}\}$. Then for any given $k$ and $l$ in $\mathbb{N}$ such that $k\neq l$ we will prove that $\|x_{n_k}-x_{n_l}\|=d$ as follows.
	\par
	Since $\{\omega\in \Omega~|~{n_k(\omega)=n_l(\omega)}\}= \emptyset$, then, for any $i$ and $j$ in $\mathbb{N}$, let $A_{ij}=\{\omega\in \Omega~|~n_k(\omega)=i$ and $n_l(\omega)=j\}$. One has $\Omega= \bigcup_{i,j\in \mathbb{N}}A_{ij}= \bigcup_{i\neq j}A_{ij}$. Further, $\|x_{n_k}- x_{n_l}\|= \|\sum_{i=1}^{\infty} \tilde{I}_{(n_k=i)} x_i- \sum_{j=1}^{\infty} \tilde{I}_{(n_l=j)} x_j\|= \sum_{i\neq j} \tilde{I}_{A_{ij}} \|x_i-x_j\|= (\sum_{i\neq j} \tilde{I}_{A_{ij}})d=d.$
	\par
	Thus it is impossible that there exists a random subsequence $\{x_{n_k}, k\in \mathbb{N}\}$ of $\{x_n, n\in \mathbb{N}\}$ such that $\{x_{n_k}, k\in \mathbb{N}\}$ converges, which contradicts the random sequential compactness of $G$.
\end{proof}
\section{Proofs of Theorems \ref{theorem1.12} and \ref{theorem1.15}}\label{section3}
\par
Before we give proofs of Theorems \ref{theorem1.12} and \ref{theorem1.15}, we will first give Lemmas \ref{lemma3.1}--\ref{lemma3.9} below as necessary auxiliaries.
\begin{lemma}\label{lemma3.1}
	Let $\mathcal{L}^0(\mathcal{F}, \mathbb{N})$ be the set of random variables from $(\Omega,\mathcal{F},P)$ to $\mathbb{N}$ and $H$ a nonempty subset of $\mathcal{L}^0(\mathcal{F}, \mathbb{N})$ such that $\sum^{\infty}_{k=1} I_{A_k} \cdot n_k$  belongs to $H$ for any sequence $\{n_k, k\in \mathbb{N}\}$ in $H$ and any $\{A_k, k\in \mathbb{N}\}$ in $\Pi_{\mathcal{F}}$, where $(\sum^{\infty}_{k=1} I_{A_k}  n_k)(\omega)= \sum^{\infty}_{k=1} I_{A_k}(\omega) \cdot n_k(\omega)$ for each $\omega\in \Omega$. Then there exists $h\in H$ such that $h(\omega)\leq h'(\omega)$ a.s. for any $h'\in H$, namely, $H$ has an essential least element.
\end{lemma}
\begin{proof}
	Let $h_0=essinf(H)$. Then $h_0\in \mathcal{L}^0(\mathcal{F}, \mathbb{N})$ (see the interpretation before Lemma \ref{lemma2.4}). For any $h_1$ and $h_2\in H,$ let $A=\{\omega\in \Omega: h_1(\omega)\leq h_2(\omega)\}$. Then $I_A\cdot h_1+ I_{A^c}\cdot h_2=min\{h_1,h_2\} $, so $min\{h_1,h_2\} \in H$, where $(min\{h_1,h_2\})(\omega)= min\{h_1(\omega), h_2(\omega)\}$ for each $\omega\in \Omega$. By the interpretation before Lemma \ref{lemma2.4}, there exists a sequence $\{h_n, n\in \mathbb{N}\}$ in $H$ such that $\{h_n, n\in \mathbb{N}\}$ converges a.s. to $h_0$ in a nonincreasing way. So by the Egoroff theorem for each $k\in \mathbb{N}$ there exists $B_k\in \mathcal{F}$ with $P(B^c_k)< \frac{1}{k}$ such that $\{h_n, n\in \mathbb{N}\}$ converges uniformly to $h_0$ on $B_k$. Consequently, there exists a sequence $\{l_k, k\in \mathbb{N}\}$ of positive integers such that $|h_{l_k}(\omega)- h_0(\omega)|< \frac{1}{2}$ for each $\omega\in B_k$ and such that $l_k< l_{k+1}$ for each $k\in \mathbb{N}$. Let $A_1=B_1$ and $A_k= B_k\backslash (\bigcup^{k-1}_{i=1} B_i)$ for each $k\geq 2$. Then $P(\bigcup ^{\infty}_{k=1} A_k)= P(\bigcup ^{\infty}_{k=1} B_k)=1$. Further, let $A_0= \Omega\backslash (\bigcup ^{\infty}_{k=1} A_k)$. Then $h:= I_{A_0\bigcup A_1}\cdot h_{l_1}+ \sum_{k=2}^{\infty} I_{A_k} \cdot h_{l_k} \in H$ and $|h(\omega)- h_0(\omega)|< \frac{1}{2}$ for each $\omega\in \bigcup_{k=1}^{\infty} A_k$, namely, $h(\omega)=h_0(\omega)$ a.s., which, of course, implies that $h(\omega)\leq h'(\omega)$ a.s. for any $h'\in H$.
\end{proof}

\begin{lemma}\label{lemma3.2}
	Let $(E,\|\cdot\|)$ be an $RN$ module with base $(\Omega,\mathcal{F}, P)$, $x_0\in E$, $G$ a $\sigma$--stable subset of $E$ and $\{x_n,n\in \mathbb{N}\}$ a sequence in $G$ such that $\{x_n,n\in \mathbb{N}\}$ converges in $\mathcal{T}_{\varepsilon,\lambda}$ to $x_0$. Then, for any given sequence $\{\varepsilon_k, k\in \mathbb{N}\}$ in $L^0_{++}(\mathcal{F})$ satisfying that $\{\varepsilon_k, k\in \mathbb{N}\}$ converges a.s. to 0 in a nonincreasing way, there exists a random subsequence $\{x_{n_k}, k\in \mathbb{N}\}$ of $\{x_n, n\in \mathbb{N}\}$ such that $\|x_{n_k}-x_0\|< \varepsilon_k$ on $\Omega$ for each $k\in \mathbb{N}$.
\end{lemma}
\begin{proof}
	Since $\overline{[\sigma(\{x_n,n\in \mathbb{N}\})]}_{\varepsilon,\lambda}= \overline{[\sigma(\{x_n,n\in \mathbb{N}\})]}_c$ by Theorem 3.12 of \cite{Guo3} and also since $x_0\in \overline{\{x_n,n\in \mathbb{N}\}}_{\varepsilon,\lambda}$ $\subset \overline{[\sigma(\{x_n,n\in \mathbb{N}\})]}_{\varepsilon,\lambda}$, there exists $y_1\in \sigma(\{x_n,n\in \mathbb{N}\})$ such that $\|y_1-x_0\|< \varepsilon_1$ on $\Omega$. Let $y_1= \sum_{n=1}^{\infty} \tilde{I}_{A_n} x_n$ for some $\{A_n, n\in \mathbb{N}\}$ in $\Pi_{\mathcal{F}}$ and further let $n_1= \sum_{n=1}^{\infty} I_{A_n}\cdot n$. Then $x_{n_1}= y_1$.
	\par
	We can assume, without loss of generality, that $P(\{\omega\in \Omega: n_1(\omega)=n\})= P(A_n)>0$ for each $n\in \mathbb{N}$. Now, for each fixed $n\in \mathbb{N}$, let $H'_n= \sigma(\{k: k\geq n+1\}): =\{\sum^{\infty}_{k=n+1} I_{D_k}\cdot k:  \{D_k: k\geq n+1\}\in \Pi_{\mathcal{F}} \}$. Then $\sigma(\{x_k, k\geq n+1\})= \{x_{m_l}: m_l\in H'_n\}$. Since $x_0\in \overline{[\sigma(\{ x_k,k\geq n+1\})]}_c= \overline{\{x_{m_l}:m_l\in H'_n \} }_c$, there exists some $m_l\in H'_n$ such that $\|x_{m_l}-x_0\|< \varepsilon_2$ on $\Omega$. Again, let $H_n= \{m_l: m_l\in H_n'$ and $\|x_{m_l}-x_0\|< \varepsilon_2$ on $\Omega \}$. It is obvious that $H_n$ satisfies the condition of Lemma \ref{lemma3.1}, and thus there exists $h_n\in H_n$ such that $h_n$ is an essential least element of $H_n$. Define $n_2=\sum_{n=1}^{\infty} I_{A_n}\cdot h_n$. We have $n_2(\omega)> n_1(\omega)$ for each $\omega\in \Omega$ and $\|x_{n_2}- x\|< \varepsilon_2$ on $\Omega$.
	\par
	Let $A_n'=(n_2=n)$ and assume, without loss of generality, that $P(A_n')>0$ for each $n\in \mathbb{N}$. We can similarly obtain an $n_3\in \mathcal{L}^0(\mathcal{F}, \mathbb{N})$ such that $n_3(\omega)> n_2(\omega)$ for each $\omega\in \Omega$ and $\|x_{n_3}-x_0\|< \varepsilon_3$ on $\Omega$ by replacing $\{A_n, n\in \mathbb{N}\}$ with $\{A'_n, n\in \mathbb{N}\}$ in the process of obtaining $n_2$. Thus, we can inductively obtain a random subsequence $\{x_{n_k}, k\in \mathbb{N}\}$ of $\{x_n, n\in \mathbb{N}\}$ such that $\|x_{n_k}- x_0\|< \varepsilon_k$ on $\Omega$ for each $k\in \mathbb{N}$.
\end{proof}
\par
The $(\varepsilon,\lambda)$--topology $\mathcal{T}_{\varepsilon,\lambda}$, as an abstract generalization of the topology of convergence in probability, is quite natural. It can help enable the verification of whether a mapping possesses some kind of continuity closely related to the $(\varepsilon,\lambda)$--topology. For example, one can easily observe that the $\mathcal{T}_{\varepsilon,\lambda}$--continuity of a mapping implies the random sequential continuity of the mapping. Compared with the $(\varepsilon,\lambda)$--topology, the locally $L^0$--convex topology $\mathcal{T}_c$ is too strong. Generally, it is not easy to determine whether a mapping is $\mathcal{T}_c$--continuous. With such an understanding, Lemma \ref{lemma3.3} below is quite surprising!

\begin{lemma}\label{lemma3.3}
	Let $(E_1, \|\cdot\|)$ and $(E_2, \|\cdot\|)$ be two $RN$ modules with base $(\Omega, \mathcal{F}, P)$, $G_1\subset E_1$ and $G_2\subset E_2$ two $\sigma$--stable subsets. Then a $\sigma$--stable mapping $T$ of $G_1$ to $G_2$ is random sequentially continuous iff it is $\mathcal{T}_c$--continuous. In particular, if a $\sigma$--stable mapping $T$ of $G_1$ to $G_2$ is $\mathcal{T}_{\varepsilon,\lambda}$--continuous, then $T$ must be $\mathcal{T}_c$--continuous.
\end{lemma}

\begin{proof}
	For any given $x_0\in G_1$ and any given $\varepsilon\in L^0_{++}(\mathcal{F})$, let $B_1(x_0,\varepsilon)= \{x\in G_1: \|x-x_0\|< \varepsilon$ on $\Omega \}$ and $B_2(T(x_0), \varepsilon)= \{y\in G_2: \|y-T(x_0)\|< \varepsilon$ on $\Omega\}$.
	\par
	If $T$ is $\mathcal{T}_c$--continuous at $x_0$, then for any positive integer $n$ there exists some $\varepsilon_n\in L^0_{++}(\mathcal{F})$ such that $T(B_1(x_0, \varepsilon_n))\subset B_2(T(x_0), \frac{1}{n})$. We can assume, without loss of generality, that $\{\varepsilon_n, n\in \mathbb{N}\}$ converges a.s. to $0$  nonincreasingly. If $\{x_n, n\in \mathbb{N}\}$ is any given sequence in $G_1$ such that $\{x_n, n\in \mathbb{N}\}$ converges in $\mathcal{T}_{\varepsilon,\lambda}$ to $x_0$, then by Lemma \ref{lemma3.2} there exists a random subsequence $\{x_{n_k}, k\in \mathbb{N}\}$ of $\{x_n, n\in \mathbb{N}\}$ such that $\|x_{n_k}-x_0\|< \varepsilon_k$ on $\Omega$ for each $k\in \mathbb{N}$, namely $x_{n_k}\in B_1(x_0, \varepsilon_k)$ for each $k\in \mathbb{N}$, which in turn implies that $\|T(x_{n_k})- T(x_0)\|< \frac{1}{k}$ on $\Omega$. It is obvious that $\{T(x_{n_k}), k\in \mathbb{N} \}$ converges in $\mathcal{T}_{\varepsilon,\lambda}$ to $T(x_0)$.
	\par
	Conversely, if $T$ is random sequentially continuous at $x_0\in G_1$, we will prove that $T$ is also $\mathcal{T}_c$--continuous at $x_0$ as follows.
	\par
	For any given $\varepsilon\in L^0_{++}(\mathcal{F})$ and any given positive integer $n\in \mathbb{N}$, let $\mathcal{A}_n= \{D\in \mathcal{F}: \tilde{I}_D T(B_1(x_0, \frac{1}{n})) $ $\subset \tilde{I}_D B_2(T(x_0), \varepsilon)\}$. Then it is easy to see that each $\mathcal{A}_n$ has the property: $\bigcup_k D_k\in \mathcal{A}_n$ if each $D_k\in \mathcal{A}_n$, from which $esssup(\mathcal{A}_n)$ $\in \mathcal{A}_n$, namely, $\mathcal{A}_n$ has an essential greatest element.
	\par
	Denote $esssup(\mathcal{A}_n)$ by $A_n$. It is obvious that $A_n\subset A_{n+1}$ a.s., namely, $P(A_n\backslash A_{n+1})=0$. We can assume, without loss of generality, that $A_n\subset A_{n+1}$ and $P(A_n)>0$ for each $n\in \mathbb{N}$. Now, we can assert that $P(\bigcup_{n=1}^{\infty} A_n)=1$ as follows: if $P(\bigcup_{n=1}^{\infty} A_n)<1$, then $P(A_0)>0$, where $A_0= \Omega\backslash (\bigcup_{n=1}^{\infty} A_n)$, and further let $\mathcal{B}_n= \{D\in \mathcal{F}: D\subset A_0, P(D)>0$ and there exists some $p\in B_1(x_0, \frac{1}{n})$ such that $\|T(p)-T(x_0)\|> \varepsilon$ on $D\}$. It is then clear that $\mathcal{B}_n \neq \emptyset$ (otherwise, $A_0\subset A_n$ a.s., leading to a contradiction). Further we prove that $\bigcup_{k=1}^{\infty} D_k$ still belongs to $\mathcal{B}_n$ if each $D_k\in \mathcal{B}_n$. In fact, let $p_k\in B_1(x_0, \frac{1}{n})$ be such that $\|T(p_k)- T(x_0)\|> \varepsilon$ on $D_k$ and $p= \tilde{I}_{(\bigcup^{\infty}_{k=1} D_k)^c} x_0+ \sum_{k=1}^{\infty} \tilde{I}_{C_k} p_k$, where $C_1=D_1$ and $C_k= D_k\backslash \bigcup^{k-1}_{i=1} D_i$ for each $k\geq 2$. Then $\|T(p)- T(x_0)\|= \|\tilde{I}_{(\bigcup^{\infty}_{k=1} D_k)^c} T(x_0)+ \sum_{k=1}^{\infty} \tilde{I}_{C_k} T(p_k)-T(x_0)\|$ (by $\sigma$--stability of $T$) $=\sum_{k=1}^{\infty} \tilde{I}_{C_k} \|T(p_k)-T(x_0)\|> \varepsilon$ on $\bigcup_{k=1}^{\infty} C_k= \bigcup_{k=1}^{\infty} D_k$, which, combined with $p\in B_1(x_0, \frac{1}{n})$ since $B_1(x_0, \frac{1}{n})$ is $\sigma$--stable, yields $\bigcup_{k=1}^{\infty} D_k\in \mathcal{B}_n$. Therefore, $B_n:= esssup(\mathcal{B}_n)\in \mathcal{B}_n$. Further, one also has $P(A_0\backslash B_n)=0$. Thus, we may regard $B_n=A_0$, since, if $P(A_0\backslash B_n)> 0$, then $A_0\backslash B_n\subset A_n $ a.s., which is again a contradiction. To sum up, for each $n\in \mathbb{N}$ there exists $p_n\in B_1(x_0, \frac{1}{n})$ such that $\|T(p_n)-T(x_0)\| >\varepsilon $ on $A_0$, namely, $\tilde{I}_{A_0} \|T(p_n)- T(p_0)\|= \|\tilde{I}_{A_0} T(p_n)- \tilde{I}_{A_0} T(x_0)\| >\varepsilon$ on $A_0$. But this is impossible: since $\|p_n-x_0\|< \frac{1}{n}$ on $\Omega$ for each $n\in \mathbb{N}$, $\{p_n, n\in \mathbb{N}\}$ clearly converges in $\mathcal{T}_{\varepsilon,\lambda}$ to $x_0$. Since $\tilde{I}_{A_0} T$ is still random sequentially continuous, there exists a random subsequence $\{p_{n_k}, k\in \mathbb{N}\}$ of $\{p_n, n\in \mathbb{N}\}$ such that $\{\|\tilde{I}_{A_0} T(p_{n_k})- \tilde{I}_{A_0} T(x_0)\|, k\in \mathbb{N}\}$ converges in probability to 0. But it is obvious that $\|\tilde{I}_{A_0} T(p_{n_k})- \tilde{I}_{A_0} T(x_0)\|> \varepsilon$ on $A_0$ for each $k\in \mathbb{N}$, which is again a contradiction.
	\par
	Now that we have proved $P(\bigcup_{n=1}^{\infty} A_n)=1$. We can assume that $\Omega= \bigcup_{n=1}^{\infty} A_n$. Let $F_1=A_1, F_n=A_n\backslash (\bigcup_{k=1}^{n-1} A_k)$ for each $n \geq 2$ and $\varepsilon_1= \sum_{n=1}^{\infty} \tilde{I}_{F_n}\cdot \frac{1}{n}$. Then $\varepsilon_1\in L^0_{++}(\mathcal{F})$. We will further show that $T(B_1(x_0, \varepsilon_1))\subset B_2(T(x_0), \varepsilon)$ as follows.
	\par
	We can assume, without loss of generality, that $P(F_n)>0$ for each $n\in \mathbb{N}$. For any given $p\in B_1(x_0, \varepsilon_1)$ and any given $n\in \mathbb{N}$, let $p_n= \tilde{I}_{F_n}p + \tilde{I}_{F_n^c} x_0$. Then it immediately follows from $p\in B_1(x_0, \varepsilon_1)$ that $\|p_n-x_0\|= \tilde{I}_{F_n} \|p-x_0\|< \frac{1}{n}$ on $\Omega$ since $\|p-x_0\|< \varepsilon_1$ on $\Omega$ (which, of course, implies that $\|p-x_0\|< \frac{1}{n}$ on $F_n$ for each $n\in \mathbb{N}$). By the definition of $A_n$, $\|\tilde{I}_{F_n} T(p_n)- \tilde{I}_{F_n} T(x_0)\|< \varepsilon$ on $\Omega$ for each $n\in \mathbb{N}$. On the other hand, it is clear that $p= \sum_{n=1}^{\infty} \tilde{I}_{F_n} \cdot p= \sum_{n=1}^{\infty} \tilde{I}_{F_n}  p_n$, and thus $\|T(p)-T(x_0)\|= \|T(\sum_{n=1}^{\infty} \tilde{I}_{A_n} p_n)- T(x_0)\|= \|\sum_{n=1}^{\infty} \tilde{I}_{F_n} T(p_n)- T(x_0)\|= \sum_{n=1}^{\infty} \tilde{I}_{F_n} \|T(p_n)- T(x_0)\|< \varepsilon$ on $\Omega$, namely, $T(p)\in B_2(T(x_0), \varepsilon)$.
	
	\par
	Finally, when a $\sigma$--stable mapping $T$ of $G_1$ to $G_2$ is $\mathcal{T}_{\varepsilon,\lambda}$--continuous, it must be random sequentially continuous, and hence $\mathcal{T}_c$--continuous.
\end{proof}
\par
Lemma \ref{lemma3.4} below is almost clear. Still, its proof is a requisite in order to assure that Lemma \ref{lemma3.5} below is indeed true.

\begin{lemma}\label{lemma3.4}
	Let $(E_1,\|\cdot\|)$ and $(E_2,\|\cdot\|)$ be two $RN$ modules with base $(\Omega,\mathcal{F},P)$, $G_1\subset E_1$ and $G_2\subset E_2$ two $\sigma$--stable subsets, and $f: G_1\rightarrow G_2$ a random sequentially continuous mapping. For a sequence $\{(x_1^m, x_2^m,\cdots, x_l^m), m\in \mathbb{N} \}$ in $G_1^l$, where $l$ is a fixed positive integer and $G^l_1$ is the $l$--th Cartesian power of $G_1$, if there exists a random subsequence $\{(x_1^{M_n^{(0)}}, x_2^{M_n^{(0)}}, \cdots, x_l^{M_n^{(0)}}),   n\in \mathbb{N}\}$ of which such that $\{x_i^{M_n^{(0)}}, n\in \mathbb{N}\}$ converges in $\mathcal{T}_{\varepsilon,\lambda}$ to some $y_i\in G_1$ for each $i\in \{1,2,\cdots, l \}$, then there exists a random subsequence $\{(x_1^{M_n}, x_2^{M_n}, \cdots, x_l^{M_n}),n\in \mathbb{N}  \}$ of $\{(x_1^{M_n^{(0)}}, x_2^{M_n^{(0)}}, \cdots, x_l^{M_n^{(0)}}),   n\in \mathbb{N}\}$ such that $\{x_i^{M_n}, n\in \mathbb{N}\}$ converges in $\mathcal{T}_{\varepsilon,\lambda}$ to $y_i$ and $\{f(x_i^{M_n}),n\in \mathbb{N} \}$ converges in $\mathcal{T}_{\varepsilon,\lambda}$ to $f(y_i)$ for each $i\in \{1,2,\cdots, l\}$.
\end{lemma}

\begin{proof}
	Since $\{x_1^{M_n^{(0)}}, n\in \mathbb{N}\}$ converges in $\mathcal{T}_{\varepsilon,\lambda}$ to $y_1$ and $f$ is random sequentially continuous at $y_1$, there exists a random subsequence $\{x_1^{M^{(0)} _{n_k}}, k\in \mathbb{N} \}$ of $\{x_1^{M_n^{(0)}}, n\in \mathbb{N}\}$ such that $\{f(x_1^{M^{(0)} _{n_k}}), k\in \mathbb{N} \}$ converges in $\mathcal{T}_{\varepsilon,\lambda}$ to $f(y_1)$, Let $M_k^{(1)}= M_{n_k}^{(0)}$ for each $k\in \mathbb{N}$, namely, $M_k^{(1)} (\omega)= M_{n_k(\omega)}^{(0)}(\omega)$ for each $\omega\in \Omega$.Then $\{(x_1^{M_n^{(1)}}, x_2^{M_n^{(1)}}, \cdots, x_l^{M_n^{(1)}}), $ $  n\in \mathbb{N}\}$ is a random subsequence of $\{(x_1^{M_n^{(0)}}, x_2^{M_n^{(0)}}, \cdots$, $x_l^{M_n^{(0)}}),   n\in \mathbb{N}\}$ such that $\{x_i^{M_n^{(1)}}, n\in \mathbb{N}\}$ still converges in $\mathcal{T}_{\varepsilon,\lambda}$ to $y_i$ for each $i\in \{1,2,\cdots, l \}$. Similarly, there exists a random subsequence $\{x_2^{M_{n_k}^{(1)}}, k\in \mathbb{N}\}$ of $\{x_2^{M_n^{(1)}}, n\in \mathbb{N}\}$ such that $\{f(x_2^{M_{n_k}^{(1)}}), k\in \mathbb{N} \}$ converges in $\mathcal{T}_{\varepsilon,\lambda}$ to $f(y_2)$. Let $M_k^{(2)}= M_{n_k}^{(1)}$ for each $k\in \mathbb{N}$, then $\{x_i^{M_n^{(2)}}, n\in \mathbb{N} \}$ still converges in $\mathcal{T}_{\varepsilon,\lambda}$ to $y_i$ for each $i\in \{1,2,\cdots, l \}$ and $\{f(x_i^{M_{n}^{(2)}}), n\in \mathbb{N} \}$ converges in $\mathcal{T}_{\varepsilon,\lambda}$ to $f(y_i)$ for each $i\in \{1,2\}$. We can inductively obtain $\{M_k^{(j)}, k\in \mathbb{N} \}$ for each $j\in \{1,2,\cdots, l \}$ such that $M_k^{(j)}= M_{n_k}^{(j-1)}, k\in \mathbb{N}$ and each $j\in \{1,2,\cdots, l \}$. Finally, taking $M_n= M_n^{(l)}$ for each $n\in \mathbb{N}$, we obtain our desired result.
\end{proof}
\par
The topological module structure under $\mathcal{T}_{\varepsilon,\lambda}$ of $L^0$--free modules of finite rank was characterized by Lemma 3.4 of \cite{GP}; the algebraic structure of finitely generated $L^0$--modules was further characterized by Theorem 1.1 of \cite{GS}; $L^0$--affine independence and $L^0$--simplexes were systematically studied in \cite{DKKS}. Now, Theorem 2.3 of \cite{DKKS} can be restated as follows: Let $S= Conv_{L^0} (X_1,\cdots , X_N)$ be an $L^0$--simplex in $L^0(\mathcal{F}, \mathbb{R}^d)$ and $f: S\rightarrow S$ a $\sigma$--stable a.s. sequentially continuous mapping. Then $f$ has a fixed point. An obvious and trivial observation was already made in \cite[Proposition 4.2]{GZWW}, which shows that Theorem 2.3 of \cite{DKKS} remains valid for a $\sigma$--stable $\mathcal{T}_{\varepsilon,\lambda}$--continuous mapping $f$. The following observation, namely, Lemma \ref{lemma3.5} below is nontrivial and also crucial for our work here since it shows that Theorem 2.3 of \cite{DKKS} remains valid for a $\sigma$--stable $\mathcal{T}_c$--continuous mapping $f$.

\begin{lemma}\label{lemma3.5}
	Let $S=Conv_{L^0}(X_1,\cdots,X_N)$ be an $L^0$--simplex in $L^0(\mathcal{F}, \mathbb{R}^d)$ and $f:S \rightarrow S$ a $\sigma$--stable random sequentially continuous mapping. Then $f$ has a fixed point.
\end{lemma}

\begin{proof}
	This proof will be finished by slightly modifying the proof of Theorem 2.3 of \cite{DKKS}. That proof consists of two parts: the first part gives an $L^0$--module version of Sperner's lemma, whose proof does not involve the continuity of $f$. Only its second part uses the a.s. sequential continuity of $f$, which, combined with the property of the sequence $\{S^m, m\in \mathbb{N}  \}$ of completely labeled $L^0$--simplexes, eventually produces a fixed point of $f$. Therefore, we only need to modify the second last paragraph of the proof of Theorem 2.3 in \cite{DKKS}. For this, for each $m\in \mathbb{N}$ let $S^m= Conv _{L^0} (V_1^m,\cdots, V_N^m)$  with $P(\{\omega: \varphi^{m-1} (V_j^m)(\omega)=j \} )=1$ for every $j=1,\cdots, N$ be the same as in the end of the third last paragraph of the proof of Theorem 2.3 of \cite{DKKS}.
	\par
	Since $S$ is random sequentially compact and $\sigma$--stable. By $\text{Lemma} ~2$ of \cite{KS} $\{V_1^m, m\in \mathbb{N}\}$, as a sequence in $S$, admits a random subsequence $\{V_1^{M_n^{(0)}}, n\in \mathbb{N}\}$ convergent a.s. to some $Y$ in $S$, where $\{M_n^{(0)}, n\in \mathbb{N}\}$ is exactly $ `` \{M_n,n\in \mathbb{N}\}$ given in the second last paragraph of the proof of Theorem 2.3 of \cite{DKKS} ''. Thus $\{V_k^{M_n^{(0)}}, n\in \mathbb{N} \}$ converges a.s. (of course, also in $\mathcal{T}_{\varepsilon,\lambda}$, namely, in probability) to $Y$ for each $k=1,\cdots,N$, as argumented in \cite{DKKS}. Now, by making use of Lemma \ref{lemma3.4} there exists a random subsequence $\{(V_1^{M_n}, V_2^{M_n},\cdots, V_N^{M_n} ), n\in \mathbb{N} \}$ of $\{(V_1^{M_n^{(0)}}, V_2^{M_n^{(0)}},\cdots, V_N^{M_n^{(0)}}),n\in \mathbb{N} \}$ such that $\{V_k^{M_n}, n\in \mathbb{N}\}$ converges in probability to $Y$ and $\{f(V_k^{M_n}), n\in \mathbb{N} \}$ also converges in probability to $f(Y)$.
	\par
	The remaining part of the proof is completely similar to the last paragraph of Theorem 2.3 of \cite{DKKS}. One can see that $f(Y)= Y$.
\end{proof}
\par
By the same method as used in the proof of Proposition 3.1 of \cite{DKKS}, we can establish Lemma \ref{lemma3.6} below. So we omit its proof, where we only need to remind the reader of the facts: ``sequentially closed'' in \cite[Proposition 3.1]{DKKS} is exactly $\mathcal{T}_{\varepsilon,\lambda}$--closed, that is, closed in the topology of convergence in probability, and ``bounded'' there is also exactly a.s. bounded. Besides, we bear in mind the fact: an $L^0$--convex and $\mathcal{T}_{\varepsilon,\lambda}$--closed subset of a $\mathcal{T}_{\varepsilon,\lambda}$--complete $RN$ module (e.g., $L^0(\mathcal{F},\mathbb{R}^d)$) is always $\sigma$--stable.
\begin{lemma}\label{lemma3.6}
	Let $G$ be an $L^0$--convex, $\mathcal{T}_{\varepsilon,\lambda}$--closed and a.s. bounded subset of $L^0(\mathcal{F},\mathbb{R}^d)$ and $f: G \rightarrow G$ a $\sigma$--stable and random sequentially continuous mapping. Then $f$ has a fixed point.
	
\end{lemma}
\begin{lemma}\label{lemma3.7}
	Let $G$ be an $L^0$--convex, $\mathcal{T}_{\varepsilon,\lambda}$--closed and a.s. bounded subset of a finitely generated $RN$ module $(E,\|\cdot\| ) $ over $\mathbb{K}$ with base $(\Omega, \mathcal{F}, P)$ and $f: G \rightarrow G$ a $\sigma$--stable and random sequentially continuous mapping, where we recall that $E$ is finitely generated means that there exists a finite subset of $E$, denoted by $\{x_1,x_2, \cdots, x_d\}$ for some $d\in \mathbb{N}$, such that $E= \{\sum^d_{i=1} \xi_i \cdot x_i: \xi_i\in L^0(\mathcal{F}, \mathbb{K}) $ for each $i\in \{1,2, \cdots, d\}\}$. Then $f$ has a fixed point.
\end{lemma}

\begin{proof}
	By Corollary 2.25 of \cite{GZWW}, $(E, \mathcal{T}_{\varepsilon,\lambda})$ is isomorphic to a closed submodule of $(L^0(\mathcal{F}, \mathbb{K}^d), \mathcal{T}_{\varepsilon,\lambda} )$ in the sense of topological modules. If the topological module isomorphism is denoted by $T$, then $T(G)$ is an $L^0$--convex,  $\mathcal{T}_{\varepsilon,\lambda}$--closed and a.s. bounded subset of $T(E)$ (and hence also of $L^0(\mathcal{F},\mathbb{K}^d)$) since $T$ preserves $L^0$--convexity, $\mathcal{T}_{\varepsilon,\lambda}$--closedness and a.s. boundedness of $G$. Now, we consider the mapping $\hat{f}: T(G) \rightarrow T(G)$ defined by $\hat{f}(y)= (T\circ f \circ T^{-1})(y)$ for each $y\in T(G)$. It is clear that $\hat{f}$ is $\sigma$--stable and random sequentially continuous. When $\mathbb{K} = \mathbb{R}$, we can directly apply Lemma \ref{lemma3.6} to $L^0(\mathcal{F},\mathbb{R}^d)$, and when $\mathbb{K} = \mathbb{C}$ we can apply  Lemma \ref{lemma3.6} to $L^0(\mathcal{F},\mathbb{R}^{2d})$. To sum up, we can have that $\hat{f}$ has a fixed point $\hat{p}$ in $T(G)$. It is obvious that $T^{-1}(\hat{p})$ is a fixed point of $f$.
\end{proof}
\par
The classical Schauder projection is defined on a finite $\varepsilon$--net in a normed space, which, combined with the classical Schauder approximation theorem, plays a crucial role in the proof of the classical Schauder fixed point theorem \cite[pp.116-117]{GD}. However, in the case of an $RN$ module $(E,\|\cdot\|)$ with base $(\Omega, \mathcal{F}, P)$, let $G= \{x_1,x_2,\cdots, x_n \}$ be an $n$--element finite subset of $E$ and $\varepsilon$ a given element in $L^0_{++}(\mathcal{F})$. We now define the corresponding Schauder projection $P_{G, \varepsilon}: (\sigma(G))_{(\varepsilon)} \rightarrow Conv_{L^0} (\{x_1,x_2,\cdots, x_n \})$ (denoted by $Conv_{L^0} (x_1,x_2,\cdots, x_n )$) as follows:
$$P_{G, \varepsilon}(x)= \frac{1}{\sum^n_{i=1} u_i(x)} \sum^n_{i=1} u_i(x) x_i,~\forall x\in (\sigma(G))_{(\varepsilon)},$$
where, for each $i\in \{1,2,\cdots, n\}, u_i: (\sigma(G))_{(\varepsilon)} \rightarrow L^0_+(\mathcal{F})$ is defined by $u_i(x)= max\{0, \varepsilon- \|x-x_i\|\}$ for each $x\in (\sigma(G))_{(\varepsilon)}$. Since $P_{G, \varepsilon}$ is defined on $(\sigma(G))_{(\varepsilon)}$ rather than on $G_{(\varepsilon)}$, we need Lemma \ref{lemma3.8} below to summarize the main properties of $P_{G, \varepsilon}$, which will be required in the remainder of this paper.

\begin{lemma}\label{lemma3.8}
	Let $(E,\|\cdot\|)$ be a $\sigma$--stable $RN$ module over $\mathbb{K}$ with base $(\Omega, \mathcal{F}, P)$, $G=\{x_1,x_2,\cdots, x_n \}$ an $n$--element finite subset of $E$, $\varepsilon\in L^0_{++}(\mathcal{F})$ and $P_{G, \varepsilon}$ the corresponding Schauder projection. Then we have the following:
	\begin{enumerate}[$(a).$]
		\item $P_{G, \varepsilon}$ is well defined, $\sigma$--stable and $\mathcal{T}_{\varepsilon,\lambda}$--continuous $($in fact, also $\mathcal{T}_c$--continuous$)$.\\
		\item $\|x- P_{G, \varepsilon}(x)\|< \varepsilon$ on $\Omega$ for each $x\in (\sigma(G))_{(\varepsilon)}$.\\
		\item If $G$ is symmetric with respect to $\theta$, i.e., $G=\{x_1,x_2,\cdots, x_k, -x_1,-x_2,\cdots, -x_k \}$ for some $k\in \mathbb{N}$, then $(\sigma(G))_{(\varepsilon)}=- (\sigma(G))_{(\varepsilon)}$ and $P_{G, \varepsilon} (-x)= -P_{G,\varepsilon}(x)$ for each $x\in (\sigma(G))_{(\varepsilon)}$.
	\end{enumerate}
\end{lemma}
\begin{proof}
	(a). We first prove that $P_{G,\varepsilon}$ is well defined. For this, we only need to check $\sum^n_{i=1} u_i(x) \in L^0_{++}(\mathcal{F})$ for each $x\in (\sigma(G))_{(\varepsilon)}$. In fact, $x\in (\sigma(G))_{(\varepsilon)}$ iff $d(x, \sigma(G))< \varepsilon$ on $\Omega$ by (2) of Lemma \ref{lemma2.4} iff $d(x, G)< \varepsilon$ on $\Omega$ by (1) of Lemma \ref{lemma2.4}. Further, since $d(x, G)= \bigwedge \{\|x-x_i\|: i=1\sim n \} = \|x-\sum^n_{i=1} \tilde{I}_{A_i} x_i\|= \sum^n_{i=1} \tilde{I}_{A_i} \|x-x_i\|$ for some finite partition $\{A_i, i=1 \sim n\}$ of $\Omega$ to $\mathcal{F} $ by (3) of Lemma \ref{lemma2.4}, so $d(x, G)< \varepsilon$ on $\Omega$ iff $\|x-x_i\|< \varepsilon$ on $A_i$ for each $i=1 \sim n$, that is, $d(x, G)< \varepsilon$ on $\Omega$ iff $u_i(x)= max \{0, \varepsilon- \|x-x_i\|\}>0$ on $A_i$ for each $i=1 \sim n$. This is equivalent to saying that $\sum^n_{i=1} u_i(x)\in L^0_{++}(\mathcal{F})$ iff $x\in (\sigma(G))_{(\varepsilon)}$.
	\par
	Since $(\sigma(G))_{(\varepsilon)}$ is clearly $\sigma$--stable and each $u_i$ is also $\sigma$--stable, let $\{B_k, k\in \mathbb{N}\}\in \Pi_{\mathcal{F}}$ and $\{y_k, k\in \mathbb{N} \}$ a sequence in $(\sigma(G))_{(\varepsilon)}$. Then
	\begin{align*}
		&P_{G,\varepsilon}(\sum_{k=1}^{\infty} \tilde{I}_{B_k} y_k )\\
		&=\frac{1}{\sum^n_{i=1 }u_i (\sum_{k=1}^{\infty} \tilde{I}_{B_k} y_k)}  \sum^n_{i=1} u_i (\sum_{k=1}^{\infty} \tilde{I}_{B_k} y_k) x_i\\
		&=\frac{1}{\sum_{k=1}^{\infty} \tilde{I}_{B_k} (\sum^n_{i=1 }u_i (y_k)) } \sum_{k=1}^{\infty} \tilde{I}_{B_k} (\sum^n_{i=1 }u_i (y_k) x_i)\\
		&=\sum_{k=1}^{\infty} \tilde{I}_{B_k} \frac{1}{\sum^n_{i=1 }u_i (y_k)} \sum^n_{i=1 }u_i (y_k) x_i\\
		&=\sum_{k=1}^{\infty} \tilde{I}_{B_k} \cdot P_{G, \varepsilon} (y_k),
	\end{align*}
	which shows that $P_{G, \varepsilon}$ is $\sigma$--stable.
	\par
	Further, since each $u_i$ is both $\mathcal{T}_{\varepsilon,\lambda}$--continuous and $\mathcal{T}_c$--continuous, and since $(E,\|\cdot\|)$  is a topological module both over the topological algebra $(L^0(\mathcal{F}, \mathbb{K}), \mathcal{T}_{\varepsilon,\lambda} )$ and over the topological ring $(L^0(\mathcal{F}, \mathbb{K}), \mathcal{T}_c )$, $P_{G,\varepsilon}$ must be both $\mathcal{T}_{\varepsilon,\lambda}$--continuous and $\mathcal{T}_c$--continuous.
	\par
	(b). For each $x\in(\sigma(G) )_{(\varepsilon)}$,  since $\|x- P_{G, \varepsilon}(x)\| \leq \frac{\sum^n_{i=1} u_i(x) \|x-x_i\|}{\sum^n_{i=1} u_i(x)}=\sum^n_{i=1} $ $\frac{u_i(x)}{\sum^n_{i=1} u_i(x)}\|x-x_i\|$, and since $u_i(x) = max\{0, \varepsilon- \|x-x_i\| \}= I_{(\|x-x_i\|< \varepsilon)} u_i(x)$ for each $i=1 \sim n$, then $\|x-P_{G, \varepsilon}(x)\|\leq \sum^n_{i=1} \frac{u_i(x)}{\sum^n_{i=1}u_i(x)} I_{(\|x-x_i\|< \varepsilon)} \|x-x_i\|< \varepsilon$ on $\Omega$ since $I_{(\|x-x_i\|< \varepsilon)} \|x-x_i\|< \varepsilon$ on $\Omega$ for each $i=1 \sim n$.
	\par
	(c) is clear.
\end{proof}

\begin{lemma}\label{lemma3.9}
	Let $(E,\|\cdot\|)$ be a $\sigma$--stable $RN$ module with base $(\Omega, \mathcal{F}, P)$, $G$ a $\sigma$--stable $L^0$--convex subset of $E$, $T: G\rightarrow G$ a $\sigma$--stable random sequentially continuous mapping and $\varepsilon\in L^0_{++}(\mathcal{F})$ such that there exists some finite subset $G^{\varepsilon}$ of $G$ satisfying $G\subset (\sigma(G^{\varepsilon}))_{(\varepsilon)}$. Then there exists $x_{\varepsilon}\in G$ such that $\|T(x_{\varepsilon})- x_{\varepsilon}\| < \varepsilon$ on $\Omega$.
\end{lemma}

\begin{proof}
	Let $P_{G^{\varepsilon}, \varepsilon}$ be the corresponding Schauder projection and $T_{\varepsilon}= P_{G^{\varepsilon}, \varepsilon} \circ T$. Then $T_{\varepsilon}$ is a well defined, $\sigma$--stable and random sequentially continuous mapping sending $Conv_{L^0} (G^{\varepsilon})$ into $Conv_{L^0} (G^{\varepsilon})$. Further, let $E_1$ be the $RN$ module generated by $G^{\varepsilon}$. Then it is obvious that $Conv_{L^0} (G^{\varepsilon})$ is $\sigma$--stable, $\mathcal{T}_{\varepsilon,\lambda}$--closed and a.s. bounded $L^0$--convex subset of $E_1$, and thus by Lemma \ref{lemma3.7} there exists some $x_{\varepsilon}\in Conv_{L^0} (G^{\varepsilon})$ such that $T_{\varepsilon} (x_{\varepsilon}) = x_{\varepsilon}$, which in turn implies that $\|x_{\varepsilon} -T(x_{\varepsilon})\| = \|P_{G^{\varepsilon}, \varepsilon} (T(x_{\varepsilon}))- T(x_{\varepsilon})\| < \varepsilon$ on $\Omega$ by Lemma \ref{lemma3.8}. Finally, it is also obvious that $x_{\varepsilon} \in G$ since $Conv_{L^0} (G^{\varepsilon})\subset G.$
\end{proof}

\par
We are now ready to prove our first main result, Theorem \ref{theorem1.12}, as follows.

\begin{proof}[\textbf{Proof of Theorem \ref{theorem1.12}}]
First, $G$ is $\mathcal{T}_{\varepsilon,\lambda}$--complete by (6) of Lemma \ref{lemma2.5} and, thus, $G$ is random totally bounded by Theorem \ref{theorem2.3}. Similar to the proof of Corollary \ref{corollary2.8}, we can also assume that $E$ is $\mathcal{T}_{\varepsilon,\lambda}$--complete. Then $E$ is $\sigma$--stable. For each fixed $k\in \mathbb{N}$, there exist a sequence $\{G_n, n\in \mathbb{N} \}$ of nonempty finite subsets of $G$ and $\{A_n,n\in \mathbb{N}\}$ in $\Pi_{\mathcal{F}}$ such that $\tilde{I}_{A_n} G\subset \tilde{I}_{A_n} (\sigma (G_n))_{(\frac{1}{k})}$ for each $n\in \mathbb{N}$.
\par
We can assume, without loss of generality, that $\theta\in G$ (otherwise, let $p_0\in G$, we can consider $G'= G-p_0$ and $T': G'\rightarrow G'$ defined by $T'(p)= T(p+p_0) - p_0$ for each $p\in G'$). Then $\tilde{I}_{A_n} G\subset G$ and
\begin{align*}
	\tilde{I}_{A_n} (\sigma (G_n))_{(\frac{1}{k})} &  = \tilde{I}_{A_n} (\sigma(G_n)+ B(\theta, \frac{1}{k}) ) \\
	& =\tilde{I}_{A_n} \sigma(G_n)+ \tilde{I}_{A_n} B(\theta, \frac{1}{k}) \\
	& \subset \tilde{I}_{A_n} \sigma(G_n) + B(\theta, \frac{1}{k}) \\
	& = \sigma (\tilde{I}_{A_n} G_n)+ B(\theta, \frac{1}{k}) \\
	& = (\sigma (\tilde{I}_{A_n} G_n))_{(\frac{1}{k})},
\end{align*}
and hence $\tilde{I}_{A_n} G\subset (\sigma (\tilde{I}_{A_n} G_n))_{(\frac{1}{k})}$ for each $n\in \mathbb{N}$. Further, for each $n\in \mathbb{N}$ we consider the mapping $T_n=$ the restriction of $\tilde{I}_{A_n} T$ to $\tilde{I}_{A_n} G$, where $\tilde{I}_{A_n} T : G \rightarrow \tilde{I}_{A_n} G$ is defined by $(\tilde{I}_{A_n} T)(g)= \tilde{I}_{A_n}(T(g))$ for each $g\in G$. Then each $\tilde{I}_{A_n} G$ is still $\sigma$--stable and random sequentially compact, and each $T_n$ remains $\sigma$--stable and random sequentially continuous. Thus, applying Lemma \ref{lemma3.9} to $T_n$, $\tilde{I}_{A_n} G$ and $\frac{1}{k}$ yields some $x_{kn} \in \tilde{I}_{A_n} G$ such that $\|T_n(x_{kn})- x_{kn}\| < \frac{1}{k}$ on $\Omega$. Further, since $x_{kn}= \tilde{I}_{A_n} x_{kn}$, we have that $\tilde{I}_{A_n} \|T(x_{kn})- x_{kn}\| = \|\tilde{I}_{A_n} T(x_{kn})- \tilde{I}_{A_n} x_{kn}\|= $ $\|T_n(x_{kn})- x_{kn}\|< \frac{1}{k}$ on $\Omega$, that is, $\|T(x_{kn})- x_{kn}\|< \frac{1}{k}$ on $A_n$ for each $n\in \mathbb{N}$.
\par
Now, let $x_k= \sum_{n=1}^{\infty} \tilde{I}_{A_n} x_{kn}$. Then $x_k\in G$ and further by the $\sigma$--stability of $T$, we have $\|T(x_k)-x_k\|= \sum_{n=1}^{\infty} $ $\tilde{I}_{A_n} \|T(x_{kn}) -x_{kn}\|< \frac{1}{k}$ on $\Omega$, which implies that $\{T(x_k)-x_k, k\in \mathbb{N} \}$ converges in $\mathcal{T}_{\varepsilon,\lambda}$ to $\theta$. Since $G$ is random sequentially compact and $\sigma$--stable, there exists a random subsequence $\{x_{n_k},k\in \mathbb{N} \}$ of $\{x_k,k\in \mathbb{N} \}$ such that $\{x_{n_k},k\in \mathbb{N} \}$ converges in $\mathcal{T}_{\varepsilon,\lambda}$ to some point $x\in G$. Further, $\{T(x_{n_k}), k\in \mathbb{N} \}$ also converges in $\mathcal{T}_{\varepsilon,\lambda}$ to $x$ since $\{T(x_{n_k})-x_{n_k}, k\in \mathbb{N} \}$ converges in $\mathcal{T}_{\varepsilon,\lambda}$ to $\theta$. Finally, since $T$ is random sequentially continuous, there exists a random subsequence $\{y_l, l\in \mathbb{N}\}$ of $\{x_{n_k},k\in \mathbb{N} \}$ such that $\{ T(y_l), l\in \mathbb{N}\}$ converges in $\mathcal{T}_{\varepsilon,\lambda}$ to $T(x)$. As $\{ T(y_l), l\in \mathbb{N}\}$ is also a random subsequence of $\{T(x_{n_k}), k\in \mathbb{N} \}$, we eventually arrive at $T(x)=x.$
\end{proof}

\par
To prove Theorem \ref{theorem1.15}, we first give Lemma \ref{lemma3.10} below.
\begin{lemma}\label{lemma3.10}
	Let $(E,\|\cdot\|)$ be a $\mathcal{T}_{\varepsilon,\lambda}$--complete $RN$ module with base $(\Omega, \mathcal{F}, P)$, $G$ a $\mathcal{T}_{\varepsilon,\lambda}$--closed $L^0$--convex subset of $E$ and $T: G\rightarrow G$ a $\mathcal{T}_{\varepsilon,\lambda}$--continuous or $\mathcal{T}_c$--continuous $\sigma$--stable mapping such that $T(G)$ is random relatively sequentially compact. Then $T$ has a fixed point.
\end{lemma}

\begin{proof}
	Since $G$ is $\sigma$--stable by (1) of Lemma \ref{lemma2.5}, and $T(G)$ is $\sigma$--stable, then $M:= (T(G))^-_{\varepsilon,\lambda}$ is random sequentially compact . It is obvious that $G':= [Conv_{L^0} (T(G))]^-_{\varepsilon,\lambda}= [Conv_{L^0} (M)]^-_{\varepsilon,\lambda}$, and thus $G'$ is a random sequentially compact $L^0$-convex set by Corollary \ref{corollary2.7}. Since $T(G')\subset G'$, there exists some $x\in G'$ such that $T(x)=x$ by Theorem \ref{theorem1.12}.
\end{proof}
\par
We are now ready to prove our second main result, Theorem \ref{theorem1.15}.

\begin{proof}[\textbf{Proof of Theorem \ref{theorem1.15}}]
Since $G$  is $\sigma$--stable by (1) of Lemma \ref{lemma2.5} and $S$ is $L^0$--Lipschitzian, $S$ is $\sigma$--stable by Lemma 2.11 of \cite{GZWG}. Denote the range of $I-S$ by $Z$, where $I$ is the identity mapping on $E$. Then $Z$ is $\sigma$--stable. Further, it is easy to check that $(I-S)^{-1}: Z \rightarrow G$ exists and is, in particular, $L^0$--Lipschitzian with the $L^0$--Lipschitzian constant $(1-\alpha)^{-1}$. Thus, $(I-S)^{-1}$ is both $\mathcal{T}_{\varepsilon,\lambda}$--continuous and $\mathcal{T}_c$--continuous. In fact, one can also have that $Z$ is $\mathcal{T}_{\varepsilon,\lambda}$--closed.
\par
Now, by condition (3) of this theorem, for any given $y\in G$ $S+T(y)$ is a random contractive mapping of $G$ into $G$. By the $\mathcal{T}_{\varepsilon,\lambda}$--completeness of $G$ it follows from \cite[Corollary 3.14]{GWYZ} that there exists a unique $x_y\in G$ such that $S(x_y)+ T(y)= x_y$. Therefore, $x_y= (I-S)^{-1} (T(y))$, namely, $T(y)\in Z$  and $(I-S)^{-1} (T(y))= x_y\in G$ for each $y\in G$. Consider the mapping $T': G\rightarrow G$ defined by $T'= (I-S)^{-1} \circ T$. We have that $T'$ is a $\mathcal{T}_{\varepsilon,\lambda}$--continuous or $\mathcal{T}_c$--continuous $\sigma$--stable mapping such that $T'(G)$ is random relatively sequentially compact. By Lemma \ref{lemma3.10} $T'$ has a fixed point $x\in G$, and $x$ is just a fixed point of $S+T$.
\end{proof}

\section{Proof of Theorem \ref{theorem1.19}}\label{section4}

\par
For the proof of Theorem \ref{theorem1.19}, we still need to study the topological structure of an $L^0$--simplex and the continuity of an $L^0$--payoff function. Besides, for future further generalization of Theorem \ref{theorem1.19}, here we also first study the basic properties of $L^0$--convex functions and the set of $L^0$--extreme points of a nonempty subset of an $L^0$--module.

\par
Let us recall the notion of a regular $L^0$--module. A left module $E$ over the algebra $L^0(\mathcal{F}, \mathbb{K})$ (briefly, an $L^0(\mathcal{F}, \mathbb{K})$--module) is said to be regular if any two elements $x$ and $y$ in $E$ are equal iff there exists some countable partition $\{A_n, n\in \mathbb{N} \}$ of $\Omega$ to $\mathcal{F}$ (namely, $\{A_n, n\in \mathbb{N} \}$ in $\Pi_{\mathcal{F}}$) such that $\tilde{I}_{A_n} x= \tilde{I}_{A_n} y$ for each $n\in \mathbb{N}$. In the subsequent part of this section, we always assume that all the $L^0$--modules occurring in the section are regular. As pointed out in \cite{Guo3,WGL,WZZ}, all $RN$ modules and random locally convex modules are regular.

\begin{definition}\label{definition4.1}
	Let $E$ be an $L^0(\mathcal{F}, \mathbb{K})$--module and $G$ an $L^0$--convex subset of $E$. A mapping $f: G\rightarrow \bar{L}^0(\mathcal{F})$  is called
	\begin{enumerate}[(1)]
		\item a proper function if $f(x)> -\infty$ on $\Omega$ for each $x\in G$ and there exists some $x_0\in G$ such that $f(x_0)< +\infty$ on $\Omega$.
		\item a proper $L^0$--convex function if $f$ is proper and $f(\lambda x+ (1-\lambda) y) \leq \lambda f(x)+ (1-\lambda) f(y)$ for any $x$ and $y$ in $G$ and any $\lambda\in L^0(\mathcal{F}, [0,1])$, where we make the convention: $0\cdot (+\infty)=0$.
		\item a proper $L^0$--quasiconvex function if $f$ is proper and $\{x\in G: f(x) \leq \alpha\}$ is $L^0$--convex for any $\alpha\in L^0(\mathcal{F})$.
		\item a proper $L^0$--concave (resp., a proper $L^0$--quasiconcave) function if $-f$ is a proper $L^0$--convex (resp., accordingly, a proper $L^0$--quasiconvex) function.
	\end{enumerate}
\end{definition}

\begin{remark}\label{remark4.2}
	Following the convention from conditional convex risk measures, according to (4) of Definition \ref{definition4.1} properness of a proper $L^0$--concave (or $L^0$--quasiconcave) function $f$ means $f(x)<+\infty$ on $\Omega$ for each $x\in G$ and there exists some $x_0\in G$ such that $f(x_0)> -\infty$ on $\Omega$. By the way, a mapping $f: G\rightarrow \bar{L}^0(\mathcal{F})$ is both a proper $L^0$--convex function and a proper $L^0$--concave function iff $f$ takes its values in $L^0(\mathcal{F})$ and $f$ is $L^0$--affine.
\end{remark}

\par
To provide the convenience for the proof of Lemma \ref{lemma4.3} below, let us recall the notion of a generalized inverse for an element $\xi$ in $L^0(\mathcal{F}, \mathbb{K} )$: for an arbitrarily chosen representative $\xi^0$ of $\xi$, define $(\xi^0)^{[-1]}: \Omega \rightarrow \mathbb{K}$ by $(\xi^0)^{[-1]} (\omega)= 1/\xi^0(\omega)$ if $\xi^0(\omega) \neq 0$ and $0$ otherwise, then the equivalence class of $(\xi^0)^{[-1]}$ is called the generalized inverse of $\xi$, denoted by $\xi^{[-1]}$. It is obvious that $\xi \cdot \xi^{[-1]} = \xi^{[-1]}\cdot \xi= I_{(\xi \neq 0)}$.

\begin{lemma}\label{lemma4.3}
	Let $G$ and $E$ be the same as in Definition \ref{definition4.1} and $f: G \rightarrow \bar{L}^0 (\mathcal{F})$ a proper $L^0$--convex function. Then we have the following:
	\begin{enumerate}[(1)]
		\item $f(\sum^n_{i=1} \lambda_i x_i) \leq \sum^n_{i=1} \lambda_i f(x_i)$ for any $n\in \mathbb{N}$, any $x_1, x_2, \cdots, x_n\in G$ and any $\lambda_1, \lambda_2, \cdots, \lambda_n\in L^0(\mathcal{F}, [0,1])$ such that $\sum^n_{i=1} \lambda_i =1$.
		\item $f$ is stable, namely $f(\sum^n_{i=1} \tilde{I}_{A_i} x_i) = \sum^n_{i=1} \tilde{I}_{A_i} f(x_i)$ for any $n\in \mathbb{N}$, any $x_1, x_2,$ $ \cdots, x_n\in G$ and any $n$--partition $\{A_1, A_2, \cdots, A_n \}$ of $\Omega$ to $\mathcal{F}$. In addition, if $G$ is $\sigma$--stable, then $f$ is also $\sigma$--stable, namely $f(\sum^{\infty}_{i=1} \tilde{I}_{A_i} x_i) = \sum^{\infty}_{i=1} \tilde{I}_{A_i} f(x_i)$ for any $\{A_i, i\in \mathbb{N} \}$ in $\Pi_{\mathcal{F}}$ and any sequence $\{x_i, i\in \mathbb{N}\}$ in $G$.
		\item If $f$ is $L^0$--affine, then $f(\sum^n_{i=1} \lambda_i x_i) = \sum^n_{i=1} \lambda_i f(x_i)$ for any $n\in \mathbb{N}$, any $x_1, x_2, \cdots, x_n$ in $G$ and any $\lambda_1, \lambda_2,\cdots, \lambda_n\in L^0(\mathcal{F}, [0,1])$ such that $\sum^n_{i=1} \lambda_i=1$.
	\end{enumerate}
\end{lemma}

\begin{proof}
	(1). It is obvious for $n=1$ or $2$. By induction, we assume (1) holds for $n=k\geq 2$. Since $\sum^{k+1}_{i=1} \lambda_i = 1$, then $I_{(1-\lambda_1\neq 0)} \lambda_i = \lambda_i$ for each $i\geq 2$, implying that $\sum^{k+1}_{i=1} \lambda_i x_i = \lambda_1 x_1+ I_{(1-\lambda_1\neq 0)} \lambda_2 x_2+ I_{(1-\lambda_1\neq 0)} \sum^{k+1}_{i=3} \lambda_i x_i$. Further, by the fact that $I_{(1-\lambda_1\neq 0)}= (1- \lambda_1) (1-\lambda_1)^{[-1]}$, we have the following:
	\begin{align*}
		& \sum^{k+1}_{i=1} \lambda_i x_i \\
		& =\lambda_1 x_1 + (1-\lambda_1) (1-\lambda_1)^{[-1]} \lambda_2 x_2+ (1-\lambda_1) \sum^{k+1}_{i=3} (1-\lambda_1)^{[-1]} \lambda_i x_i \\
		& = \lambda_1 x_1 + (1-\lambda_1) I_{(1-\lambda_1=0)} x_2 +  (1-\lambda_1) (1-\lambda_1)^{[-1]} \lambda_2 x_2\\
		&~~\ + (1-\lambda_1) \sum^{k+1}_{i=3} (1-\lambda_1)^{[-1]} \lambda_i x_i ~\text{(by noting}~(1-\lambda_1) I_{(1-\lambda_1=0)}=0)\\
		& = \lambda_1 x_1+ (1-\lambda_1) ((I_{(1-\lambda_1=0)}+ (1-\lambda_1)^{[-1]} \lambda_2 )x_2+ \sum^{k+1}_{i=3}(1-\lambda_1)^{[-1]} \lambda_i x_i).
	\end{align*}
	Thus, by induction assumption we have the following:
	\begin{align*}
		& f(\sum^{k+1}_{i=1} \lambda_i x_i) \\
		& \leq \lambda_1 f(x_1) + (1-\lambda_1) f((I_{(1-\lambda_1=0)}+ (1-\lambda_1)^{[-1]}\lambda_2 )x_2+ \sum^{k+1}_{i=3} (1-\lambda_1)^{[-1]} \lambda_i x_i ) \\
		& \leq \lambda_1 f(x_1) + (1-\lambda_1) ((I_{(1-\lambda_1=0)}+ (1-\lambda_1)^{[-1]}\lambda_2)f(x_2)+ \sum^{k+1}_{i=3} (1-\lambda_1)^{[-1]} \lambda_i f(x_i) )\\
		&~~\ \text{(by noting that}~ I_{(1-\lambda_1=0)}+ (1-\lambda_1)^{[-1]}\lambda_2 + \sum^{k+1}_{i=3} (1-\lambda_1)^{[-1]} \lambda_i = I_{(1-\lambda_1=0)}\\
		&~~\ + I_{(1-\lambda_1 \neq 0)} =1)\\
		& = \lambda_1 f(x_1) + I_{(1-\lambda_1 \neq 0)} \lambda_2 f(x_2) + \sum^{k+1}_{i=3} I_{(1-\lambda_1 \neq 0)} \lambda_i f(x_i) \\
		& = \lambda_1 f(x_1) +\lambda_2 f(x_2) + \cdots + \lambda_{k+1} f(x_{k+1}) (\text{again by}~I_{(1-\lambda_1 \neq 0)} \lambda_i =\lambda_i~\text{for each}~i\geq 2).
	\end{align*}
	
	\par
	(2). For any given $x_0\in G$, consider the function $f_1: G-x_0: =\{x-x_0: x\in G \} \rightarrow \bar{L}^0(\mathcal{F})$ defined by $f_1(x-x_0)= f(x)$ for any $x\in G$. It is easy to see that $f_1$ is a proper $L^0$--convex function and $f$ is stable ($\sigma$--stable) iff $f_1$ is stable ($\sigma$--stable). Therefore, we can, without loss of generality, assume that $G$ contains the null element $\theta$. We have $\tilde{I}_A x = \tilde{I}_A x+ \tilde{I}_{A^c} \theta$ still belongs to $G$ for any $x\in G$ and any $A\in \mathcal{F}$. Since $f(\tilde{I}_A x)= f(\tilde{I}_A x+ \tilde{I}_{A^c} \theta)\leq \tilde{I}_A f(x)+ \tilde{I}_{A^c} f(\theta)$, so $\tilde{I}_A f(\tilde{I}_A x)\leq \tilde{I}_A f(x)$. On the other hand, $\tilde{I}_A f(x) = \tilde{I}_A f(\tilde{I}_Ax+ \tilde{I}_{A^c} x)= \tilde{I}_A f(\tilde{I}_A \tilde{I}_A x+ \tilde{I}_{A^c} \tilde{I}_{A^c} x)\leq \tilde{I}_A f(\tilde{I}_A x).$ Thus, $\tilde{I}_A f(x) = \tilde{I}_A f(\tilde{I}_A x)$ for any $x\in G$ and any $A\in \mathcal{F}$. We have
	\begin{align*}
		& f(\sum^n_{i=1} \tilde{I}_{A_i} x_i) \\
		& = (\sum^n_{i=1} \tilde{I}_{A_i}) f(\sum^n_{i=1} \tilde{I}_{A_i} x_i) \\
		& = \sum^n_{i=1} \tilde{I}_{A_i}  f(\sum^n_{i=1} \tilde{I}_{A_i} x_i) \\
		& = \sum^n_{i=1} \tilde{I}_{A_i}  f(\tilde{I}_{A_i} x_i) \\
		& = \sum^n_{i=1} \tilde{I}_{A_i} f(x_i).
	\end{align*}
	
	\par
	If $G$ is $\sigma$--stable, for any $\{ A_i, i\in \mathbb{N}\}$ in $\Pi_{\mathcal{F}}$ and any sequence $\{ x_i, i\in \mathbb{N}\}$ in $G$, we similarly have:
	\begin{align*}
		& f(\sum^{\infty}_{i=1} \tilde{I}_{A_i} x_i) \\
		& = (\sum^{\infty}_{i=1} \tilde{I}_{A_i}) f(\sum^{\infty}_{i=1} \tilde{I}_{A_i} x_i) \\
		& = \sum^{\infty}_{i=1} \tilde{I}_{A_i}  f(\sum^{\infty}_{i=1} \tilde{I}_{A_i} x_i) \\
		& = \sum^{\infty}_{i=1} \tilde{I}_{A_i}  f(\tilde{I}_{A_i} x_i) \\
		& = \sum^{\infty}_{i=1} \tilde{I}_{A_i} f(x_i).
	\end{align*}
	
	\par
	(3). If $f$ is $L^0$--affine, then $f$ is both proper $L^0$--convex and proper $L^0$--concave. So we have $f(\sum^n_{i=1} \lambda_i x_i ) \leq \sum^n_{i=1} \lambda_i f(x_i)$ and $f(\sum^n_{i=1} \lambda_i x_i) \geq \sum^n_{i=1} \lambda_i f(x_i)$, namely, (3) holds.
\end{proof}

\par
To convey the gist of the proof of Theorem \ref{theorem1.19}, we give Lemma \ref{lemma4.4} below, which shows the simplest and most fundamental fact on $L^0$--extreme points, as well as motivates a systematic investigation into the Krein--Milman theorem in random locally convex modules \cite{GWT23}. For this, let us first recall the notion of an $L^0$--extreme point: let $E$ be an $L^0(\mathcal{F}, \mathbb{K})$--module, $G$ and $H$ two nonempty subsets of $E$ such that $H \subset G$. $H$ is called an $L^0$--extreme set of $G$ if both $x$ and $y$ belong to $H$ whenever $x$ and $y$ belong to $G$ and there exists $\lambda\in L^0(\mathcal{F}, (0,1))$ such that $\lambda x+(1-\lambda)y\in H$. Further, if $H=\{h\}$ is a singleton and $L^0$--extreme set of $G$, then $h$ is called an $L^0$--extreme point of $G$. We always denote by $ext_{L^0} (G)$ the set of $L^0$--extreme points of $G$.

\begin{lemma}\label{lemma4.4}
	Let $G$ be a nonempty subset of an $L^0(\mathcal{F}, \mathbb{K})$--module $E$, then we have the following statements:
	\begin{enumerate}[(1)]
		\item If $G$ is $\sigma$--stable and $ext_{L^0}(G) \neq \emptyset$, then $ext_{L^0} (G)$ is also $\sigma$--stable.
		\item If $G= Conv_{L^0} (\{x_1,x_2,\cdots, x_n \})$ is an $(n-1)$--dimensional $L^0$--simplex, then $ext_{L^0}(G)= \sigma(\{x_1,x_2, \cdots, x_n \})$ (Similar to Definition \ref{definition1.16}, the notion of an $L^0$--simplex can also be introduced in an $L^0(\mathcal{F}, \mathbb{K})$--module.)
	\end{enumerate}
\end{lemma}

\begin{proof}
	\par
	(1). For any $\{A_n, n\in \mathbb{N} \}$ in $\Pi_{\mathcal{F}}$ and any sequence $\{x_n, n\in \mathbb{N} \}$ in $ext_{L^0} (G)$, we will prove $\sum_n \tilde{I}_{A_n} x_n\in ext_{L^0} (G)$ as follows.
	
	\par
	Assume $\sum_n \tilde{I}_{A_n} x_n = \xi x+(1-\xi) y$ for some $x$ and $y$ in $G$ and for some $\xi \in L^0(\mathcal{F}, (0,1))$, we want to prove that $x=y$. For this, we only need to prove $\tilde{I}_{A_n} x= \tilde{I}_{A_n} y$ for each $n\in \mathbb{N}$. Otherwise, there exists some $n_0\in \mathbb{N}$ such that $\tilde{I}_{A_{n_0}} x \neq \tilde{I}_{A_{n_0}} y$. Now, let $\hat{x} = \tilde{I}_{A_{n_0}} x+ \tilde{I}_{A^c_{n_0}} x_{n_0}$ and $\hat{y} = \tilde{I}_{A_{n_0}} y+ \tilde{I}_{A^c_{n_0}} x_{n_0}$. Then $\hat{x} \neq \hat{y}$, and $\hat{x}$ and $\hat{y}$ belong to $G$. Since $\tilde{I}_{A_n} x_n = \tilde{I}_{A_n} (\sum^{\infty}_{k=1} \tilde{I}_{A_k} x_k)= \tilde{I}_{A_n} (\xi x+ (1-\xi)y) = \tilde{I}_{A_n} \xi x + \tilde{I}_{A_n} (1-\xi) y$, it immediately follows from the definition of $\hat{x}$ and $\hat{y}$ that $\tilde{I}_{A_{n_0}} x_{n_0} = \tilde{I}_{A_{n_0}} (\xi \hat{x} + (1-\xi) \hat{y})$ and $\tilde{I}_{A_n} x_{n_0} = \tilde{I}_{A_n} \hat{x} = \tilde{I}_{A_n} \hat{y}$ for any $n \neq n_0$. By noting that $\xi+ (1-\xi)=1$, we also have that $\tilde{I}_{A_n} x_{n_0} = \tilde{I}_{A_n} (\xi \hat{x} + (1-\xi) \hat{y})$ for any $n \neq n_0$. To sum up, $x_{n_0}= \xi \hat{x} + (1-\xi) \hat{y}$, which contradicts the fact of $x_{n_0}\in ext _{L^0} (G)$. So we have competed the proof of (1).
	
	\par
	(2). It is obvious for $n=1$. When $n> 2$, it is clear that $\{x_1, x_2, \cdots, x_n \}\subset ext_{L^0}(G)$. By (1) and the fact that an $L^0$--simplex is always $\sigma$--stable, we can have $\sigma(\{x_1, x_2, \cdots, x_n \} ) \subset ext_{L^0} (G)$.
	
	\par
	Conversely, if there exists $x\in ext_{L^0} (G) \backslash \sigma(\{x_1, x_2, \cdots, x_n \} )$, let $x= \sum^n_{i=1} \xi_i x_i$ for some $\xi_1, \xi_2, \cdots$ and $\xi_n\in L^0(\mathcal{F}, [0,1])$ such that $\sum^n_{i=1} \xi_i =1$. Then there must be some $i_0\in \{1,2,\cdots,n \}$ such that $P\{\omega\in \Omega: 0< \xi^0_{i_0} (\omega)<1 \}>0$, where $\xi^0_{i_0}$ is an arbitrarily chosen representative of $\xi_{i_0}$. We can, without loss of generality, assume that $i_0=1$, and further let $A=\{\omega\in \Omega: 0< \xi^0_1(\omega)<1 \}$. Then it is obvious that $\sum^n_{i=2} \xi_i = 1-\xi_1$ satisfies $0< 1-\xi_1<1$ on $A$, and it is also clear that $\tilde{I}_A x= \tilde{I}_A\xi_1 x_1+ \tilde{I}_A (1-\xi_1) ((1-\xi_1)^{[-1]} \sum^n_{i=2} \xi_i x_i )$ and $\tilde{I}_A x_1 \neq \tilde{I}_A (1-\xi_1)^{[-1]} \sum^n_{i=2} \xi_i x_i$. (Otherwise, $\tilde{I}_A$ must equal to $0$, which contradicts $P(A)>0$). Now, we define $\xi$ and $\eta$ in $L^0(\mathcal{F}, (0,1))$ and $\hat{x}$ and $\hat{y}$ in $G$ as follows:
	
	\par
	\noindent
	$\xi = \tilde{I}_A \xi_1 + \frac{1}{2} \tilde{I}_{A^c}$,
	
	\par
	\noindent
	$\eta = \tilde{I}_A (1-\xi_1) + \frac{1}{2} \tilde{I}_{A^c}$,
	\par
	\noindent
	$\hat{x} = \tilde{I}_A x_1 + \tilde{I}_{A^c} x,$
	\par
	\noindent
	$\hat{y} = \tilde{I}_A (1-\xi_1)^{[-1]} \sum^n_{i=2} \xi_i x_i + \tilde{I}_{A^c} x$.\\
	Then $\xi+ \eta =1$ and $x= \xi \hat{x} + \eta \hat{y}$ with $\hat{x} \neq \hat{y}$, which contradicts the fact that $x\in ext_{L^0}(G)$.
\end{proof}

\par
Lemma \ref{lemma4.5} characterizes the topological structure of an $L^0$--simplex in an $RN$ module (in fact, Lemma \ref{lemma4.5} remains valid for a random locally convex module from the arguments of its proof).

\begin{lemma}\label{lemma4.5}
	Let $(E, \|\cdot\|)$ be an $RN$ module over $\mathbb{K}$ with base $(\Omega, \mathcal{F}, P)$ and $G= Conv_{L^0} (\{x_1,x_2,\cdots, x_n \})$ an $(n-1)$--dimensional $L^0$--simplex. Further, let $\Delta _{n-1} = Conv_{L^0} (\{ e_1, e_2, \cdots , e_n\} )$ be the standard $(n-1)$--dimensional $L^0$--simplex in $L^0(\mathcal{F}, \mathbb{R}^n)$. Define the mapping $T: G \rightarrow \Delta _{n-1}$ by $T(x)= (\lambda_1, \lambda_2,\cdots, \lambda_n)$ for any $x= \sum^n_{i=1} \lambda_i x_i\in G$, where $\{e_i: i=1\sim n \}$ is the orthonormal basis for $\mathbb{R}^n$. Then $T$ is an $L^0$--affine $\mathcal{T}_{\varepsilon,\lambda}$--isomorphism form $G$ onto $\Delta _{n-1}$.
\end{lemma}

\begin{proof}
	It is obvious when $n=1$, we only need to prove the case for $n\geq 2$. Since $\{x_i-x_1: i=2\sim n \}$ is $L^0$--independent in $E$ and $\{e_i-e_1: i=2\sim n \}$ is $L^0$--independent in $L^0(\mathcal{F}, \mathbb{K}^n)$, by Lemma 3.4 of \cite{GP} both $M_1 = \{\sum^n _{i=2} \xi_i (x_i- x_1): \xi_i\in  L^0(\mathcal{F}, \mathbb{K}) $ for each $i=2\sim n\}$ as an $L^0(\mathcal{F}, \mathbb{K})$--free submodule in $E$ of $(n-1)$ rank and $M_2 = \{\sum^n _{i=2} \xi_i (e_i- e_1): \xi_i\in  L^0(\mathcal{F}, \mathbb{K}) $ for each $i=2\sim n\}$ as an $L^0(\mathcal{F}, \mathbb{K})$--free submodule in $L^0(\mathcal{F}, \mathbb{K}^n)$ of $(n-1)$ rank are $\mathcal{T}_{\varepsilon,\lambda}$--isomorphic onto
	$L^0(\mathcal{F}, \mathbb{K}^{n-1})$ in the sense of topological modules. So $M_1$ and $M_2$ are also $\mathcal{T}_{\varepsilon,\lambda}$--isomorphic under the isomorphism $T_1$ defined by $T_1(\sum^n_{i=2} \xi_i (x_i-x_1) ) = \sum^n_{i=2} \xi_i (e_i-e_1)$ for any $\xi_2, \xi_3, \cdots$ and $\xi_n$ in $L^0(\mathcal{F}, \mathbb{K})$. Further, let $S_{n-1} = \{(\lambda_2, \lambda_3, \cdots, \lambda_n )\in L^0(\mathcal{F}, [0,1])^{n-1}: \sum^n_{i=2} \lambda_i \leq 1 \}$, and denote by $T_1'$ the restriction of $T_1$ to $G':= \{\sum^n_{i=2} \lambda_i (x_i-x_1): (\lambda_2, \lambda_3,\cdots, \lambda_n)\in S_{n-1} \}$. It is obvious that $T_1'$ is an $L^0$--affine $\mathcal{T}_{\varepsilon,\lambda}$--isomorphism from $G_1'$ onto $\Delta'_{n-1}:= \{ \sum^n_{i=2} \lambda_i (e_i-e_1): (\lambda_2, \lambda_3,\cdots, \lambda_n)\in S_{n-1} \}$. Finally, by observing $G= x+ G'$ and $\Delta_{n-1} = e_1+ \Delta'_{n-1}$, further from $T(\sum^n_{i=1} \lambda_i x_i ) = (\lambda_1, \lambda_2, \cdots, \lambda_n) = \lambda_1 e_1+ \lambda_2 e_2 + \cdots + \lambda_n e_n = e_1 + \sum^n_{i=2} \lambda_i (e_i-e_1) = T(x_1) + T_1'(\sum^n_{i=2} \lambda_i (x_i-x_1) )$, we see that $T$ is also an $L^0$--affine $\mathcal{T}_{\varepsilon,\lambda}$--isomorphism from $G$ onto $\Delta_{n-1}$.
\end{proof}

\begin{remark}\label{remark4.6}
	Lemma \ref{lemma4.5} means both $T$ and $T^{-1}$ are $L^0$--affine and $\mathcal{T}_{\varepsilon,\lambda}$--continuous. So by Lemma \ref{lemma3.3}, both $T$ and $T^{-1}$ are $\mathcal{T}_c$--continuous, namely, $T$ is also $L^0$--affinely $\mathcal{T}_c$--isomorphic. Further, Lemma \ref{lemma4.5} also means that the sequence $\{y_k = \sum^n_{i=1} \lambda^k_i x_i,$ $ k\in \mathbb{N} \}$ in $G$ converges in $\mathcal{T}_{\varepsilon,\lambda}$ to $y_0 = \sum^n_{i=1} \lambda^0_i x_i$ in $G$ iff $\{\lambda^k_i, k\in \mathbb{N} \}$ converges in $\mathcal{T}_{\varepsilon,\lambda}$ (namely, in probability $P$) to $\lambda^0_i$ for each $i=1\sim n$.
\end{remark}

\par
We are now in a position to prove Theorem \ref{theorem1.19}. First, recall the hypothesis of a $n$--person conditional game: let the $RN$ module $(E_i, \|\cdot\|_i)$ over $\mathbb{K}$ with base $(\Omega, \mathcal{F}, P)$, the $(n_i-1)$--dimensional $L^0$--simplex $S_i = Conv_{L^0} (\{\pi_{i\alpha}: \alpha= 1\sim n_i \})$ and the $n$--fold $L^0$--affine function $p_i: S:= \Pi^n_{i=1} S_i \rightarrow L^0(\mathcal{F})$, respectively, denote the space of random payoffs, the set of mixed strategies and the $L^0$--payoff function corresponding to the player $i$. Similar to the proof of Lemma \ref{lemma3.4}, one can prove that $S$ is a random sequentially compact $L^0$--convex subset of the $RN$ module $(E, \|\cdot\|)$ since each $S_i$ is random sequentially compact and $L^0$--convex, where $(E, \|\cdot\|)$ is the product space of $\{(E_i, \|\cdot\|_i): i=1\sim n \}$. Further, by Lemma \ref{lemma4.5} each $p_i$ is $\mathcal{T}_{\varepsilon,\lambda}$--continuous.

\par
Even though our proof of Theorem \ref{theorem1.19} follows the construction and method used in the proof of Theorem 1 of \cite{Nash}, our proof of Theorem \ref{theorem1.19} is based on our Theorem \ref{theorem1.12} or Lemma \ref{lemma3.7}, whereas the proof of Theorem 1 of \cite{Nash} is based on the classical Brouwer fixed point theorem. Besides, each $s_i\in S_i$ can be uniquely represented as $s_i = \sum^{n_i}_{\alpha=1} c_{i\alpha} \pi_{i\alpha}$, the coefficients $c_{i\alpha}$ in the $L^0$--convex combination are elements in $L^0(\mathcal{F}, [0,1] )$. Then there exist rather complicated measurability problems to be treated in our proof of Theorem \ref{theorem1.19}, which makes our proof of Theorem \ref{theorem1.19} somewhat more involved than Theorem 1 of \cite{Nash}.

\begin{proof}[\textbf{Proof of Theorem \ref{theorem1.19}}]
	Given a point $s= (s_1, s_2, \cdots, s_n)$ in $S: = \Pi^n_{i=1} S_i$, since each $s_i$ has a unique representation $s_i = \sum^{n_i} _{\alpha =1} c_{i\alpha } \pi_{i\alpha}$ with $c_{i\alpha } \in L^0(\mathcal{F}, [0,1])$ such that $\sum^{n_i}_{\alpha=1} c_{i\alpha} = 1$, then by (3) of Lemma \ref{lemma4.3} we have the following representation for $p_i(s)$:\\
	%\begin{align*}
	%  p_i & = \sum^{n_i}_{\alpha=1} c_{i\alpha} p_i(s_1,s_2,\cdots, s_{i-1}, \pi_{i\alpha}, s_{i+1}, \cdots, s_n ) \notag\\
	%      & = \sum^{n_i}_{\alpha=1} c_{i\alpha} p_{i\alpha}(s),
	%\end{align*}
	$p_i  = \sum^{n_i}_{\alpha=1} c_{i\alpha} p_i(s_1,s_2,\cdots, s_{i-1}, \pi_{i\alpha}, s_{i+1}, \cdots, s_n )= \sum^{n_i}_{\alpha=1} c_{i\alpha} p_{i\alpha}(s)${\hfill (1) \\}
	where, $p_{i\alpha}(s) = p_i(s_1,s_2,\cdots, s_{i-1}, \pi_{i\alpha}, s_{i+1}, \cdots, s_n )$.
	
	\par
	It immediately follows from (1) that $s$ is a Nash equilibrium point iff $p_i(s)= \max_{1\leq \alpha \leq n_i} p_{i\alpha}(s)$.{\hfill (2) \par}
	
	\par
	Now, we define the mapping $T: S\rightarrow S$ as follows:\\
	$s_i'=\frac{s_i+ \sum^{n_i}_{\alpha=1} \varphi_{i\alpha} (s) \pi_{i\alpha}}{1+ \sum^{n_i}_{\alpha=1} \varphi_{i\alpha} (s)}$ for each $i=1\sim n$, {\hfill (3) \par}
	\noindent
	where, we use $(s_1', s_2', \cdots, s'_n)$ for $T(s)$ for each $s\in S$, and $\varphi_{i\alpha} (s)$ for $\max(0, p_{i\alpha}(s)-p_i(s))$ for each $i=1\sim n$ and each $\alpha=1\sim n_i$.
	
	\par
	It is obvious from (2) that if $s$ is a Nash equilibrium point, then $s$ must be a fixed point of $T$. It remains to prove that a fixed point $s$ of $T$ is also a Nash equilibrium point, which is just the complicated part of the proof of the theorem.
	
	\par
	If $s$ is a fixed point of $T$, then $s'_i=s_i$ for each $i=1\sim n$ and we obtain the following relations:\\
	$c_{i\alpha} = \frac{c_{i\alpha}+ \varphi_{i\alpha} }{1+ \sum^{n_i}_{\alpha=1} \varphi_{i\alpha} }$ for each $i=1\sim n$ and each $\alpha=1\sim n_i$.{\hfill (4) \par}
	
	According to (1) and (4), we can choose a representative $p_i^0(s)$ of $p_i(s)$ for each $i=1\sim n$, a representative $c^0_{i\alpha}$ of $c_{i\alpha}$ and a representative $p^0_{i\alpha}(s)$ of $p_{i\alpha}(s)$ for each $i=1\sim n$ and each $\alpha=1\sim n_i$ such that the following relations (5) and (6)  hold for some $\Omega'\in \mathcal{F}$ with $P(\Omega')=1$:\\
	$p^0_i(s) = \sum^{n_i}_{\alpha=1} c^0_{i\alpha}(\omega) p^0_{i\alpha}(s)(\omega)$ for each $i=1\sim n$ and each $\omega\in \Omega'$. {\hfill (5) \\}
	$c^0_{i\alpha}(\omega) = \frac{c^0_{i\alpha}(\omega) + \varphi^0_{i\alpha}(s)(\omega) }{1+ \sum^{n_i}_{\alpha=1}\varphi^0_{i\alpha}(s)(\omega) }$ for each $i=1\sim n$, each $\alpha=1\sim n_i$ and each $\omega\in \Omega'$. {\hfill (6) \\}
	Here $\varphi^0_{i\alpha}(s) = \max(0, p^0_{i\alpha}(s)-p^0_i(s))$.
	\par
	We can, without loss of generality, assume $\Omega'=\Omega$. Once we have chosen these representatives as above, then for an arbitrary but fixed $i\in \{1,2,\cdots, n \}$, we first define a $\mathcal{F}$--measurable random variable $\tau: \Omega\rightarrow \{1,2,\cdots, n_i\}$ as follows:\\
	$\tau = \sum^{n_i}_{\alpha=1} I_{(c^0_{i\alpha}>0)}$, where $I_{(c^0_{i\alpha}>0)}$ stands for the characteristic function of $\{\omega\in \Omega: c^0_{i\alpha}(\omega)>0 \}$ for each $\alpha=1\sim n_i$, namely, for each $\omega\in \Omega$ $\tau(\omega)$ is the cardinal number of $\{\alpha\in \{1,2,\cdots, n_i \}: c_{i\alpha}(\omega)>0 \}$.
	
	\par
	We can, without loss of generality, assume that $P(A_k)>0$ for each $k=1\sim n_i$, where $A_k=\{\omega\in \Omega: \tau(\omega)=k \}$ (otherwise, we omit this $A_k$). For each $k\in \{ 1,2,\cdots, n_i\}$, first let $\tau_{k,0}: A_k \rightarrow \mathbb{R}$ be the function defined by $\tau_{k,0}(\omega) = 0$ for each $\omega\in A_k$. Then we can recursively define $A_k\cap \mathcal{F}$--measurable functions $\tau_{k,1}, \tau_{k,2},\cdots $ and $\tau_{k,k}$ as follows:\\
	$\tau_{k,m}: A_k \rightarrow \{1,2,\cdots, n_i \}$ by $\tau_{k,m}(\omega) = \min\{l\in \{1,2,\cdots,n_i\}: \tau_{k,m-1}(\omega)<l\leq n_i ~\text{and}~c^0_{il}(\omega)>0\}$ for each $m\in \{1,2,\cdots,k\}$ and each $\omega\in A_k$, where $A_k \cap \mathcal{F}:= \{A_k\cap A: A\in  \mathcal{F}\}$.
	
	\par
	It is clear that $\tau_{k,1}(\omega)< \tau_{k,2}(\omega)< \cdots < \tau_{k,m-1}(\omega) < \tau_{k,m}(\omega)< \cdots < \tau_{k,k}(\omega)$ for each $\omega\in A_k$, and (5) becomes the following form:\\
	$p^0_i(s)(\omega) = \sum^k_{m=1} c^0_{i\tau_{k,m}(\omega)}(\omega) p^0_{i\tau_{k,m}(\omega) }(s)(\omega)$ for each $\omega\in A_k$.{\hfill (7) \par}
	Now, we define an $\mathcal{F}$--measurable random variable $\sigma: \Omega= \sum^{n_i}_{k=1} A_k \rightarrow \{1,2,\cdots, $ $n_i \}$ as follows:\\
	$\sigma(\omega)= \min\{\tau_{k,m}(\omega): m=1,2,\cdots,k$ and $p^0_{i \tau_{k,m}(\omega)}(s) (\omega) =\min\{p^0_{i \tau_{k,l}(\omega)}(s) (\omega): l=1,2,\cdots,k \} \}$ when $\omega\in A_k$ for some $k\in \{1,2,\cdots,n_i\}$.
	
	\par
	It is obvious that $c^0_{i\sigma (\omega)} (\omega) \neq 0$ and $p^0_i(s) (\omega) \geq p^0_{i \sigma(\omega)} (s) (\omega)$ for each $\omega\in \Omega$. Further, let $B_q= \{\omega\in \Omega: \sigma(\omega)=q \}$ for each $q\in \{1,2,\cdots, n_i \}$.Then $B_q\in \mathcal{F}$ and $\sum^{n_i}_{q=1} B_q = \Omega$. We can, without loss of generality, again assume that $P(B_q)>0$ for each $q=1\sim n_i$ (otherwise, we can omit this $B_q$). Now, by (6) we have the following:\\
	$c^0_{iq}(\omega) = \frac{c^0_{iq}(\omega)}{1+ \sum^{n_i}_{i=1} \varphi^0_{i\alpha }(s)(\omega) }$ for each $\omega\in B_q$ since $\varphi^0_{iq} (s) (\omega) = 0$ for each $\omega\in B_q$.\\
	Thus $\sum^{n_i}_{i=1} \varphi^0_{i\alpha }(s)(\omega) =0$ for each $\omega\in B_q$ since $c^0_{iq}(\omega) \neq 0$ for each $\omega\in B_q$, which also means $\sum^{n_i}_{i=1} \varphi^0_{i\alpha }(s)(\omega) =0$ for each $\omega\in \Omega$ since $\sum^{n_i}_{q=1} B_q = \Omega$, namely, $\sum^{n_i}_{i=1} \varphi_{i\alpha }(s) =0$, yielding $\varphi_{i\alpha} (s)=0$ for each $\alpha =1\sim n_i$ (that is to say, $p_{i\alpha} (s) \leq p_i(s)$ for each $\alpha = 1\sim n_i$).
	
	\par
	Since $i$ is arbitrary, we can get that $p_i(s) = \max\{p_{i\alpha}(s): \alpha=1\sim n_i \}$ for each $i=1\sim n$. Then by (2) $s$ is a Nash equilibrium point. Similar to the proof of (a) of Lemma \ref{lemma3.8} one can check that $T$ is $\sigma$--stable. Further, since $T$ is also $\mathcal{T}_{\varepsilon,\lambda}$--continuous,  by Theorem \ref{theorem1.12} or Lemma \ref{lemma3.7} $T$ always has a fixed point. So every $n$--person conditional game in turn always has a Nash equilibrium point.
\end{proof}

\begin{remark}\label{remark4.7}
	In the proof of Theorem \ref{theorem1.19}, we have used the $n$--fold $L^0$--affine property of an $L^0$--payoff function to show that $s\in S$ is a Nash equilibrium point iff $p_i(s) = \max\{p_{i\alpha}(s): \alpha=1\sim n_i \}$. It is easy to see that there is an $n_i$--partition $\{A_1,A_2,\cdots, A_{n_i} \}$ of $\Omega$ to $\mathcal{F}$ such that $\max\{p_{i\alpha}(s): \alpha=1\sim n_i \} = p_i(s_1,s_2,\cdots, s_{i-1}, \sum^{n_i}_{\alpha=1} \tilde{I}_{A_{\alpha}} \pi_{i\alpha}, s_{i+1}, \cdots, s_n )$, namely, $p_i(s_1,s_2,\cdots, s_{i-1},\cdot, s_{i+1}, \cdots, $ $s_n )$ attains its maximum value at an $L^0$--extreme point $\sum^{n_i}_{\alpha=1} \tilde{I}_{A_{\alpha}} \pi_{i\alpha}$ of $S_i$. From (2) of Lemma \ref{lemma4.3}, one can easily see that an $L^0$--concave function defined on an $L^0$--convex and $\sigma$--stable set must be $\sigma$--stable (but the assertion does not necessarily hold for an $L^0$--quasiconcave function). If we only assume that each $L^0$--payoff function $p_i: S\rightarrow L^0(\mathcal{F})$ satisfies the following two conditions:
	\begin{enumerate}[(A)]
		\item Each $p_i$ is random sequentially continuous and $\sigma$--stable.
		\item For each given $s\in S$, the function $p_i(s_1,s_2,\cdots,s_{i-1}, \cdot, s_{i+1}, \cdots, s_n): S_i \rightarrow L^0(\mathcal{F})$ is $L^0$--quasiconcave.
	\end{enumerate}
	It is easy to check that if $p_i$ is $\sigma$--stable, then the function in (B) above is also $\sigma$--stable. At the present time, we do not know whether Theorem \ref{theorem1.19} still holds under both hypotheses (A) and (B).  We do not know even whether the function attains its maximum value at an $L^0$--extreme point of $S_i$.
\end{remark}

\par
When $(\Omega, \mathcal{F}, P)$ is trivial, namely $\mathcal{F} = \{\Omega, \emptyset \}$, then the conditional finite game reduces to the static finite game originally considered by Nash in \cite{Nash}. In such a situation, the static case (A) and (B) reduce to the following, respectively:
\begin{enumerate}[(A$'$)]
	\item Each $p_i$ is a continuous real--valued function.
	\item For each given $s\in S$, the function $p_i(s_1, s_2, \cdots, s_{i-1}, \cdot, s_{i+1}, \cdots, s_n): S_i \rightarrow \mathbb{R}$ is quasiconcave.
\end{enumerate}
It is well known from \cite{Bor,Fan,Mau} that Theorem 1 of \cite{Nash} remains valid under both hypotheses (A$'$) and (B$'$). Nevertheless, the method used to prove the general case is by means of Ky Fan's coincidence theorem \cite{Fan} rather than the Brouwer fixed point theorem. Since Ky Fan's coincidence theorem \cite{Fan} can be deduced from the Knaster--Kuratowski--Mazurkiewicz theorem (briefly, KKM theorem), this requires us to establish the noncompact KKM theorem and Ky Fan's coincidence theorem in $RN$ modules (or, more general random locally convex modules), just as we have established the noncompact Schauder fixed point theorem in Section \ref{section3} of this paper. Therefore, we can eventually solve the existence problem of the conditional Nash equilibrium points under (A) and (B). We are going to study those closely related problems in a forthcoming work.

\section{Concluding remarks}\label{section5}

The classical Tychonoff fixed point theorem \cite{Ty} states that a compact convex set of a Hausdorff locally convex space has the fixed point property. From there, one is naturally concerned with the following fundamental problem---Problem \ref{problem5.1}. The reader may refer to \cite{Guo3} for the notion of a random locally convex module $(E, \mathcal{P})$ over $\mathbb{K}$ with base $(\Omega, \mathcal{F}, P)$ together with the $(\varepsilon,\lambda)$--topology $\mathcal{T}_{\varepsilon,\lambda}$ and the locally $L^0$--convex topology $\mathcal{T}_c$ induced by the family $\mathcal{P}$ of $L^0$--seminorms. The central purpose of this section is to point out that Problem \ref{problem5.1} is still open and our Theorem \ref{theorem1.12} is just a partial answer to this problem.

\begin{problem}\label{problem5.1}
	Let $(E, \mathcal{P})$ be a random locally convex module over $\mathbb{K}$ with base $(\Omega, \mathcal{F}, P)$ and $G$ a stably compact $L^0$--convex subset of $E$. Does a $\sigma$--stable $\mathcal{T}_{\varepsilon,\lambda}$--continuous or $\mathcal{T}_c$--continuous mapping $T$ from $G$ to $G$ have a fixed point?
\end{problem}

\par
To help the reader understand Problem \ref{problem5.1}, let us first introduce the idea of stable compactness for a $\sigma$--stable set. Based on the work \cite{DJKK} on conditional compactness, Jamneshan and Zapata \cite{JZ} introduced  stable compactness as a kind of conditional compactness in the classical set--theoretic framework. Let $(E, \mathcal{P})$ be the same as in Problem \ref{problem5.1} and $G$ a  $\sigma$--stable subset of $E$. A nonempty family $\mathcal{E}$ of nonempty subsets of $G$ is said to be $\sigma$--stable if $\sum^{\infty}_{n=1} \tilde{I}_{A_n} G_n: = \{\sum^{\infty} _{n=1} \tilde{I}_{A_n} g_n : g_n\in G_n$ for each $n\in \mathbb{N}\}$ still lies in $\mathcal{E}$ for any $\{A_n, n\in \mathbb{N} \}$ in $\Pi_{\mathcal{F}}$ and any countable subfamily $\{G_n, n\in \mathbb{N} \}$ of $\mathcal{E}$. It is readily pointed out in \cite{GWT23} that a $\sigma$--stable subset $G$ of $E$ is stably compact in the sense of \cite{JZ} iff every $\sigma$--stable family $\mathcal{E}$ of $\sigma$--stable $\mathcal{T}_{\varepsilon,\lambda}$--closed subsets of $G$ has a nonempty intersection whenever $\mathcal{E}$ has the finite intersection property.

\par
To clarify the fact that a $\sigma$--stable set of an $RN$ module is random sequentially compact iff it is stably compact, we denote by $L^0(\mathcal{F}, \mathbb{N})$ the set of equivalence classes of random variables from $(\Omega, \mathcal{F}, P)$ to $\mathbb{N}$ (obviously, $L^0(\mathcal{F}, \mathbb{N})$ is a $\sigma$--stable subset of $L^0(\mathcal{F})$). Further recall the following.

\begin{definition}[\cite{JZ}]\label{definition5.2}
	Let $(E, \|\cdot\|)$ be an $RN$ module over $\mathbb{K}$ with base $(\Omega, \mathcal{F}, P)$ and $G$ a $\sigma$--stable subset of $E$. A $\sigma$--stable function $x: L^0(\mathcal{F}, \mathbb{N}) \rightarrow G$ is called a stable sequence in $G$, denoted by $\{x_n, n\in L^0(\mathcal{F}, \mathbb{N}) \}$. A stable sequence $\{x_n, n\in \mathbb{N} \}$ in $G$ is said to be $\mathcal{T}_c$--Cauchy if for each $\varepsilon\in L^0_{++}(\mathcal{F})$, there exists $l\in  L^0(\mathcal{F}, \mathbb{N})$ such that $\|x_n-x_m\|\leq \varepsilon$ whenever $n\geq l$ and $m\geq l$. Further $G$ is said to be $\mathcal{T}_c$--stably sequentially complete if every $\mathcal{T}_c$--Cauchy stable sequence in $G$ converges in $\mathcal{T}_c$ to a point in $G$. A stable sequence $\{y_m, m\in L^0(\mathcal{F}, \mathbb{N})\}$ in $G$ is called a stable subsequence of a stable sequence $\{x_n, n\in L^0(\mathcal{F}, \mathbb{N}) \}$ if $y_m =x_{n_m}$ for each $m\in L^0(\mathcal{F}, \mathbb{N})$ for some $\sigma$--stable function $n: L^0(\mathcal{F}, \mathbb{N}) \rightarrow L^0(\mathcal{F}, \mathbb{N})$ (denoted $n_m=n(m)$) such that $\{y_m, m\in L^0(\mathcal{F}, \mathbb{N})\}$ is also a subnet of $\{x_n, n\in L^0(\mathcal{F}, \mathbb{N}) \}$. Finally, $G$ is said to be stably sequentially compact if every stable sequence in $G$ has a stable subsequence in $G$ convergent in $\mathcal{T}_c$ to a point in $G$.
\end{definition}

\par
The locally $L^0$--convex topology $\mathcal{T}_c$ on an $RN$ module is not metrizable in general and a stable sequence is merely a special net. Thus, on the surface the notion of $\mathcal{T}_c$--stably sequential completeness defined as in Definition \ref{definition5.2} is not that of the standard $\mathcal{T}_c$--completeness. However, Theorem \ref{theorem5.3} below unifies the three types of completeness.

\begin{theorem}\label{theorem5.3}
	Let $(E, \|\cdot\|)$ be an $RN$ module over $\mathbb{K}$ with base $(\Omega, \mathcal{F}, P)$ and $G$ a $\sigma$--stable subset of $E$. Then the following are equivalent to each other:
	\begin{enumerate}[(1).]
		\item $G$ is $\mathcal{T}_{\varepsilon,\lambda}$--complete.
		\item $G$ is $\mathcal{T}_c$--complete.
		\item $G$ is $\mathcal{T}_c$--stably sequentially complete.
	\end{enumerate}
\end{theorem}

\begin{proof}
	(1) $\Leftrightarrow$ (2) is known from Remark \ref{remark1.7}. Since $\mathcal{T}_{\varepsilon,\lambda}$ is metrizable and the usual sequences can be conveniently connected with stable sequences. So we only need prove (1) $\Leftrightarrow$ (3) as follows.
	
	\par
	(1) $\Rightarrow$ (3). Let $\{ x_n, n\in L^0(\mathcal{F}, \mathbb{N})\}$ be a $\mathcal{T}_c$--Cauchy stable sequence in $G$, then for any given $\varepsilon\in L^0_{++} (\mathcal{F})$ there exists $n_{\varepsilon} \in L^0(\mathcal{F}, \mathbb{N})$ such that $\|x_m-x_n\|\leq \varepsilon$ for any $m$ and $n$ in $L^0(\mathcal{F}, \mathbb{N})$ satisfying $m,n\geq n_{\varepsilon}$. Since $\mathcal{T}_c$ is stronger than $\mathcal{T}_{\varepsilon,\lambda}$, $\{ x_n, n\in L^0(\mathcal{F}, \mathbb{N})\}$ is also a $\mathcal{T}_{\varepsilon,\lambda}$--Cauchy net. Thus, there exists $x\in G$ such that $\{ x_n, n\in L^0(\mathcal{F}, \mathbb{N})\}$ converges in $\mathcal{T}_{\varepsilon,\lambda}$ to $x$. Now, for any fixed $m\in L^0(\mathcal{F}, \mathbb{N})$ with $m\geq n_{\varepsilon}$, it is easy to see that $\{y\in G: \|x_m-y\|\leq \varepsilon \}$ is $\mathcal{T}_{\varepsilon,\lambda}$--closed since $G$ is $\mathcal{T}_{\varepsilon,\lambda}$--closed according to the $\mathcal{T}_{\varepsilon,\lambda}$--completeness of $G$, then $\|x_m-x\| \leq \varepsilon$, namely, (3) holds.
	
	\par
	(3) $\Rightarrow$ (1). Since $(G, d)$ is a random metric space with base $(\Omega, \mathcal{F}, P)$, where $d(x,y)= \|x-y\|$ for any $x$ and $y$ in $G$, as in the case of an ordinary metric space, $(G, d)$ has a completion $(\widetilde{G}_{\varepsilon,\lambda}, \widetilde{d})$ with respect to the $(\varepsilon,\lambda)$--uniform structure (so we can still denote $\widetilde{d}$ by $d$), see \cite{GWYZ} for the $(\varepsilon,\lambda)$--uniform structure induced by a random metric. We will prove $G$ is $\mathcal{T}_{\varepsilon,\lambda}$--complete by verifying $G= \widetilde{G}_{\varepsilon,\lambda}$ as follows. Let $\tilde{x}$ be any given element in
	$\widetilde{G}_{\varepsilon,\lambda}$. Then there exists some $\mathcal{T}_{\varepsilon,\lambda}$--Cauchy sequence $\{x_k, k\in \mathbb{N} \}$ in $G$ such that $\{\|x_k- \tilde{x}\|, k\in \mathbb{N} \}$ converges in probability to $0$, we can, without loss of generality, assume that $\{\|x_k- \tilde{x}\|, k\in \mathbb{N} \}$ converges a.s. to $0$. Given an $n\in L^0(\mathcal{F}, \mathbb{N})$ written as $n = \sum^{\infty}_{k=1} \tilde{I}_{A_k} \cdot k$ for some $\{A_k, k\in \mathbb{N} \}$ in $\Pi_{\mathcal{F}}$, let $x_n = \sum^{\infty}_{k=1} \tilde{I}_{A_k} \cdot x_k$. Then it is easy to see that $\{x_n, n\in  L^0(\mathcal{F}, \mathbb{N})\}$ is a stable sequence in $G$. We will prove that $\{\|x_n- \tilde{x}\|, n\in L^0(\mathcal{F}, \mathbb{N}) \}$ converges in $\mathcal{T}_c$ to $0$ as follows, namely, for any given $\varepsilon \in L^0_{++}(\mathcal{F})$ there exists $n_{\varepsilon} \in L^0(\mathcal{F}, \mathbb{N})$ such that $\|x_n-\tilde{x}\|\leq \varepsilon$ for any $n\in L^0(\mathcal{F}, \mathbb{N})$ with $n\geq n_{\varepsilon}$.
	
	\par
	In fact, since $\varepsilon \in L^0_{++}(\mathcal{F})$ and $\mathbb{N}$ is countable, we can choose a representative $\|x_k-\tilde{x}\|^0$ of $\|x_k-\tilde{x}\|$ for each $k \in \mathbb{N}$ and a representative $\varepsilon^0$ of $\varepsilon$ such that $\varepsilon^0(\omega)>0$ and $\{\|x_k-\tilde{x}\|^0(\omega), k\in \mathbb{N} \}$ converges to $0$ for each $\omega\in \Omega$. Now, define $n^0_{\varepsilon}: \Omega \rightarrow \mathbb{N}$ as follows:\\
	$n^0_{\varepsilon}(\omega)= \min\{l\in \mathbb{N}: \sup\{\|x_k-\tilde{x}\|^0(\omega): k\in \mathbb{N}~\text{and}~k\geq l \}\leq \varepsilon^0(\omega) \}$ for each $\omega\in \Omega$.\\
	Then it is easy to verify that $n_{\varepsilon}$ is $\mathcal{F}$--measurable and $\|x_k-\tilde{x}\|^0(\omega) \leq \varepsilon^0(\omega)$ for each $\omega\in \Omega$ and any $k \in \mathbb{N}$ satisfying $k\geq n^0_{\varepsilon}(\omega)$. Further, let $n_{\varepsilon}\in L^0(\mathcal{F}, \mathbb{N})$ be the equivalence class of $n^0_{\varepsilon}$. It is then clear that $\|x_n- \tilde{x}\| \leq \varepsilon$ for any $n\in L^0(\mathcal{F}, \mathbb{N})$ with $n\geq n_{\varepsilon}$, which also means that $\{x_n, n\in L^0(\mathcal{F}, \mathbb{N})\}$ is a $\mathcal{T}_c$--Cauchy stable sequence in $G$. Hence there exists $x$ in $G$ such that $\{x_n, n\in L^0(\mathcal{F}, \mathbb{N})\}$ converges in $\mathcal{T}_c$ to $x$ by the $\mathcal{T}_c$--stable sequential completeness of $G$. Finally, since for any $\varepsilon\in L^0_{++} (\mathcal{F})$ there exists $n_0\in L^0(\mathcal{F}, \mathbb{N})$ such that $\|x_n-\tilde{x}\|\leq \frac{\varepsilon}{2}$ and $\|x_n-x\|\leq \frac{\varepsilon}{2}$ for any $n\in L^0(\mathcal{F}, \mathbb{N})$ with $n\geq n_0$. Therefore, $\|x-\tilde{x}\|\leq \|x_n-x\|+ \|x_n-\tilde{x}\| \leq \varepsilon$, which yields $\tilde{x} = x\in G$ by the arbitrariness of $\varepsilon$.
\end{proof}

\par
It immediately follows form Theorem 5.10 of \cite{JZ} that a $\sigma$--stable subset $G$ of an $RN$ module $(E, \|\cdot\|)$ is stably compact iff $G$ is $\mathcal{T}_c$--stably sequentially compact iff $G$ is both random totally bounded and $G$ is $\mathcal{T}_c$--stably sequentially complete. Further, by Theorem \ref{theorem5.3} $G$ is both random totally bounded and $\mathcal{T}_c$--stably sequentially complete iff $G$ is both random totally bounded and $\mathcal{T}_{\varepsilon,\lambda}$--complete. Finally, by Theorem \ref{theorem2.3} we have that a $\sigma$--stable subset of an $RN$ module is stably compact iff it is random sequentially compact. Hence, our Theorem \ref{theorem1.12} has provided a partial answer to Problem \ref{problem5.1}.

\par
Finally, we should also point out that the authors of \cite{AZ} first studied Problem \ref{problem5.1}. They attempted to give an affirmative answer to Problem \ref{problem5.1} by making use of the transfer principle in Boolean valued analysis---Theorem 2.4.1 of \cite[p.67]{KK}, which states that every theorem of ZFC holds inside $V^{(B)}$, where $V^{(B)}$ stands for the Boolean valued universe over a complete Boolean algebra $B$. This transfer principle claims that any theorem in classical set theory has a transcription in the Boolean valued setting, which is also true with probability one. But, unfortunately the authors of \cite{AZ} went too far with the use of the transfer principle in \cite{AZ}, where they employed this principle together with some simple discussions on locally $L^0$--convex modules to reach the assertion in \cite{AZ} that ``any known results on classical locally convex spaces have a transcription in locally $L^0$--convex modules, which is also true with probability one''. Based on such a logical reasoning, they said in \cite{AZ} that they could obtain any theorem in locally $L^0$--convex modules without any proof, in particular, they obtained Theorem 2.8 of \cite{AZ}, namely, they asserted in \cite{AZ} that Problem \ref{problem5.1} had an affirmative answer for the case of a $\sigma$--stable $\mathcal{T}_c$--continuous mapping $T$. Proposition \ref{proposition5.4} below, as a special case of Theorem 3.6 of \cite{Guo3} or Proposition 4.1 of \cite{GZZ1}, has been known for over ten years, which shows that the classical separation theorem in a normed space does not hold true with probability one in an $RN$ module, namely, the reasoning used in \cite{AZ} to obtain Theorem 2.8 of \cite{AZ} does not work. Therefore, Problem \ref{problem5.1} remains open!

\begin{proposition}[\cite{Guo3,GZZ1}]\label{proposition5.4}
	Let $(E, \|\cdot\|)$ be an $RN$ module over $\mathbb{K}$ with base $(\Omega, \mathcal{F}, P)$, $G$ a $\sigma$--stable and $\mathcal{T}_c$--closed $L^0$--convex subset of $E$ and $x\in E\backslash G$. Denote $d(x,G)= \bigwedge \{\|x-y\|: y\in G\}$ (it is easy to see that $d(x, G)>0$ since $x\notin G$) and further let $A$ be an arbitrary representative of $(d(x,G)>0)$, then there exists $f\in E^*$ such that the following are satisfied:
	\begin{enumerate}[(1).]
		\item $Re(f(x)) > \bigvee \{ Re(f(y)): y\in G\}$ on $A$;
		\item $Re(f(x)) = \bigvee\{Re(f(y)): y\in G\}$ on $A^c: = \Omega \backslash A$.
	\end{enumerate}
\end{proposition}

\par
It immediately follows from Proposition \ref{proposition5.4} that it is impossible to strictly separate $x$ from $G$ with probability one when $0< P(A)< 1$.
%\paragraph{\bf Acknowledgements} The first two authors are supported by the NNSF of China (NO.11571369); the third author is supported by the NNSF of %China (NO.U1811461) and the Australian Research Council/Siscovery Project (DP200100124); the fourth author is supported by the NNSF of China %(NO.U1811462 and NO.71971031).

%\begin{thebibliography}{99}
%%% Format for Journal Reference
% etc
%\section*{Declarations}
%\section*{Acknowledgements}
%\bmhead{Acknowledgements}
%The first author was supported in part by the NNSF of China Grant \#12371141 and \# 11971483. The third author was supported in part by NNSF of China Grant \#U1811461 and the Australian Research Council/Discovery Project \#DP200100124. The fourth author was supported in part by the NNSF of China Grant \#U1811462 and \#71971031.

\bibliographystyle{amsplain}


\end{document}